\documentclass[11pt]{article}

\usepackage[utf8]{inputenc}
\usepackage{graphicx}
\usepackage{mathtools}
\usepackage{amsmath}
\RequirePackage{amssymb}
\RequirePackage{amsthm}
\usepackage{amsmath,scalefnt}
\usepackage{footnote}
\usepackage{fmtcount}
\usepackage{threeparttable}
\usepackage{caption}
\usepackage{subcaption}
\usepackage{color}
\usepackage{tikz}
\usetikzlibrary{positioning,chains,fit,shapes,calc}
\usetikzlibrary{arrows}
\tikzstyle{vertex}=[circle, draw, inner sep=0pt, minimum size=6pt]
 \usepackage{enumitem}
 \usepackage{multicol}
 \usepackage{wrapfig}
\usepackage{enumitem}
\usepackage{bm} 

	\newtheorem{definition}{Definition}
	\newtheorem{claim}{Claim}
	\newtheorem{lemma}{Lemma}

	 \newtheorem{assumption}{Assumption}
	 \newtheorem{proof of proposition}{Proof of Proposition}
	\newtheorem{proof of claim}{Proof of Claim}
	\newtheorem{proposition}{Proposition} 
	\newtheorem{theorem}{Theorem}
	\newtheorem{notation}{Notation}
        \newtheorem{corollary}{Corollary}
        \newtheorem{remark}{Remark}

\newcommand{\diag}{\mbox{diag}\,}
\newcommand{\sech}{\mbox{sech}\,}

\renewcommand{\r}{{\mathbb R}}
\renewcommand{\a}{{\alpha}}

\renewcommand{\o}{{\omega}}

\newcommand{\w}{{w}}   
\renewcommand{\c}{{c}}  
\newcommand{\x}{{\xi}} 
\newcommand{\et}{{\eta}} 
\newcommand{\Ii} {{I_{ext}^i(t)}}
\newcommand{\I} {{\mathbf{I(t)}}}
\newcommand{\dd}{{\delta}}

\newcommand{\ee}{\end{equation}}
\newcommand{\bal}{\begin{aligned}}
\newcommand{\eal}{\end{aligned}}
\newcommand{\bi}{\begin{itemize}}
\newcommand{\ei}{\end{itemize}}
\newcommand{\ben}{\begin{enumerate}}
\newcommand{\een}{\end{enumerate}}
\newcommand{\beqn}{\begin{eqnarray*}}
\newcommand{\eeqn}{\end{eqnarray*}}

\newcommand{\be}[1]{\begin{equation}\label{#1}}
\newcommand{\bp}{\begin{proof}}
\newcommand{\ep}{\end{proof}}
\newcommand{\bremark}{\begin{remark}\rm } 
\newcommand{\eremark}{\end{remark}}
\newcommand{\blem}{\begin{lemma}}
\newcommand{\elem}{\end{lemma}}
\newcommand{\bclaim}{\begin{claim}}
\newcommand{\eclaim}{\end{claim}}
\newcommand{\bnote}{\begin{notation}}
\newcommand{\enote}{\end{notation}}

\newcommand{\bthm}{\begin{theorem}}
\newcommand{\ethm}{\end{theorem}}
\newcommand{\bprop}{\begin{proposition}}
\newcommand{\eprop}{\end{proposition}}
\newcommand{\bcor}{\begin{corollary}}
\newcommand{\ecor}{\end{corollary}}
\newcommand{\dis}{\displaystyle}
\newcommand{\lt}{\left}
\newcommand{\rt}{\right}
\newcommand{\FR}{A_{FR} }
\newcommand{\FL}{A_{FL} }
\newcommand{\BR}{A_{BR} }
\newcommand{\BL}{A_{BL} }
\newcommand{\T}{A_{Tri} }

\graphicspath{{./Hetero_Figs/}}

\title{Heterogeneous inputs to central pattern generators\\ can shape insect gaits.}
\author{Zahra Aminzare%
  \thanks{The Program in Applied and Computational Mathematics,  Princeton University, NJ, USA.
  Email: aminzare@math.princeton.edu (Zahra Aminzare).}
 \and Philip Holmes%
  \thanks{The Program in Applied and Computational Mathematics,  
  Department of Mechanical and Aerospace Engineering,  
  and Princeton Neuroscience Institute,  Princeton University, NJ, USA.
 Email: pholmes@math.princeton.edu (Philip Holmes).}
 }
\date{}                                           

\begin{document}
\maketitle

\begin{abstract}
In our previous work \cite{SIADS2018}, we studied an interconnected bursting neuron model for 
insect locomotion, and its corresponding phase oscillator model, which at high speed can generate 
stable tripod gaits with three legs off the ground simultaneously in swing, and at low speed can 
generate stable tetrapod gaits with two legs off the ground simultaneously in swing. 
However, at low speed several other stable locomotion patterns, that are not typically observed 
as insect gaits, may coexist. In the present paper, by adding heterogeneous external input to each oscillator, 
we modify the bursting neuron model so that  its corresponding phase oscillator model 
produces only one stable gait at each speed, specifically: a unique stable tetrapod gait 
at low speed, a unique stable tripod gait at high speed, and a unique branch
of stable transition gaits connecting them. This suggests that control signals originating 
in the brain and central nervous system can modify gait patterns.
\end{abstract}



\section{Introduction}
\label{s.introduction}

This paper is based on our previous work \cite{SIADS2018}, in which we studied the 
effect of stepping frequency on transitions from multiple tetrapod insect gaits with two legs 
off the ground simultaneously in swing, to tripod gaits with three legs off the ground 
simultaneously in swing. In that paper, we used an ion-channel bursting 
neuron model  to describe each of six mutually inhibitory units that form the central 
pattern generator (CPG) located in the insect's thorax. We assumed that each unit 
drives one leg of the insect and that the units are connected to 
their nearest neighbors in a homogeneous network as shown in Figure~\ref{6-oscillator}
below, but where input currents are identical, $I_i = I_{ext}$.
 
Employing phase reduction, we collapsed the network of bursting neurons represented by 
24 ordinary differential equations to 6 coupled nonlinear phase oscillators, each corresponding 
to a sub-network of neurons controlling one leg.
Assuming that the left  and right legs maintain constant phase differences (contralateral 
symmetry), we reduced from 6 equations to 3, allowing analysis of a  dynamical system 
with 2 phase differences defined on a 2-dimensional torus. 

With certain balance conditions on the coupling strengths 
among the homogeneous oscillators, described in Section~\ref{Reduced_phase_section} below,
we showed that at low speeds, the phase differences model on the torus 
can generate multiple fixed points, including stable tetrapod and unstable tripod gaits. 
In contrast, at high speeds, it generates a unique stable tripod gait. 
Moreover, as speed increases, the gait transition occurs  through 
degenerate bifurcations, at which a subset of the multiple fixed 
points merge to produce a unique stable fixed point \cite[Figure 23]{SIADS2018}.

In the current paper, we study this degenerate bifurcation in the phase difference model, 
by unfolding the original system. 
To this end we relax the condition of homogeneous (identical) ion-channel 
bursting neuron models in the network of CPGs  and allow heterogeneous (nonidentical) 
 models by adding different external inputs to each oscillator. 
We subsequently show that this heterogeneity is equivalent to perturbing the coupling functions 
or the contralateral coupling strengths in a phase reduced oscillator model; i.e. different
heterogeneities can have the same effects on dynamics.
 (see Section \ref{Equivalent_Perturbations} below).
 
The paper is organized as follows. 
In Section~\ref{Bursting neuron model}, we review the ion-channel model for bursting neurons 
which was developed in \cite{SIAM2,SIAM1}, and the influence of 
its parameters on speed, which was studied in \cite{SIADS2018}. 
In Section~\ref{phase_reduction}, we review the derivation of  phase equations for  
 heterogeneous networks and apply these techniques to the interconnected bursting neuron model. 
In Section~\ref{Reduced_phase_section}, we define approximate tetrapod,
tripod and transition gaits for  
heterogeneous networks and then, by assuming constant phase differences 
between left- and right-hand oscillators, as in the homogeneous case,
we reduce the 6 phase equations to 2 phase difference equations.
In Section~\ref{Main_Result}, we describe the main results of this paper. 
By choosing appropriate heterogeneous external inputs, we show that the phase 
differences model generates only one stable fixed point, which at low speed,
corresponds to a tetrapod gait and at high speed, corresponds to a tripod gait. 
Interpreting the heterogeneities as small bifurcation parameters, we find cases
in which two or three saddle-node 
bifurcations occur as heterogeneity increases and a unique tetrapod gait 
emerges from multiple tetrapod gaits and other, ill-defined gaits. This shows that
specific fixed points (gaits) can be preserved, or removed, by small external input currents.
In Section~\ref{Equivalent_Perturbations}, 
we show that our heterogeneities are equivalent to perturbing 
the coupling functions or the contralateral coupling strengths in a phase reduced
oscillator model. In Section~\ref{Conclusion}, we conclude.  

\section{A network of weakly interconnected bursting neurons}
\label{Bursting neuron model}

 In \cite[Section 2.1]{SIADS2018}, we employed an ion-channel bursting neuron model 
 for an insect central pattern generator which was developed
in \cite{SIAM2,SIAM1}. The bursting neuron model of each unit of the CPG contains a system of 4 ODEs 
describing trans-membrane cell voltages, slow and fast ionic gates, and the dynamics of 
neurotransmitter release at synapses, as follows:
\begin{subequations}\label{BN}
\begin{align}
      C\dot {v}&= -\left\{I_{Ca}(v)+I_K(v,m)+I_{KS}(v,\w)+I_L(v)\right\}  + I_{ext},\label{BN1}\\
         \dot {m}&= \dis\frac{\gamma}{\tau_m(v)} [m_{\infty} (v)-m],\label{BN2}\\
         \dot {\w}&= \dis\frac{\delta}{\tau_{\w}(v)} [\w_{\infty} (v)-\w],\label{BN3}\\
         \dot {s}&= \dis\frac{1}{\tau_{s}}[s_{\infty} (v)(1-s)-s], 
\label{BN4}
\end{align}
\end{subequations}
where the ionic currents are of the following forms
\be{currents}
\bal
      &I_{Ca} (v)\;=\; \bar g_{Ca} n_{\infty}(v) (v-E_{Ca}),\qquad
        I_{K} (v,m)\;=\; \bar g_{K} \;m \;(v-E_{K}),\\
      &I_{KS} (v,\w)\;=\; \bar g_{KS} \w \;(v-E_{KS}), \qquad
        I_L(v)\;\;=\;\;\bar g_L(v-E_L).
\eal
\ee
The steady state gating variables associated with ion channels and their
time scales take the forms
\be{steady_states}
\bal 
        &m_{\infty}(v) \;=\; \frac{1}{1+e^{-2k_{K}(v-v_K)}}\;,\qquad \qquad
          \w_{\infty}(v) \;=\; \frac{1}{1+e^{-2k_{KS}(v-v_{KS})}}\;,\\
        &n_{\infty} (v)\;=\;  \frac{1}{1+e^{-2k_{Ca}(v-v_{Ca})}}\;,\qquad\qquad
          s_{\infty} (v)\;=\;  \frac{a}{1+e^{-2k_s(v-E^{pre}_{s})}}\;,
\eal
\ee 
 and 
\be{time_scales}
\bal 
       \tau_{m} (v)= \sech( k_{K}(v-v_{K})),\quad
       \tau_{\w} (v)= \sech( k_{KS}(v-v_{KS})).
\eal
\ee

The external current $I_{ext}$, which represents input from the central nervous system and brain, 
 varies between $35.65$ and $37.7$ as speed increases. 
 Other parameters are generally fixed as specified in 
Table \ref{table_of_parameters} and are chosen such that the model (\ref{BN}) possesses an
 attracting hyperbolic limit cycle $\Gamma$.
 Most of the parameter values are taken 
from \cite{SIAM1}, but some of our notations are different.  
See \cite[Section 2.1]{SIADS2018}
for further details of the model and its parameters. 

\begin{table}[h!]
\begin{center}
 \begin{tabular}{|c c c c c c  c c c c c  c c |} 
 \hline
  & $I_{ext}$  
  & $\bar g_{Ca}$ 
  & $\bar g_K$ 
  & $\bar g_{KS}$ 
  & $\bar g_{L}$ 
  & $\bar g_{syn}$ 
  & $E_{Ca}$ 
  & $ E_K$ 
  &$ E_{KS}$ 
  &  $ E_{L}$ 
  & $E^{post}_s$  
  &  $E^{pre}_s$ 
 \\ 
  [0.5ex] 
 \hline\hline
 & varies 
 &4.4 
 &9.0 
 &0.5 
 &2.0
 &0.01
 &120
 &-80 
 &-80
 &-60
 &-70
  & 2
 \\
  [1ex] 
 \hline
 \end{tabular}
 
 \begin{tabular}{|c c c c c  c c c c c c c c  |} 
 \hline
   & $k_{Ca}$ 
   & $k_K$
   & $k_{KS}$
   & $k_s$ 
   & $v_{Ca}$ 
   & $v_K$
   & $v_{KS}$
   & a 
   & $C$ 
   & $\gamma$ 
   & $\tau_s$
     & $\delta$  
   \\ 
  [0.5ex] 
 \hline\hline
 & 0.056
 & 0.1 
 &0.8 
 & 0.11
 &-1.2
 &2
 &-26
 &444.48
 &1.2
 &5.0
 &5.56
  & 0.027 
  \\
  [1ex] 
 \hline
 \end{tabular}
\caption{{The constant parameters in the bursting neuron model.}}
\label{table_of_parameters}
\end{center}
\end{table}

\begin{figure}
\centering 
         \includegraphics[width=.4\textwidth]{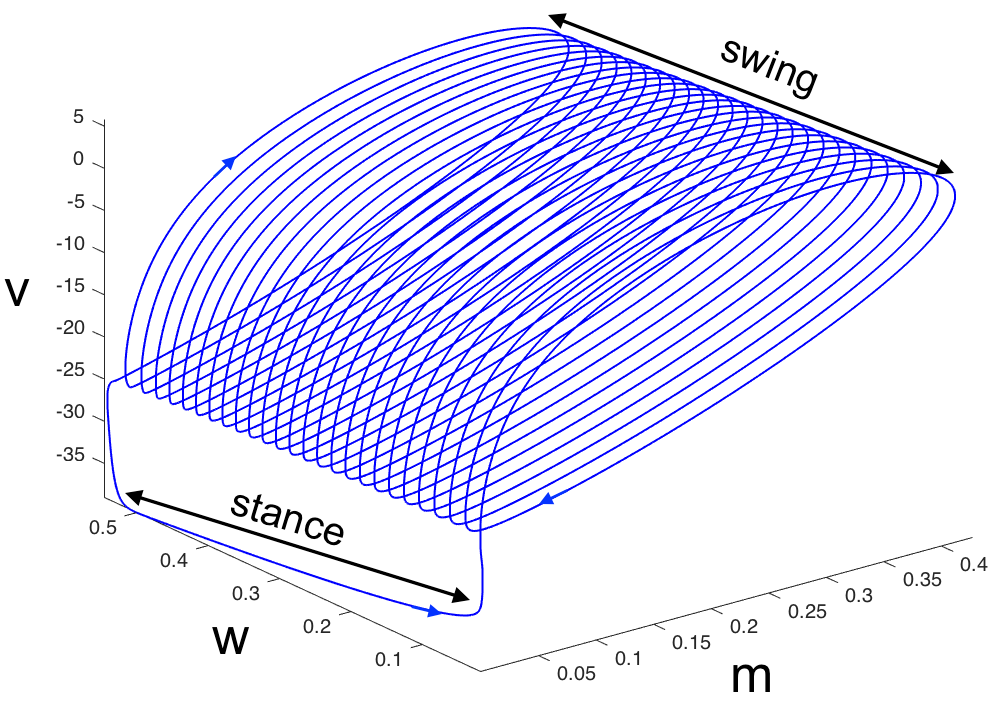}\qquad\qquad
           \includegraphics[width=.3\textwidth]{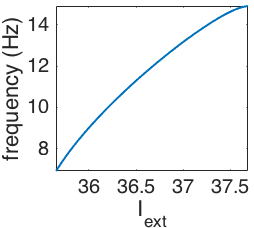}
\caption{ Left: A periodic orbit of the bursting neuron model, Equation~(\ref{BN}),  in $(v,m,\w)$ space. 
               Right: The effect of $I_{ext}$ on the stepping frequency of the periodic orbit.  }
\label{bursting_vmw}
\end{figure}

As shown in Figure~\ref{bursting_vmw}(left), the periodic orbit in $(v,m,\w)$ space contains a 
sequence of spikes (a burst) followed by a quiescent phase, which correspond respectively to 
the swing and stance durations of one leg. The burst from the CPG inhibits depressor 
 motoneurons and excites levator motoneurons, 
allowing the swing leg to lift from the ground \cite[Figure 2]{SIAM2} and 
\cite[Figure 11]{pearson_Iles_1973} (see also  \cite{pearson_Iles_1970,pearson_1972}). 
We denote the period of the periodic orbit by $T$, i.e., it takes $T$
time units (ms here) to complete the stance and swing cycle of each leg.  
The number of  steps completed by one leg per unit of time is the
stepping frequency and is equal to $\omega = 2\pi/T$. 
In \cite[Figure 2]{SIADS2018}, we observed that as one of the two parameters in the bursting
neuron model, either the slow time scale $\delta$ or the external current $I_{ext}$, increases, 
the period of the periodic orbit decreases, 
primarily by decreasing stance duration, and so the insect's speed increases. 
There, we used these parameters as speed parameters, denoted by $\x$, and studied transitions from 
tetrapod to tripod gaits as $\x$ increases. Gait transition is not the main focus of 
 the present paper, although here we will  show gait transitions
using only $I_{ext}$ as a speed parameter. To see how $I_{ext}$ affects 
the frequency of the periodic orbit, see Figure~\ref{bursting_vmw}(right). 
 
\begin{figure}
\centering 
         \includegraphics[width=.25\textwidth]{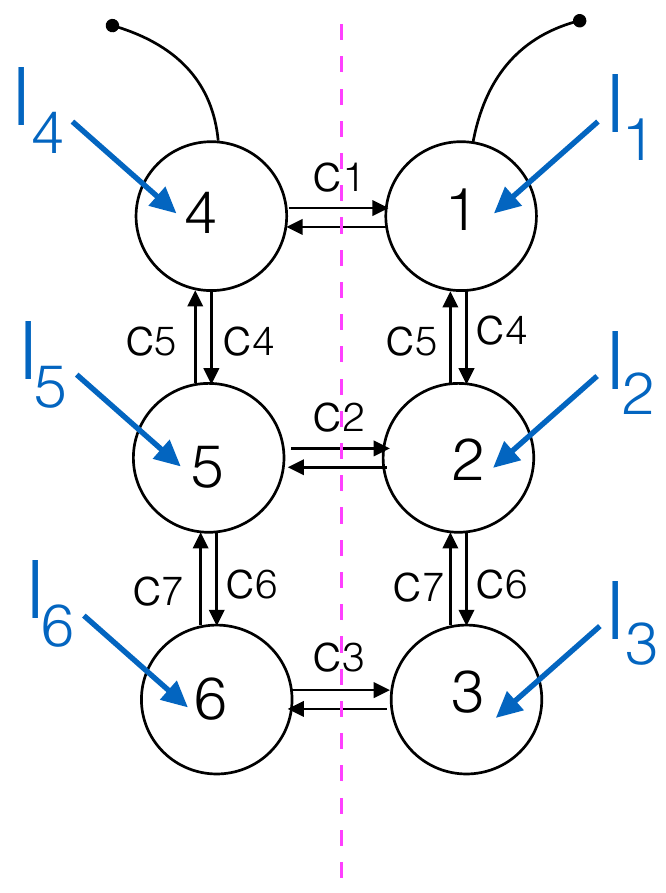}
\caption{The network of heterogeneous CPGs with different external inputs, 
$I_i = I_{ext} +  I^i_{ext}(t)$.}
\label{6-oscillator}
\end{figure}

In \cite{SIADS2018}, we  assumed that inhibitory coupling is achieved via synapses 
that produce negative postsynaptic currents. The synapse variable $s$ enters the 
postsynaptic cell in Equation~(\ref{BN1}) as an additional term, $I_{syn}$, 
\be{synapse_v}
\bal 
          C\dot {v}_i&= -\left\{I_{Ca}+I_K+I_{KS}+I_L\right\}  + I_{ext} + I_{syn}\;,
\eal
\ee 
where  
\be{I_syn}
\bal 
          I_{syn}=\dis\sum_{j\in \mathcal{N}_i} I_{syn}(v_i,s_j) 
          = \dis\sum_{j\in \mathcal{N}_i} - \bar\c_{ji} \bar g_{syn} s_j \lt(v_i- E_s^{post}\rt), 
 \eal
\ee  
$\bar g_{syn}$ denotes the synaptic strength, and $\mathcal{N}_i$
 denotes the set of the nodes adjacent to node $i$.   
The multiplicative factor $\bar\c_{ji}$ accounts for the fact that multiple
bursting neurons are interconnected in the insects,  and 
$- \bar\c_{ji} \bar g_{syn}$ represents an overall coupling strength between 
hemi-segments. Following \cite{Fuchs14} we assumed contralateral symmetry
and included only nearest neighbor coupling, so that there are three
contralateral coupling strengths $\c_1, \c_2, \c_3$ and four ipsilateral
coupling strengths $\c_4, \c_5, \c_6,$ and $\c_7$; see Figure \ref{6-oscillator}. 
For example, $\bar\c_{21} = \c_5$, $\bar\c_{41} = \c_1,$ etc.  
We chose reversal potentials $E^{post}_s$ that make all synaptic connections inhibitory;
 this implies that the $c_i$'s are positive.
 
The following system of $24$ ordinary differential equations (ODEs)
describes the dynamics of the $6$ coupled
cells in the network as shown in Figure~\ref{6-oscillator}. We assume that
each cell, which is governed by Equation (\ref{BN}), represents one leg of
the insect. Cells $1,2$, and $3$ represent right front, middle, and hind legs,
and cells $4, 5$, and $6$ represent left front, middle, and hind legs, respectively:
 \begin{equation} \label{24ODE_closed}
\begin{array}{ll}
& \dot{x}_1  \;=\; f(x_1) + \c_1g(x_1,x_4) + \c_5g(x_1,x_2), \\
& \dot{x}_2  \;=\; f(x_2) + \c_2g(x_2,x_5) + \c_4g(x_2,x_1) + \c_7g(x_2,x_3), \\
& \dot{x}_3  \;=\; f(x_3) + \c_3g(x_3,x_6) + \c_6g(x_3,x_2), \\
& \dot{x}_4  \;=\; f(x_4) + \c_1g(x_4,x_1) + \c_5g(x_4,x_5), \\
& \dot{x}_5  \;=\; f(x_5) + \c_2g(x_5,x_2) + \c_4g(x_5,x_4) + \c_7g(x_5,x_6), \\
& \dot{x}_6  \;=\; f(x_6) + \c_3g(x_6,x_3) + \c_6g(x_6,x_5), 
\end{array}
\end{equation}
where $x_i= (v_i,m_i,\w_i,s_i)^\top$,  $f(x_i)$ is as the right hand side of Equations~(\ref{BN})
and 
\be{g:interconnected:bursting:neuron}
g (x_i,x_j)= \lt(-  \bar g_{syn} s_j \lt(v_i-E_s^{post}\rt) , 0, 0, 0 \rt)^\top, 
\ee
is the coupling function with a small synaptic coupling strength $\bar g_{syn}$. 
This assumption of weak coupling is necessary for the use of phase reduction in Section~\ref{phase_reduction}. 

This 6-bursting neuron model was used to drive agonist-antagonist muscle pairs in 
a neuro-mechanical model with jointed legs that reproduced the dynamics of 
freely-running cockroaches \cite{KukProc09}, also see \cite{kukillaya.09}. 
These papers and subsequent phase-reduced models \cite{ProctorKH10, Fuchs14} 
support our belief that the bursting neuron model is capable of producing realistic 
inputs to muscles in insects. 
In \cite[Figures 5 and 6]{SIADS2018}, we showed that the 24 ODEs
coupled bursting neuron model with 
small $I_{ext}$ (or $\delta$) can produce a tetrapod gait with two legs lifted off the 
ground simultaneously in swing; and as $I_{ext}$ (or $\delta$) increases, it can produce a 
tripod gait with three legs lifted off the ground simultaneously in swing. 

In \cite[Section 2.2]{SIADS2018} we considered a network of six identical mutually inhibiting 
homogeneous units, representing the hemi-segmental CPG networks contained in the 
insect's thorax. 
In the present paper,  
we assume that in addition to $I_{ext}$, each unit receives a different external 
input denoted by $\Ii$, as shown in Figure~\ref{6-oscillator}, where $\Ii$ is a periodic 
function with frequency $\omega_i$ close to $\omega$ and a magnitude of order 
$\mathcal{O}(\bar g_{syn})$, and $\bar  g_{syn}$ is the synaptic strength. 
Therefore, in all the previous equations, $I_{ext}$ will be replaced by $I_{ext}+\Ii$, 
where $i$ denotes the leg number. 
Later we also assume that  for $i=1,2,3$, $\Ii=I_{ext}^{i+3}$, 
to preserve the contralateral symmetry  condition;
see Assumption \ref{external_input_assumption} in Section~\ref{Reduced_phase_section}. 

To analyze the gait transition mathematically, in  \cite[Section 3]{SIADS2018}, 
we applied the theory of weakly coupled oscillators to the coupled bursting neuron 
models to reduce the 24 ODEs to 6 phase oscillator equations.  
In the following section we apply the phase reduction technique again and derive 6  
phase oscillator equations in the presence of heterogeneous external inputs. 
 
\section{A phase oscillator model}
\label{phase_reduction}

Let the ODE
\be{single:neuron}
\dot x = f(x), \qquad  x \in \r^n, 
\ee
describe the dynamics of a single oscillator.  Assume that
Equation~(\ref{single:neuron}) has an attracting hyperbolic limit cycle
{$\Gamma = \Gamma(t)$}, with period $T$ and frequency $\o=2\pi/T$. 
The phase of an oscillator, denoted by $\phi$, is the time that has elapsed as its state
moves around $\Gamma$, starting from an arbitrary reference point
in the cycle, called relative phase.  
In this section, we derive the phase equations of weakly coupled oscillators 
with heterogeneous dynamics, i.e., coupled oscillators with different frequencies.  
 We develop the theory in greater generality than our specific applications
will demand, allowing periodic input currents $I^i_{ext}(t)$. 

\subsection {A pair of weakly coupled heterogeneous oscillators}\label{pair_coupled}

Consider a pair of weakly coupled  heterogeneous oscillators 
\be{coupled:neurons}
\bal
&\dot x_1= f_1(x_1) + \epsilon g(x_1,x_2),\\
&\dot x_2= f_2(x_2) + \epsilon g(x_2,x_1),
\eal
\ee
where $f_i$ describes the intrinsic dynamics of each oscillator,  
$0<\epsilon\ll1$ is the coupling strength, and $g$ is the coupling function.
 For each oscillator, the phase equation can be written as follows. 
 For more details see \cite[Section 3]{SIADS2018}. 
 \be{take:average:33}
\dis\frac{d\phi_i}{dt}(t) =\omega_i+ \epsilon H_i(\phi_j(t) -\phi_i(t))+\mathcal{O}(\epsilon^2),
\ee
where 
\[
H_i = H_i(\theta) =  \frac{1}{T}\int_0^T Z_i (\Gamma_i(\tilde t)) \cdot 
g(\Gamma_i(\tilde t),\Gamma_i(\tilde t +\theta )) \; d\tilde t,
\]
is the coupling function: the convolution  of the coupling $g$ 
and the oscillator's infinitesimal phase response curve (iPRC), $Z_i$,
and $\omega_i$ is the frequency 
of each oscillator described by $\dot x_i= f_i(x_i)$. Under the weak coupling assumption,
the iPRC captures the local dynamics  of each oscillator in a neighborhood of 
its limit cycle $\Gamma_i$,  \cite{Schwemmer2012}. 

Equation~(\ref{take:average:33}) is a general phase equation for a pair of weakly coupled heterogeneous oscillators where the
heterogeneity is of any arbitrary size. This means that the oscillators' frequencies can be very different
from each other. But if the frequencies are close to each other, i.e, the heterogeneities are 
small and in particular are of order of the coupling strength $\epsilon$, 
then one can approximate Equation~(\ref{take:average:33}) as follows \cite[Chapter 5]{Kuramoto_book}. 

Assume that $f_i = f +  \tilde f_i$, where the heterogeneity $\tilde f_i$ is periodic with period close to the period of $f$ 
and $ \tilde f_i$ is of order $\epsilon, \mathcal{O}(\epsilon)$. 
This is equivalent to having identical oscillators with dynamics $f$ and 
non-identical coupling functions $g_i =  g + \tilde f_i /\epsilon$.
Then Equation~(\ref{take:average:33}) can be approximated by the following phase equations:
\be{take:average:3}
\dis\frac{d\phi_i}{dt}(t) =\omega + \tilde\omega_i+ \epsilon H(\phi_j(t) -\phi_i(t))+\mathcal{O}(\epsilon^2),
\ee
where 
\[
H = H(\theta) =  \frac{1}{T}\int_0^T Z (\Gamma(\tilde t)) \cdot g(\Gamma(\tilde t),\Gamma(\tilde t +\theta )) \; d\tilde t,
\]
is a coupling function: specifically, the convolution of the synaptic coupling $g$
and the oscillator's infinitesimal phase response curve (iPRC), $Z$. 
Here $Z$ is computed for the limit cycle of  $\dot x= f(x)$, and the frequency
differences are constant shifts of $\mathcal{O}(\epsilon)$:
\be{heterogenous_w}
\tilde\omega_i = \frac{1}{T}\int_0^T Z (\Gamma(\tilde t)) \cdot \tilde f_i(\Gamma(\tilde t)) \; d\tilde t.
\ee  
The advantage of this
decomposition is that only one iPRC and so only one coupling function need to be computed.
 
 In what follows we generalize the approximation of Equation~(\ref{take:average:3}) to a 
 network of weakly coupled heterogeneous neurons with multiple connections \cite[Chapter 5]{Kuramoto_book}. 
 
 \subsection {A network of weakly coupled heterogeneous oscillators}\label{network_coupled}
 
Now consider a network of $N$ heterogeneous oscillators with intrinsic dynamics 
$\dot x_i= f_i(x_i)$ and corresponding frequencies $\omega_i$. For $i=1,\ldots, N$, let
\be{network:coupled:neurons}
\bal
& \dot x_i= f_i(x_i) + \sum_{j\in\mathcal{N}_i} \epsilon_j g(x_i,x_j),\\
\eal
\ee
describe the dynamics of each $x_i$ in the weakly coupled network. Here
$\mathcal{N}_i$ denotes the neighbors of oscillator $i$, and $\epsilon_j$ 
denotes the coupling strengths, which are all of $\mathcal{O}(\epsilon)$
for some $0<\epsilon\ll1$.  
As in the case of a pair of coupled oscillators, one can derive phase equations from  
Equation~(\ref{network:coupled:neurons}) as follows 
\be{network:phase}
\bal
&\dot \phi_i= \omega_i +\sum_{j\in\mathcal{N}_i} \epsilon_j H_j(\phi_j- \phi_i)+\mathcal{O}(\epsilon^2),
\eal
\ee
where the coupling function $H_j$ is the convolution of 
the coupling function $g$ and each oscillator's iPRC, $Z_i$.  
 
Now, similar to Equation~(\ref{take:average:3}), we approximate
Equation~(\ref{network:phase}) such that all the coupling functions can be computed 
from a single iPRC which corresponds to the limit cycle of $\dot x= f(x)$. 
Assume that $f_i  = f +  \tilde f_i$, where $\tilde f_i$ is periodic with period close
to the period of $f$, and $\tilde f_i$ is of $\mathcal{O}(\epsilon)$.
As in the case of a pair of coupled neurons, since the perturbation $\tilde f_i=\mathcal{O}(\epsilon)$, 
we can approximate each limit cycle by the limit cycle of  $\dot x= f(x)$ and
consider $\tilde f_i$ as a perturbation to the coupling function $g$. 
Therefore, Equation~(\ref{network:phase}) can be written as 
\be{network:phase:approx}
\bal
&\dot \phi_i= \omega + \tilde \omega_i + \sum_{j\in\mathcal{N}_i} \epsilon_j H(\phi_j- \phi_i)+\mathcal{O}(\epsilon^2),
\eal
\ee
 where $\omega$ is the frequency of the limit cycle of  $\dot x=f (x)$, 
 $H$  is the convolution of the coupling function $g$ and $Z$, 
  the iPRC of the limit cycle of $\dot x=f (x)$, and the frequency differences are
 \be{heterogenous_w_network}
 \tilde\omega_i \;=\;  \frac{1}{T}\int_0^T Z (\Gamma(\tilde t)) \cdot \tilde f_i(\Gamma(\tilde t)) \; d\tilde t. 
  \ee
 
\subsection {Phase equations for six weakly coupled heterogeneous bursting neuron model}

We now apply the techniques from Section \ref{network_coupled} to
six heterogeneous units in the coupled bursting neuron model 
 and derive the 6-coupled phase oscillator model via phase reduction.

In the interconnected bursting neuron model, Equation~(\ref{24ODE_closed}), 
the intrinsic dynamics of each hemi-segmental unit, described by $\dot x= f(x)$, is homogeneous. 
We now assume that each hemi-segmental unit receives a small heterogeneous external input, i.e., 
each unit can be described by 
\[\dot x_i= f_i(x_i) = f(x_i) + \Ii\cdot\mathbf{e}_1,\] 
where $\dot x_i= f_i(x_i)$ has an  attracting hyperbolic limit cycle with frequency close to the 
 attracting hyperbolic limit cycle of $\dot x= f(x)$. 
 For $i=1,\ldots, 6$, $\Ii$ is the additional small external input to each unit and represents the 
weak heterogeneity of the corresponding  unit such that  $\Ii=\mathcal{O}(\bar g_{syn})$ and 
 $\mathbf{e}_1 = (1, 0 , 0, 0)^\top$, i.e., only the voltage equations are heterogeneous. 
  
\begin{figure}[h!]
\begin{center}
 \includegraphics[scale=.4]{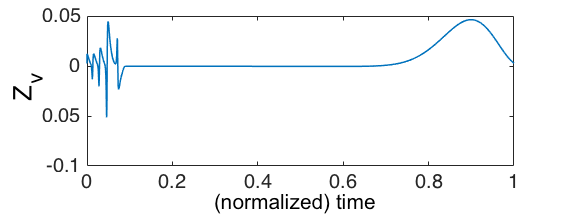}
  \includegraphics[scale=.4]{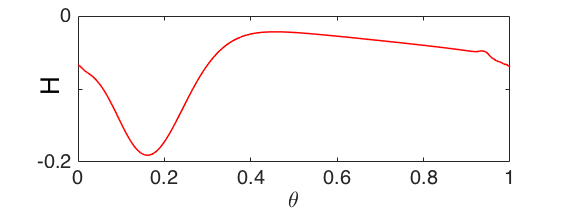}
\caption{\cite[Figure 8]{SIADS2018}  iPRC (in the direction of $v$) (left) 
and  the coupling function $H(\theta)$ (right) for $I_{ext}=35.9$. 
Phase $\theta=0$ is defined to be the onset of the burst.}
\label{H_PRC_versus_Iext}
\end{center}
\end{figure}

Recalling Equation~(\ref{network:phase:approx}), we can derive approximate phase equations for 
 the coupled bursting neuron model of Figure~\ref{6-oscillator} as follows. 
\begin{equation} \label{eq.osc1}
\begin{array}{ll}
& \dot{\phi}_1  \;=\; \omega + \tilde\omega_1+ \c_1H(\phi_4 - \phi_1) + \c_5H(\phi_2 - \phi_1), \\
& \dot{\phi}_2  \;=\; \omega + \tilde\omega_2+ \c_2H(\phi_5 - \phi_2) + \c_4H(\phi_1 - \phi_2) + \c_7H(\phi_3 - \phi_2), \\
& \dot{\phi}_3  \;=\; \omega + \tilde\omega_3+ \c_3H(\phi_6 - \phi_3) + \c_6H(\phi_2 - \phi_3), \\
& \dot{\phi}_4  \;=\; \omega + \tilde\omega_4+ \c_1H(\phi_1 - \phi_4) + \c_5H(\phi_5 - \phi_4), \\
& \dot{\phi}_5  \;=\; \omega + \tilde\omega_5+  \c_2H(\phi_2 - \phi_5) + \c_4H(\phi_4 - \phi_5) + \c_7H(\phi_6 - \phi_5), \\
& \dot{\phi}_6  \;=\; \omega + \tilde\omega_6+ \c_3H(\phi_3 - \phi_6) + \c_6H(\phi_5 - \phi_6), 
\end{array}
\end{equation}
where
\be{tilde_omega_i}
\tilde\omega_i  = \frac{1}{T}\int_0^T Z_v (\Gamma(t)) I_{ext}^i(t) dt,  
\ee
and  $Z_v$ is the iPRC of the limit cycle of $\dot x = f(x)$ in the direction of voltage.
In Figure~\ref{H_PRC_versus_Iext}
(left), we show  $Z_v$  for $I_{ext}=35.9$. 
Note that  the averaging theorem and convolution integral
in phase reduction eliminates time dependence in $I^i_{ext}$.
Also, the coupling function $H$ takes the following form:
\be{H:interconnected:bursting:neuron}
H(\theta )=   -\frac{\bar g_{syn}}{T}\int_0^T Z_v (\Gamma(t))  \lt(v_i(\Gamma(t))  -E_s^{post}\rt) s_j\lt(\Gamma(t +\theta )\rt)  \; dt. 
\ee

In Figure~\ref{H_PRC_versus_Iext}
(right), we show the coupling function $H$ derived in 
Equation~(\ref{H:interconnected:bursting:neuron}) for $I_{ext}=35.9$. 
Note that $H(\theta)<0$
over most of its range, and in particular over the interval $[1/3,2/3]$, 
which we will show contains the tetrapod, tripod and transition gaits. 
 
To simplify the notations, for the remainder of the paper,  $T=1$ and 
all the phases and the coupling functions are considered in the domain
of $[0, 1]$ instead of $[0, 2\pi]$.

 \section{Reduced phase equations }
\label{Reduced_phase_section}

In this section, the goal is to reduce 
the 6 equations (\ref{eq.osc1}) to 2 equations on a 2-torus. 
Although we are interested in gaits generated by the bursting neuron model and 
its phase reduction equations (\ref{eq.osc1}), we prove our
results for a more general case.
To this end, we assume the following conditions 
for the coupling function $H$ and the external inputs $\Ii$. 
We let the coupling function $H$ and the frequency
$\omega$ depend on the speed parameter $\x$ and write $H=H(\phi;\x)$ and $\omega=\omega(\x)$. 

\begin{assumption}
\label{eta_assumption}
Let $H = H(\theta; \x)$ be a differentiable function, defined on
$\r\times[\x_1,\x_2]$ which is $1$-periodic on its first argument and has 
the following property. For any fixed $\x \in [\x_1,\x_2]$, 
\be{eta}
H\lt(\frac{2}{3} - \et; \x\rt) = H\lt(\frac{1}{3}+\et; \x\rt),
\ee
has a unique solution $\eta(\x)$ such that $\eta = \eta(\x):[\x_1,\x_2] \to [0,1/6]$
is an onto and non-decreasing function. Note that Equation~(\ref{eta}) is also
trivially satisfied by the constant solution $\et= 1/6$. 
 \end{assumption}
Assumption \ref{eta_assumption} defines a class of coupling functions that exhibit
the gait transitions studied here and in \cite{SIADS2018}. 
The coupling functions derived from the 
bursting neuron model satisfy  and motivate this assumption, see 
 \cite[Figure 9]{SIADS2018}. 
For the rest of the paper, we assume that the coupling function $H$ satisfies 
Assumption \ref{eta_assumption}. 
In \cite[Proposition 11]{SIADS2018}, we provided sufficient conditions for
Assumption \ref{eta_assumption} to hold for more general classes of functions. 

\begin{assumption}
\label{external_input_assumption}
For $i=1,2,3$, let $I_{ext}^{i+3} (t) = \Ii$. 
This assumption maintains contralateral symmetry.  
In addition, we assume that  for $i=1,2,3$, $\Ii$ are not equal, 
otherwise the system becomes homogeneous.  
 \end{assumption}

\begin{assumption}\label{assu:balance}
Let the coupling strengths satisfy the following balance condition 
\be{balance}
\c_1 + \c_5 = \c_2 + \c_4 + \c_7 = \c_3 + \c_6. 
\ee
 \end{assumption}
 
In \cite[Proposition 3 and Corollary 4]{SIADS2018} we proved that when the
coupling strengths satisfy the balance condition, in the homogeneous case $\tilde\omega_i=0$, 
Equation~(\ref{eq.osc1}) admits 
tetrapod gaits at low speeds and tripod gaits at high speeds. 

\subsection{Gait definitions: Generalization to heterogeneous systems}

In \cite[Definition 1]{SIADS2018} we defined four versions of tetrapod gaits and a tripod gait.   
Each gait corresponded
to a $1$-periodic solution of Equation~(\ref{eq.osc1}) with $\tilde \omega_i=0$. 
In what follows, we generalize  those definitions to heterogeneous models, i.e., 
Equation~(\ref{eq.osc1}) with at least one $\tilde \omega_i\neq0$.   

\begin{definition}[Approximate tetrapod and tripod gaits] \label{defn1}

 The approximate gaits, denoted by $A_T^p$, 
 are $1$-periodic solutions of Equation~(\ref{eq.osc1}) with at least one $\tilde \omega_i\neq0$:
{\small 
\[
 A_T^p:=
\lt(\hat\omega t+\psi_1 + \dd_1, \;  \hat\omega t+\psi_2+ \dd_2, \;    \hat\omega t+\psi_3+ \dd_3;
  \;\;    \hat\omega t+\psi_1+ \dd_1+\psi, \;   \hat\omega  t+\psi_2+ \dd_2+\psi, \; \hat\omega   t+\psi_3+ \dd_3+\psi \rt)^\top,
  \]
  }
where $\hat\omega$ is a  coupled stepping frequency that all six oscillators share 
\be{omega_hat}
\bal
\hat\omega &= \omega(\x) + \c_1H\lt(\psi; \x\rt) +\c_5 H(\psi_2 - \psi_1;\x)\\
&= \omega(\x)+ \c_2 H\lt(\psi; \x\rt) +\c_4 H(\psi_1-\psi_2;\x)+\c_7 H\lt(\psi_3-\psi_2; \x\rt)\\
&= \omega(\x) + \c_3H\lt(\psi; \x\rt) +\c_6 H\lt(\psi_2-\psi_3; \x\rt), 
\eal
\ee
 $\psi_1, \psi_2, \psi_3$ are corresponding relative phases, and 
$\psi$ is the corresponding constant contralateral phase difference in approximate gaits. 
Note that the equalities in $\hat\omega$ hold by Assumptions \ref{eta_assumption} and \ref{assu:balance}. 

The $\dd_i$'s are perturbations to the legs' phases due to the heterogeneity and are the solutions of 
\[
\left(\begin{array}{c} \tilde\omega_1 \\ \tilde\omega_2 \\ \tilde\omega_3\end{array}\right) = 
\mathcal{L}(\psi_1, \psi_2, \psi_3)
\left(\begin{array}{c}\dd_1 \\\dd_2 \\\dd_3\end{array}\right), 
\]
where 
\[
\mathcal{L} (\psi_1,\psi_2,\psi_3)=  \left(\begin{array}{ccc}
 \c_5 H'(\psi_2-\psi_1;\x) & -\c_5 H'(\psi_2-\psi_1;\x) & 0 \\
-\c_4 H'(\psi_1- \psi_2;\x) & \c_4 H'(\psi_1-\psi_2;\x)+\c_7 H'(\psi_3-\psi_2;\x) & -\c_7 H'(\psi_3-\psi_2;\x)\\
 0 & -\c_6 H'(\psi_2-\psi_3;\x) & \c_6 H'(\psi_2-\psi_3;\x)\end{array}\right), 
\]
and $H'$ denotes the derivative of $H$ w.r.t. its first argument. 
The matrix $\mathcal{L}$ is singular, so 
we let $(\delta_1,\delta_2,\delta_3)^\top=\mathcal{L}^+( \tilde\omega_1, \tilde\omega_2, \tilde\omega_3)^\top$, 
where $\mathcal{L}^+$ is the generalized inverse (pseudoinverse) of $\mathcal{L}$, see \cite{ginverse1971}. 

The following choices of the relative phases $\psi_1, \psi_2, \psi_3$, and the contralateral phase difference $\psi$ give 
four different versions of approximate tetrapod gaits and an approximate tripod gait. 
 \begin{enumerate}[leftmargin=*]
 
 \item  The approximate forward right tetrapod gait, denoted by $\FR^p$,  corresponds to $A_T^p$ with 
  $\psi_1=2/3, \psi_2= 0 , \psi_3=1/3$, and $\psi = 2/3$. 
 
 \item The approximate forward left tetrapod gait, denoted by $\FL^p$, corresponds to $A_T^p$ with
  $\psi_1=2/3, \psi_2= 0 , \psi_3=1/3$, and $\psi = 1/3$. 

\item The approximate backward right tetrapod gait, denoted by $\BR^p$, corresponds to $A_T^p$ with
  $\psi_1=1/3, \psi_2= 0 , \psi_3=2/3$, and $\psi = 1/3$. 

\item The approximate backward left tetrapod gait, denoted by $\BL^p$, corresponds to $A_T^p$ with
  $\psi_1=1/3, \psi_2= 0 , \psi_3=2/3$, and $\psi = 2/3$. 

\item The approximate tripod gait, denoted by $\T^p$,  corresponds to $A_T^p$ with
  $\psi_1=1/2, \psi_2= 0 , \psi_3=1/2$, and $\psi = 1/2$. 
\end{enumerate}

 \end{definition}
 In forward tetrapod gaits a wave of swing phases runs from hind to front legs and 
 in backward tetrapod gaits the swing phases run from front to hind legs, see \cite[Figure 4]{SIADS2018}. 
 
The matrix  $\mathcal{L}$ in Definition \ref{defn1} can be derived by 
substituting $A_T^p$ into Equation~(\ref{eq.osc1}) and 
approximating $H$ by the first two terms of its Taylor expansion. 
For instance, substituting $A_T^p$ into the first equation of (\ref{eq.osc1}), we get
\be{}
\bal
 \hat\omega = \dot\phi_1 &= \omega + \tilde\omega_1+ \c_1 H (\psi;\x) + \c_5 H(\psi_2+\dd_2 - \psi_1 -\dd_1;\x)\\
&=\omega + \tilde\omega_1+ \c_1 H (\psi;\x) + \c_5 H(\psi_2- \psi_1;\x) +\c_5 H' (\psi_2-\psi_1;\x) (\dd_2-\dd_1) + \mathcal{O}(\dd_2-\dd_1)^2. \nonumber
\eal
\ee
By substituting $\hat\omega=\omega(\x) + \c_1H\lt(\psi; \x\rt) +\c_5 H(\psi_2 - \psi_1;\x)$ into the above equation, 
$\tilde\omega_1$ can be approximated by $-\c_5 H' (\psi_2-\psi_1;\x) (\dd_2-\dd_1)$, 
which gives the first row of $\mathcal{L}$. The other rows are found in the same way.
 
Note that when $\dd_i=0$, i.e., in the homogeneous system,  two  (resp. three) legs 
swing simultaneously in tetrapod (resp. tripod) gaits but when $\dd_i\neq0$,
the corresponding legs do not swing exactly together due to the small perturbations
$\dd_i$, so we call them \textit{approximate tetrapod (resp. tripod) gaits}. 
 
In  \cite{SIADS2018} we showed  that Equation~(\ref{eq.osc1}) admits a 
solution at  a tetrapod gait, when the speed parameter $\x$ is small, 
and a solution at a tripod gait, when  $\x$ is large.
To connect tetrapod gaits to tripod gaits, we defined transition gaits \cite[Definition 2]{SIADS2018}. 
In what follows we generalize those definitions for heterogeneous models to connect approximate tetrapod gaits to approximate tripod gaits. 

 \begin{definition} [Approximate transition gaits] \label{defn2}
For any fixed number $\et\in[0,1/6]$ 
 \begin{enumerate}[leftmargin=*]
 
 \item  The approximate forward right transition gait, denoted by $\FR^p(\et)$,  corresponds to $A_T^p$ with 
  $\psi_1=2/3 - \et, \psi_2= 0 , \psi_3=1/3+\et$, and $\psi = 2/3-\et$. 
 
 \item The approximate forward left transition gait, denoted by $\FL^p(\et)$, corresponds to $A_T^p$ with
  $\psi_1=2/3-\et, \psi_2= 0 , \psi_3=1/3+\et$, and $\psi = 1/3+\et$. 

\item The approximate backward right transition gait, denoted by $\BR^p(\et)$, corresponds to $A_T^p$ with
  $\psi_1=1/3+\et, \psi_2= 0 , \psi_3=2/3-\et$, and $\psi = 1/3+\et$. 

\item The approximate backward left transition gait, denoted by $\BL^p(\et)$, corresponds to $A_T^p$ with
  $\psi_1=1/3+\et, \psi_2= 0 , \psi_3=2/3-\et$, and $\psi = 2/3-\et$.
  
\end{enumerate}
\end{definition}

As $\x$ increases, $\et=\et(\xi)$, the solution of Equation~\eqref{eta}, varies from $0$ to $1/6$. 
Therefore, at low speeds, when $\et=0$, $\FR^p(\et)$ (resp. $\FL^p(\et)$, $\BR^p(\et)$, and $\BL^p(\et)$) 
corresponds to  the approximate forward right 
(resp. forward left, backward right, and backward left) transition gait
and as $\x$ increases and $\et$ approaches $1/6$, 
all the approximate transition gaits tend to an approximate tripod gait. 
 
In what follows, we will see how certain properties of $H$ allow us 
to reduce 6 phase equations to 3 ipsilateral equations. 

In both approximate tetrapod and tripod gaits, the phase difference between
the left and right legs, denoted by $\psi$, is constant, and is either equal to
$\psi= 1/3$ or $\psi=2/3$ (in tetrapod gaits) or $\psi=1/2$ (in the tripod gait).   
In addition, we observe that the phase differences between the left and  right legs
in approximate transition gaits are  constant and equal to $2/3-\et$ or $1/3+\et$. 
For steady states, this assumption is supported by experiments for tripod gaits
\cite{Fuchs14}, and by simulations for tripod and tetrapod gaits in the bursting
neuron model \cite[Figures 4 and 5]{SIADS2018}. 

We make a further simplifying assumption that the steady state contralateral
phase differences remain constant for all $t$.
\begin{assumption}
\label{constant_contralateral_phase_diff}
The phase differences between the left and right legs are constant. For $i=1,2,3$, 
 \[ \phi_{i+3} - \phi_i=  {2}/{3} - \et \quad \mbox{or} \quad  \phi_{i+3} - \phi_i = {1}/{3} + \et. \] 
\end{assumption}

As discussed earlier, the coupling function computed from the bursting neuron model satisfies 
Assumption~\ref{eta_assumption} (Equation (\refeq{eta})) and thus allows reduction to 3
ipsilateral equations, as we now describe. 

\subsection{Phase differences model}
\label{phase_diffs}

In this section, the goal is to reduce 
the 6 equations (\ref{eq.osc1}) to 2 equations on a 2-torus. 

By Assumptions~\ref{eta_assumption}, \ref{external_input_assumption}, and \ref{constant_contralateral_phase_diff},  
Equation~(\ref{eq.osc1}) can
be reduced to the following 3 equations describing the right legs' motions:
\begin{subequations}\label{eq.osc5}
\begin{align}
& \dot{\phi}_1 \;=\;\omega(\x) + \tilde\omega_1+ \c_1H\lt(\frac{2}{3} - \et; \x\rt) + \c_5H(\phi_2 - \phi_1; \x), \\
& \dot{\phi}_2  \;=\; \omega(\x) + \tilde\omega_2+ \c_2H\lt(\frac{2}{3} -\et; \x\rt) + \c_4H(\phi_1 - \phi_2; \x) + \c_7H(\phi_3 - \phi_2; \x), \\
& \dot{\phi}_3  \;=\; \omega(\x) + \tilde\omega_3+ \c_3H\lt(\frac{2}{3} -\et; \x\rt)  + \c_6H(\phi_2 - \phi_3; \x) .
\end{align}
\end{subequations}

Because only phase differences appear in the vector field, we may define
\[
\theta_1:=\phi_1-\phi_2 \quad \mbox{and}\quad \theta_2:=\phi_3-\phi_2,
\]
so that, from Equations~(\ref{eq.osc5}), the following equations describe
the dynamics of $\theta_1$ and $\theta_2$:
\begin{subequations} \label{torus:equation}
\begin{align}
&\dot\theta_1=   \tilde\omega_1 -  \tilde\omega_2 + (\c_1 - \c_2)H\lt(\frac{2}{3} - \et; \x\rt) + \c_5 H(-\theta_1; \x) - \c_4 H(\theta_1; \x) - \c_7 H(\theta_2; \x), \\ 
&\dot\theta_2=   \tilde\omega_3 -  \tilde\omega_2 + (\c_3 - \c_2)H\lt(\frac{2}{3} - \et; \x\rt) + \c_6 H(-\theta_2; \x) - \c_4 H(\theta_1; \x) - \c_7 H(\theta_2; \x),
\end{align}
\end{subequations}
where the $\tilde\omega_i$'s and $H$ are defined in Equations~(\ref{tilde_omega_i}) 
and (\ref{H:interconnected:bursting:neuron}), respectively. 

Note that  Equations~(\ref{torus:equation}) are $1$-periodic in both variables,
i.e., $(\theta_1, \theta_2) \in \mathbb{T}^2$, where $\mathbb{T}^2$ is a 2-torus. 

In Equations~(\ref{torus:equation}), the approximate tripod gait $\T^p$ corresponds approximately 
to the fixed point $(1/2, 1/2)$, the approximate forward tetrapod gaits, $\FR^p$ and $\FL^p$,
correspond approximately to the fixed point $(2/3, 1/3)$, 
the approximate backward tetrapod gaits, $\BR^p$ and $\BL^p$,
correspond approximately to the fixed point $(1/3, 2/3)$, and the approximate transition gaits, 
$\FR^p(\et)$ and $\FL^p(\et)$ (resp.  $\BR^p(\et)$ and $\BL^p(\et)$), 
correspond approximately to $(2/3-\et, 1/3+\et)$ (resp. $(1/3+\et, 2/3-\et)$). 
See \cite{Yeldesbay2018} for similar definitions of tetrapod and tripod gaits on a torus. 

In  \cite{SIADS2018} we observed that 
when $\x$ is small, the forward tetrapod gaits are not the only solutions and there exist 
some other stable and unstable solutions  (e.g. stable or unstable backward tetrapod and unstable tripod gaits). 
We showed that as $\x$ increases, one stable tripod gait emerges, through a degenerate bifurcation. 
In the present work, we show how heterogeneity,  $\Ii$,  can break the
degenerate bifurcation into separate saddle-node bifurcations such that at low
speed, only one stable  (either forward or backward tetrapod) gait exists. 
We are primarily interested in the existence of  approximate {\it forward} tetrapod gaits, 
since they have been observed widely in insects (Section \ref{forward} below).
However, backward tetrapod gaits have also been 
observed in backward-walking fruit flies
\cite[Supplementary Materials, Figure S1]{Bidaye-Science14},  and so 
in Section \ref{backward} we show that $\Ii$ can be chosen such that 
only approximate backward tetrapod gaits exist at low speed. 

\section {Main results}
\label{Main_Result}

In this section, we fix a low speed parameter (e.g., $I_{ext} = 35.65$) where 
the bursting neuron model (\ref{BN}) can generate tetrapod gaits. Also, we
 assume that the balance condition holds, so several fixed points including
 the forward $(2/3,1/3)$ and backward $(1/3,2/3)$ tetrapod gaits exist. 
The main goal is to show how adding small heterogeneous external 
currents $\Ii$ can successively remove fixed points 
on the torus while preserving the stable forward  (see Section \ref{forward} below) 
or stable backward   (see Section \ref{backward} below) tetrapod gaits, respectively. 

\begin{figure}[h!]
\begin{center}
\includegraphics[scale=.1]{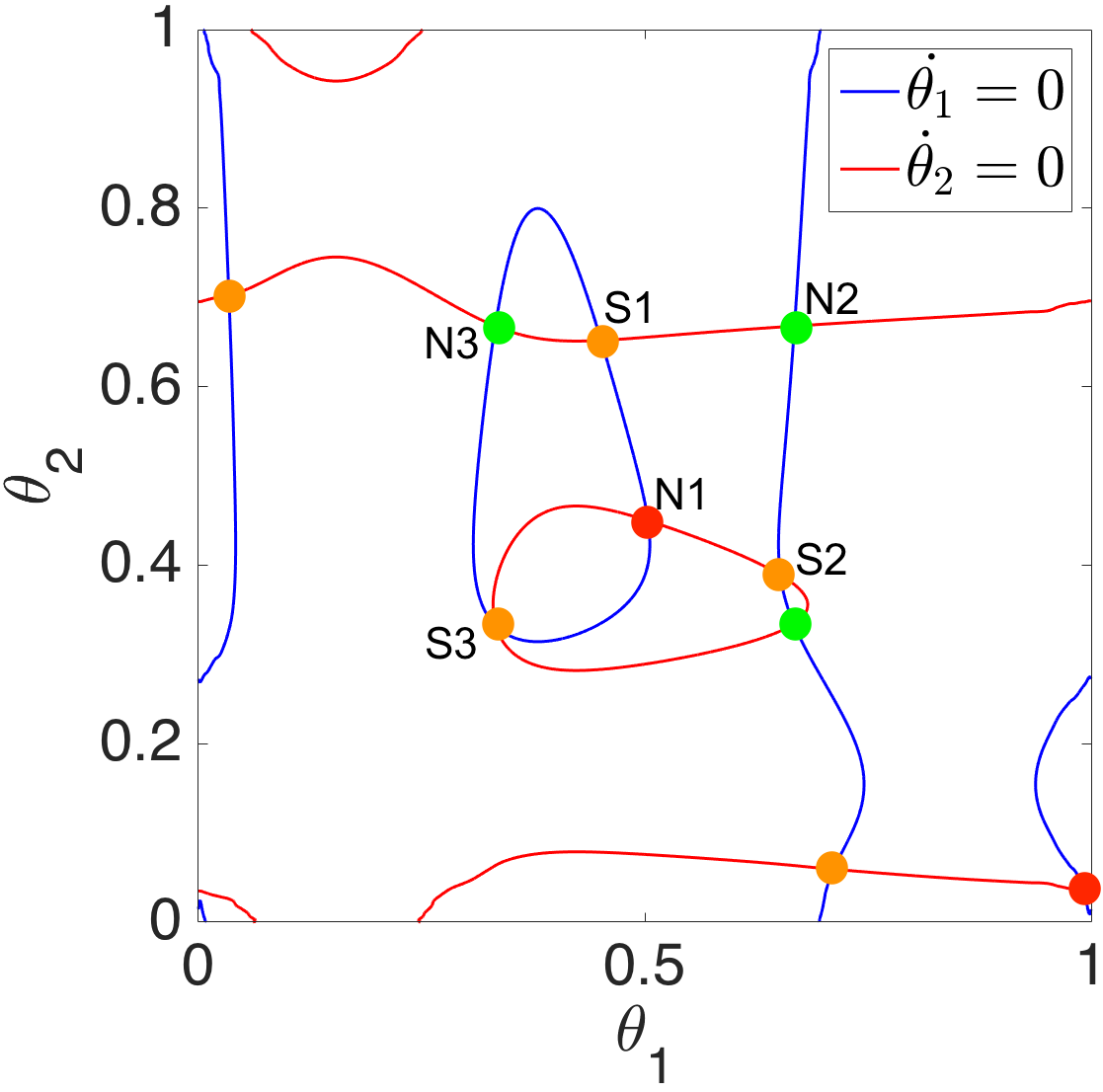}\;
\includegraphics[scale=.1]{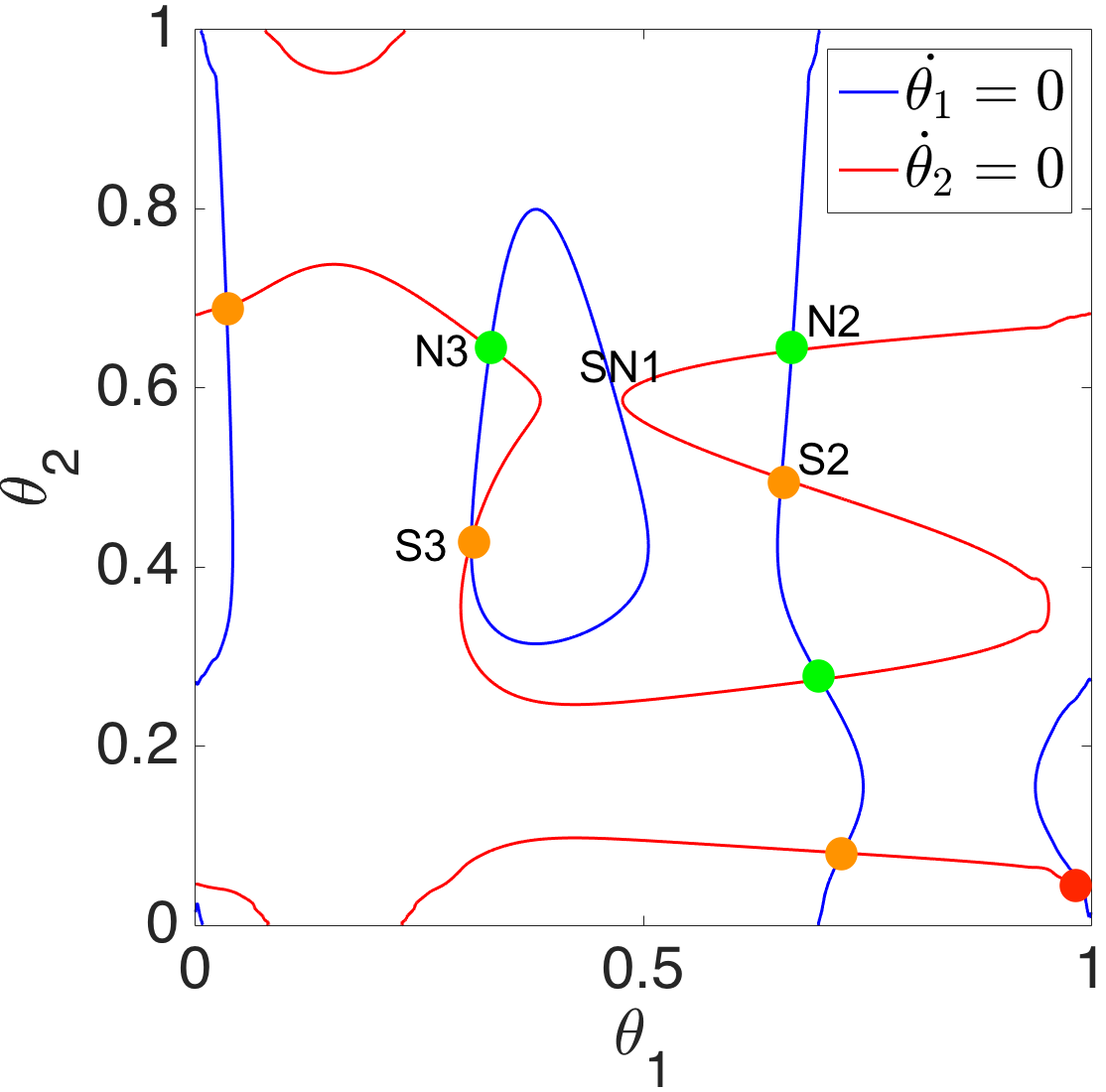}\;
\includegraphics[scale=.1]{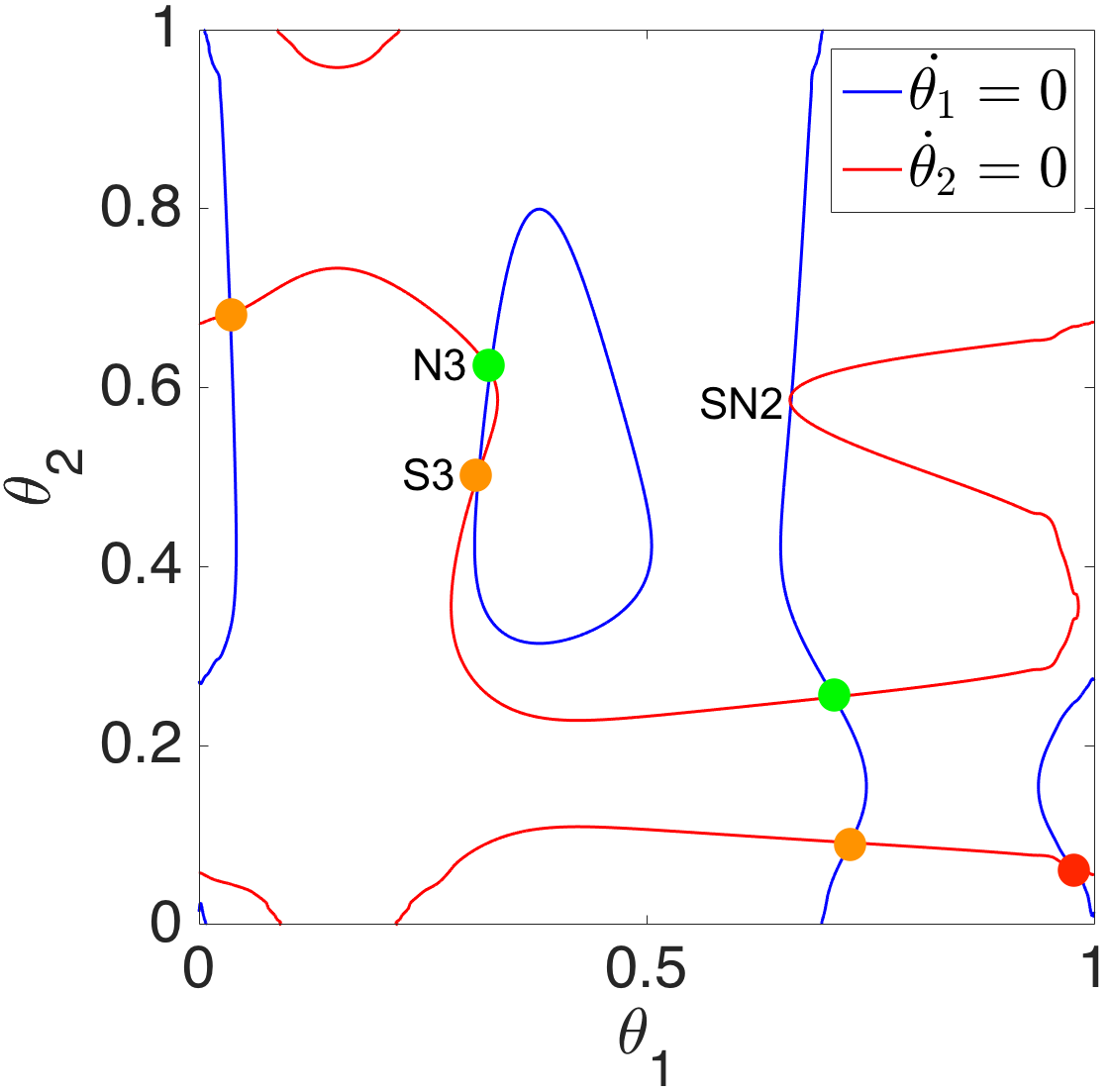}\;
\includegraphics[scale=.1]{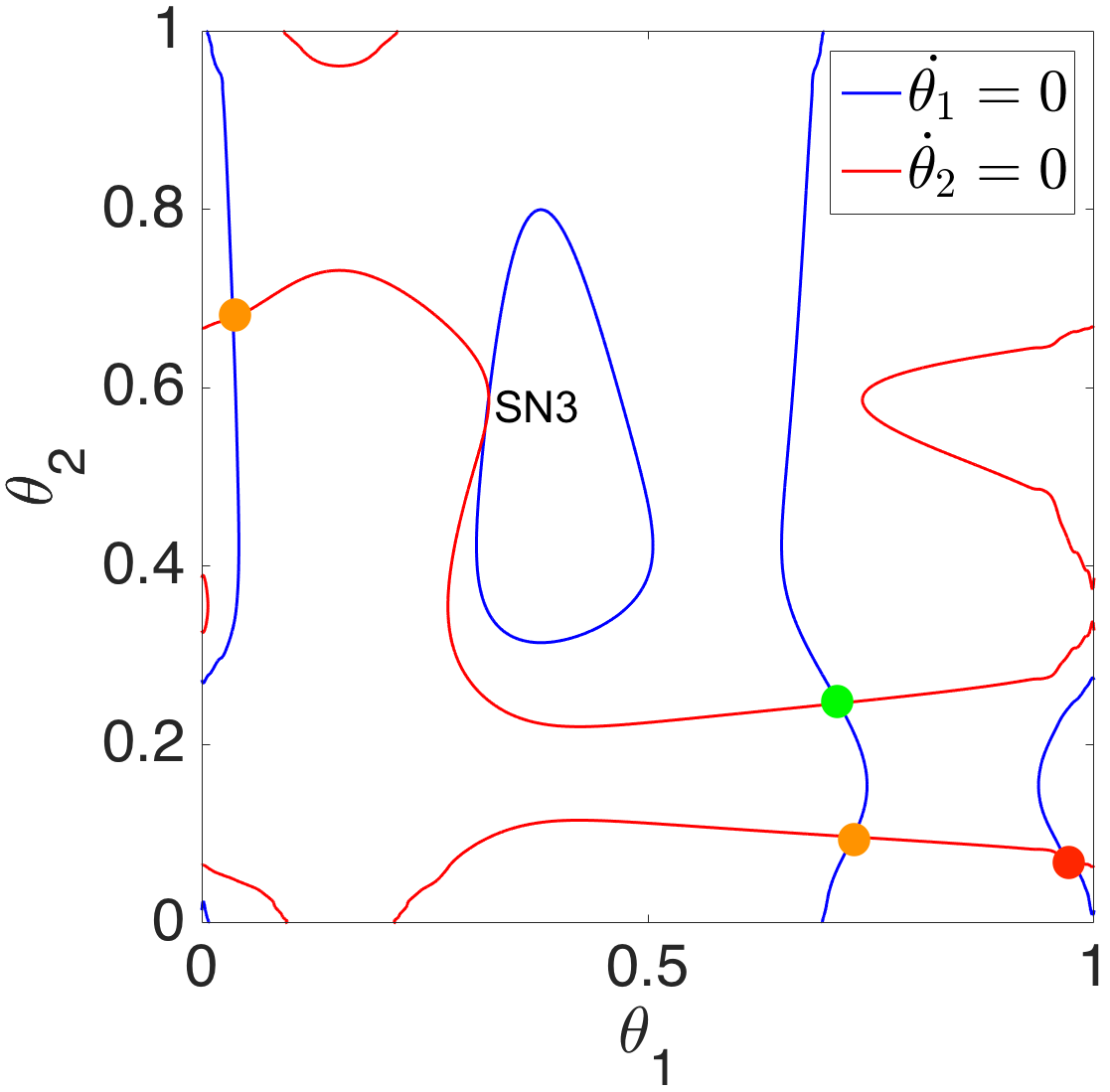}
\end{center}
\caption{(Left to right) Nullclines of Equations~(\ref{torus:equation})
when $\c_i$'s satisfy Equation~(\ref{random_balance_coupling}), $I_{ext} = 35.65$,
and $\delta I \approx 0, 0.02, 0.032, 0.038$, respectively. $\dot\theta_1 = 0$
nullcline shown in blue; $\dot\theta_2 = 0$ nullcline shown in red.
 Green dots indicate sinks, red dots are sources, and orange dots are saddle points. 
At  $\delta I \approx 0.02$, the first saddle-node bifurcation (shown by SN1) occurs 
and the unstable tripod gait (shown by N1) disappears together with a saddle point (shown by S1).
At $\delta I \approx 0.032$, the second saddle-node bifurcation (shown by SN2) occurs and a stable 
fixed point (shown by N2) disappears together with a saddle point (shown by S2).
Finally, at $\delta I \approx 0.038$, the third saddle-node bifurcation (shown by SN3) occurs and 
the stable backward tetrapod gait (shown by N3) disappears together with a saddle point (shown by S3). 
A single stable approximate forward tetrapod gait remains together with a source and two saddle points.  
 }
\label{NC_hetero_random_bifurcations}
\end{figure}

For example, consider the following randomly generated coupling strengths $\c_i$
that satisfy the balance condition
\be{random_balance_coupling}
\c_1=0.8147, \; \c_2=0.9058,\; \c_3=0.1270,\;  \c_4=0.9134,\; \c_5=1.6368, \; \c_6=2.3245,\;  \c_7=0.6324.
\ee
For any $t$, let $\delta I :=I_{ext}^1(t) = I_{ext}^2(t),$ and $I_{ext}^3(t) = 0$, and vary the heterogeneity parameter $\delta I$ from 0. 
Figure \ref{NC_hetero_random_bifurcations} (left to right) shows the 
nullclines of Equation~(\ref{torus:equation})
with $\delta I \approx 0, 0.02 , 0.032 , 0.038$, respectively. 
In this example, at $\delta I=0$, for which the model is homogeneous, 
there exist 3 stable sinks, 2 unstable sources 
and 5 saddle points. As the heterogeneity parameter $\delta I$ increases, 
3 saddle-node bifurcations occur at
approximately $\delta I \approx  0.02, 0.032 , 0.038$ and
one stable fixed point located at $\approx (0.71,0.25)$ remains, which
 corresponds to a stable approximate forward tetrapod gait. The other 3
remaining fixed points are a source and 2 saddle points.

 Note that the balance condition is sufficient for
the existence of forward and backward tetrapod gaits but does not guarantee
their stability. 
Assuming that the balance condition holds, in \cite[Proposition 6 (resp. Proposition 7)]{SIADS2018}
we proved that the forward (resp. backward) tetrapod gait is always
stable if $\c_1=\c_2=\c_3$, and  $\alpha:= \frac{c_4}{c_4+c_7} 
< \alpha_{max}$ (resp. $\alpha > \alpha_{min}$), where $\alpha_{max}$
(resp. $\alpha_{min}$) can be computed from the derivatives of $H$:
\be{alpha_alphamax} 
\alpha_{max} (\x):=\frac{H'\lt(\frac{1}{3}; \x\rt)}{H'\lt(\frac{1}{3}; \x\rt)-H'\lt(\frac{2}{3}; \x\rt)}, \quad
\alpha_{min} (\x):= \frac{H'\lt(\frac{2}{3}; \x\rt)}{H'\lt(\frac{2}{3}; \x\rt)-H'\lt(\frac{1}{3}; \x\rt)}. 
\ee

Recall that here we fixed the speed parameter $\x = I_{ext}$, so $\alpha_{max}$ and $\alpha_{min}$
are constant. 
 
Without loss of generality, we can assume that one of the coupling strengths is equal to 1. 
For the rest of the paper we assume that $\c_4 = 1$, the balance condition \eqref{balance} holds, 
and $\c_1=\c_2=\c_3$. 
Therefore, by making a change of time variable that eliminates $\c_5=\c_6=1+\c_7 = 1/\alpha$, 
Equations~(\ref{torus:equation}) can be written as 
 \begin{subequations} \label{eq.osc.simplified}
\begin{align}
&\dot\theta_1=  \alpha(\tilde\omega_1 -  \tilde\omega_2) +H(-\theta_1; \x) - \alpha  H(\theta_1; \x)- (1-\alpha)  H(\theta_2; \x), \\ 
&\dot\theta_2 = \alpha(\tilde\omega_3 -  \tilde\omega_2) + H(-\theta_2; \x) -  \alpha H(\theta_1; \x) - (1-\alpha) H(\theta_2; \x), 
\end{align}
\end{subequations}
which possess both forward and backward tetrapod gaits with stabilities dependent on the value of $\a$, i.e., for
\bi
\item $\a<\a_{min}$, the forward tetrapod gait is stable while the backward tetrapod gait is a saddle; 
 \item $\a_{min}<\a<\a_{max}$, both backward and forward tetrapod gaits are stable; 
 \item $\a>\a_{max}$, the backward tetrapod gait is stable while the forward tetrapod gait is a saddle.  
\ei
In Sections~\ref{forward} and \ref{backward}, we do
the following.  
\begin{description}[leftmargin=*]

\item[Section~\ref{forward}] 
We assume $\a<\a_{max}$ and let $\delta I_f \;:=\; I_{ext}^1(t) = I_{ext}^2(t) \;\geq\; 0,$ and $I_{ext}^3(t) = 0$. 
We show that for some small value of the heterogeneity parameter $\delta I_f$, Equations~(\ref{eq.osc.simplified}) 
possess only one stable forward tetrapod gait (together with a source and 2 saddle points). 

\item[Section~\ref{backward}] 
We assume $\a>\a_{min}$ and let $\delta I_b \;:=\; I_{ext}^2(t) = I_{ext}^3(t) \;\geq\; 0,$ and $ I_{ext}^1(t) = 0$. 
We show that for some small value of the heterogeneity parameter $\delta I_b$, Equations~(\ref{eq.osc.simplified}) 
possess only one stable backward tetrapod gait (together with a source and 2 saddle points). 

\end{description}

\subsection {Emergence of a unique stable forward tetrapod gait at low speed}
\label{forward}

We assume $\a<\a_{max}$ so that the forward tetrapod gait, $(2/3,1/3)$, 
is stable while the backward tetrapod gait can be either stable or a saddle, as described above.  
 For any $t$, let
 \be{external_Ii}
 \bal
\delta I_f \;:=\; I_{ext}^1(t) = I_{ext}^2(t) \;\geq\; 0, \quad I_{ext}^3(t) = 0, 
 \eal
 \ee
and consider the heterogeneity parameter  $\delta I_f$ as a bifurcation parameter. 

Choosing $\Ii$ as in Equation~(\ref{external_Ii}) implies $\tilde\omega_1 -  \tilde\omega_2  = 0$, and  
$ \tilde\omega_3 -  \tilde\omega_2 = -\delta I_f  \bar Z\leq0 $, 
where $\bar Z= \frac{1}{T}\int_0^T Z_v (\Gamma(t)) \; d t>0$ is 
the average of the phase response curve. Therefore, Equations~(\ref{eq.osc.simplified}) become
\begin{subequations} \label{eq.osc.simplified:forward}
\begin{align}
&\dot\theta_1=  H(-\theta_1; \x) - \alpha  H(\theta_1; \x)- (1-\alpha)  H(\theta_2; \x), \\ 
&\dot\theta_2 = - \alpha\delta I_f  \bar Z+ H(-\theta_2; \x) -  \alpha H(\theta_1; \x) - (1-\alpha) H(\theta_2; \x). 
\end{align}
\end{subequations}

Equations~(\ref{eq.osc.simplified:forward}) can also be obtained if 
 $ I_{ext}^1(t) = I_{ext}^2(t) = 0, \; \delta I_f = I_{ext}^3(t)< 0.$

We will show that when $\alpha_{min}<\alpha<\alpha_{max}$ (resp. $0<\alpha < \alpha_{min}$), 
as the bifurcation parameter $\delta I_f$ increases, 
Equations~\eqref{eq.osc.simplified:forward} lose 6 (resp. 4) fixed points 
through 3 (resp. 2) saddle-node bifurcations and keep only one stable approximate forward tetrapod gait. 
To show this, we consider two topologically different cases:

\begin{description}[leftmargin=*]

\item[\bm{$\alpha_{min}<\alpha<\alpha_{max}$}] At $\delta I_f = 0$, 
Equations~\eqref{eq.osc.simplified:forward} admit 10 fixed points (5 saddle points, 3 sinks, and 2 sources). 
As $\delta I_f$ increases, 3 saddle-node bifurcations occur and one sink (corresponding to the approximate forward tetrapod gait), 
one source and 2 saddle points remain.  

\item[\bm{$0<\alpha < \alpha_{min}$}] At $\delta I_f = 0$, 
Equations~\eqref{eq.osc.simplified:forward} admit 8 fixed points (4 saddle points, 2 sinks, and 2 sources). 
As $\delta I_f$ increases, 2 saddle-node bifurcations occur and one sink (corresponding to the approximate forward tetrapod gait), 
one source and 2 saddle points remain.  

\end{description}

\begin{figure}[h!]
\begin{center}
\includegraphics[scale=.1]{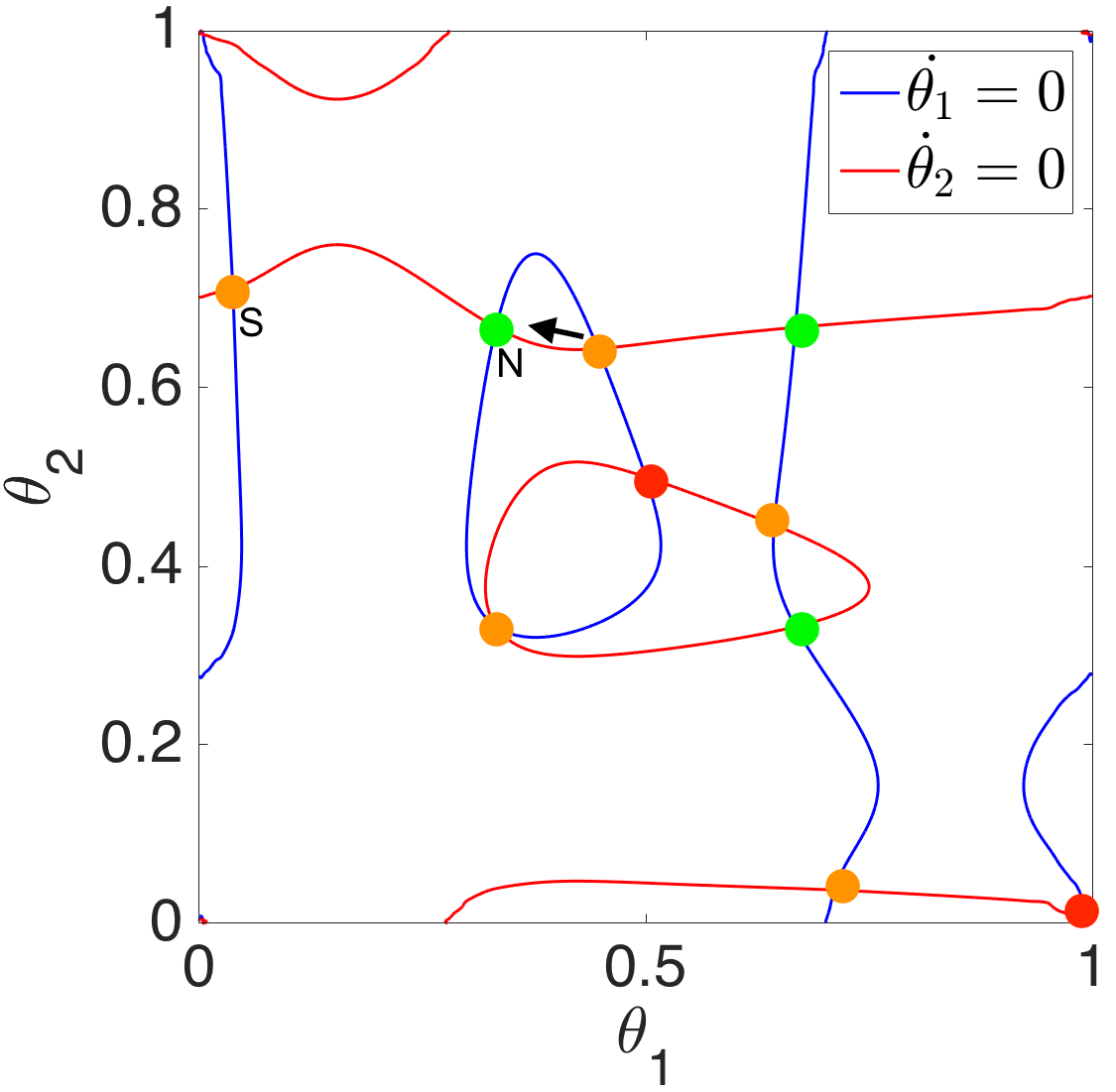}\;
\includegraphics[scale=.1]{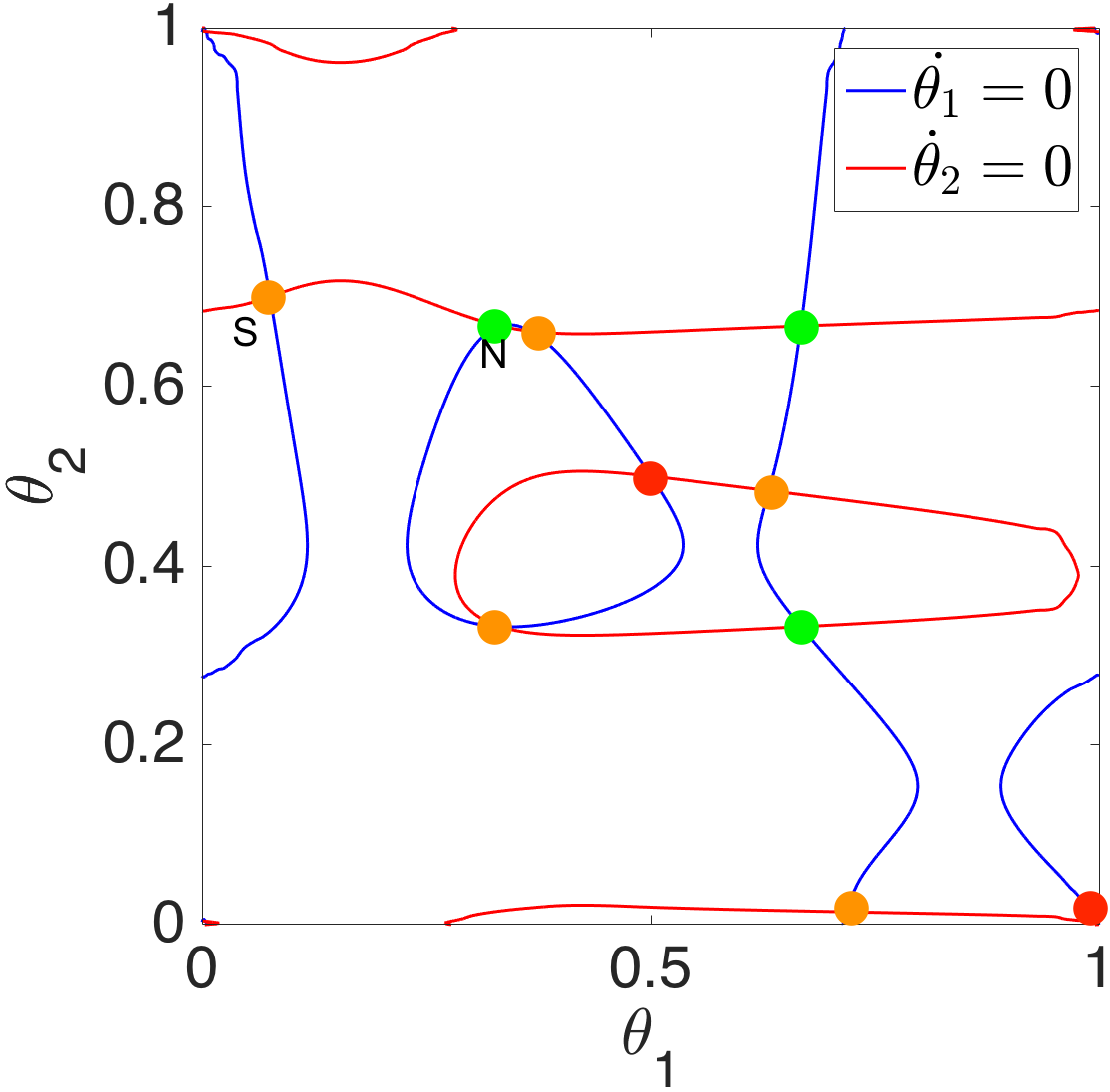}\;
\includegraphics[scale=.1]{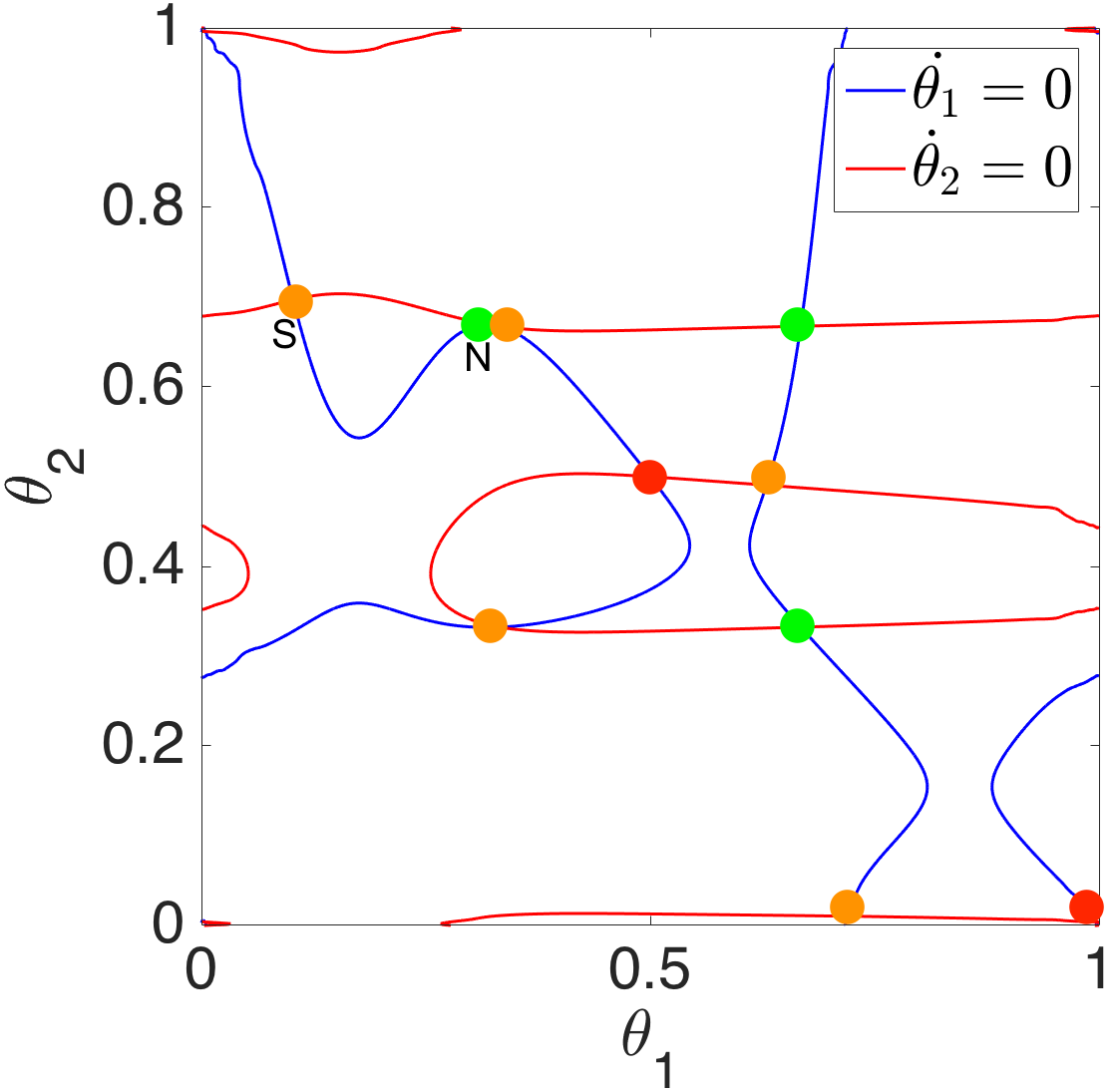}\;
\includegraphics[scale=.1]{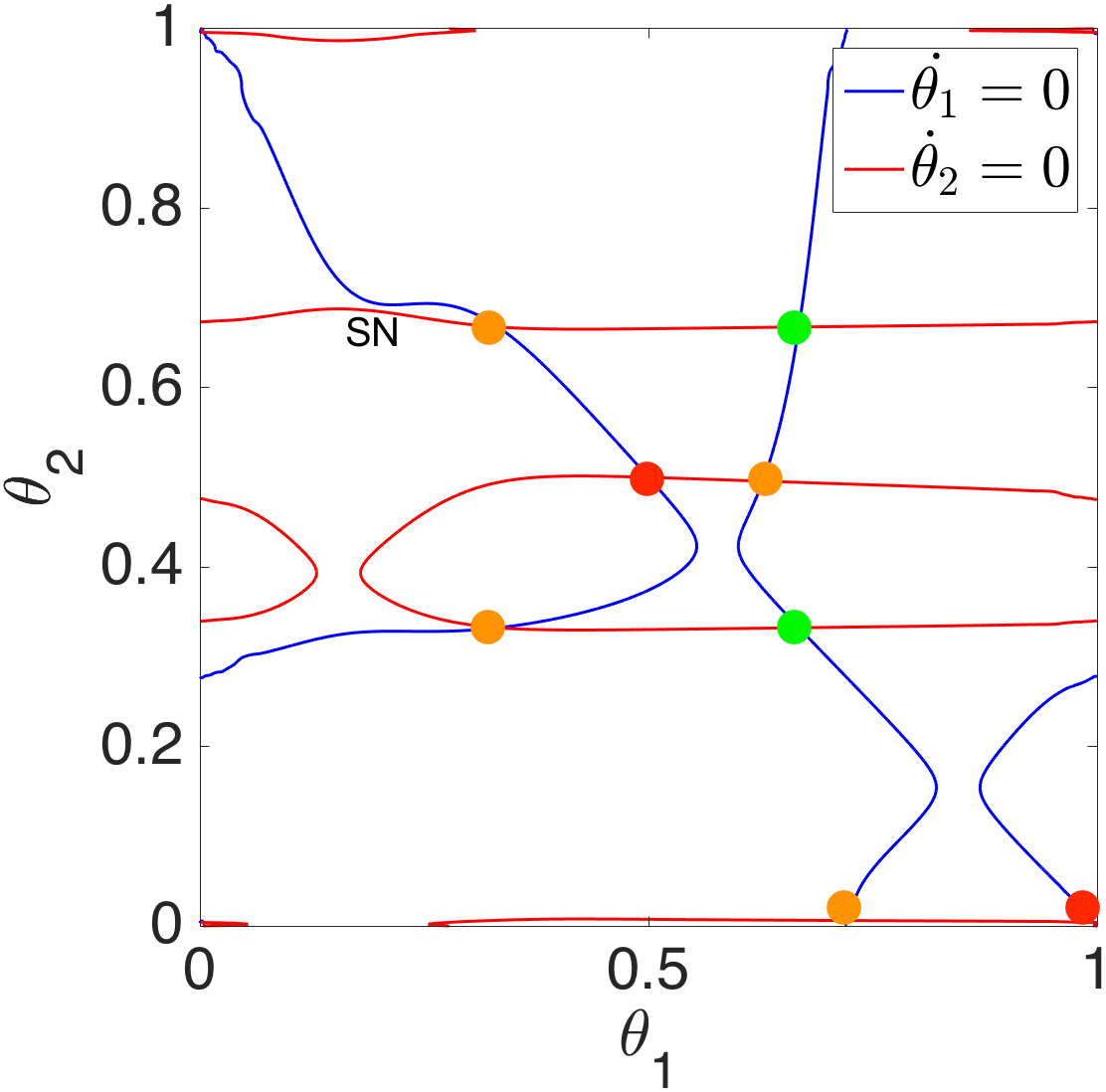}
\end{center}
\caption{(Left to right) Nullclines of 
Equations~(\ref{eq.osc.simplified:forward})  with $\alpha = 1/2, 1/5, 1/8,1/17$, 
$I_{ext} = 35.65$, and $\delta I_f =0$, respectively. 
At $\alpha = 1/2$, there exist 10 fixed points including one stable backward tetrapod gait shown by N. 
As $\alpha$ decreases, through a saddle-node bifurcation  (shown by SN)  the stable backward tetrapod gait disappears together with 
a saddle point shown by S and another saddle point becomes a backward tetrapod gait. 
Nullclines and fixed points are indicated as in Figure~\ref{NC_hetero_random_bifurcations}.}
\label{saddle_node_alpha_forward}
\end{figure}

Figure~\ref{saddle_node_alpha_forward} shows that as $\a$ decreases and approaches $\a_{min}$, 
through a  saddle-node bifurcation the stable backward tetrapod gait disappears 
 and one of the  saddle points moves toward the position of the  backward tetrapod gait, which is  shown by an arrow in 
 Figure~\ref{saddle_node_alpha_forward} (left). 
As $\alpha$ decreases,  the isolated $\dot\theta_1=0$ nullcline combines with the $\dot\theta_1=0$ nullcline
that encircles the torus and thereafter the number of fixed points reduces to 8 from 10.
 Therefore, when $0<\a<\a_{min}$, there are only 8 fixed points. 

\subsubsection*{3 saddle-node bifurcations: \bm{$\alpha_{min}<\alpha<\alpha_{max}$}}\label{forward_3SN}

\begin{figure}[h!]
\begin{center}
\includegraphics[scale=.1]{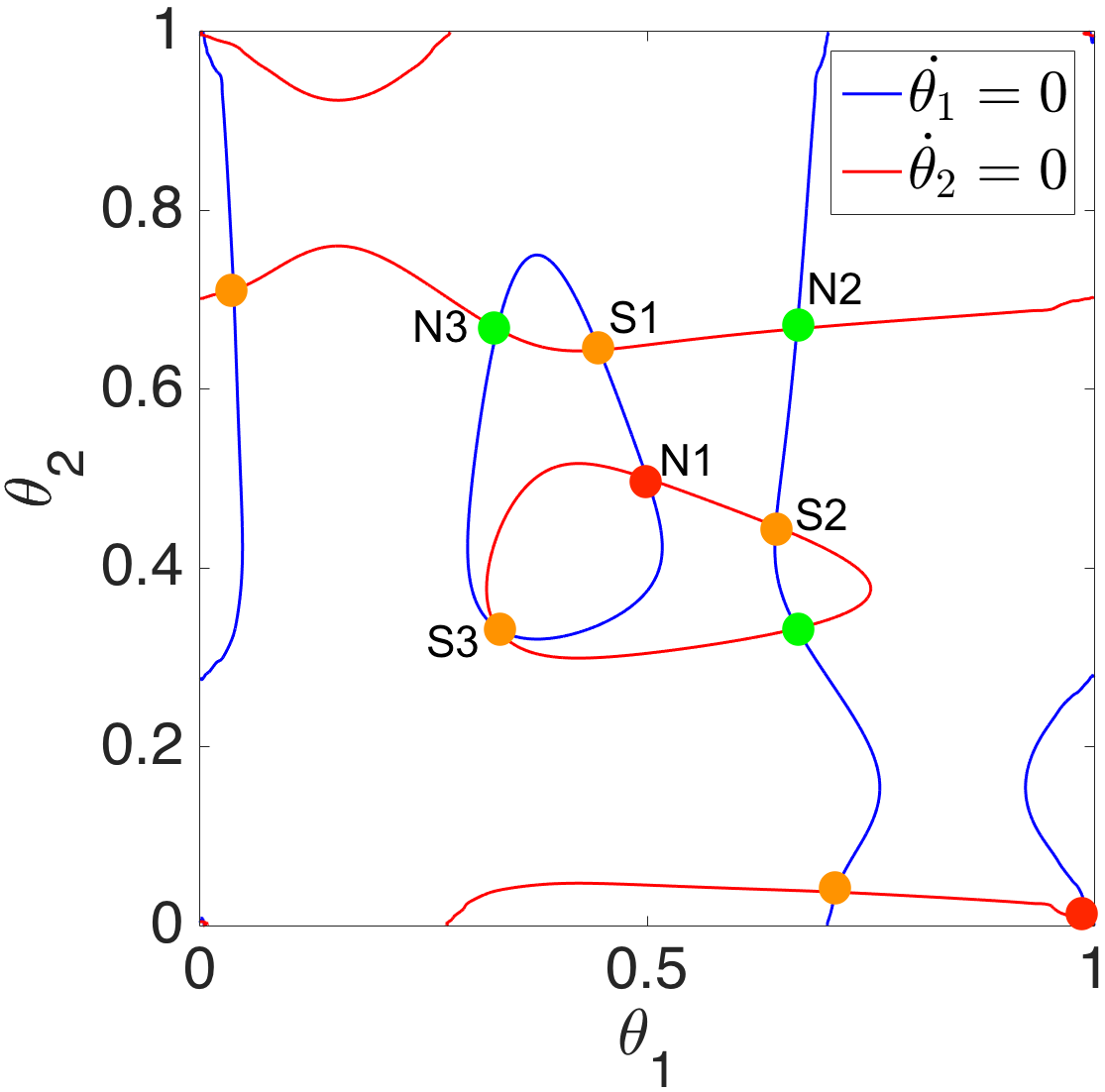}\;
\includegraphics[scale=.1]{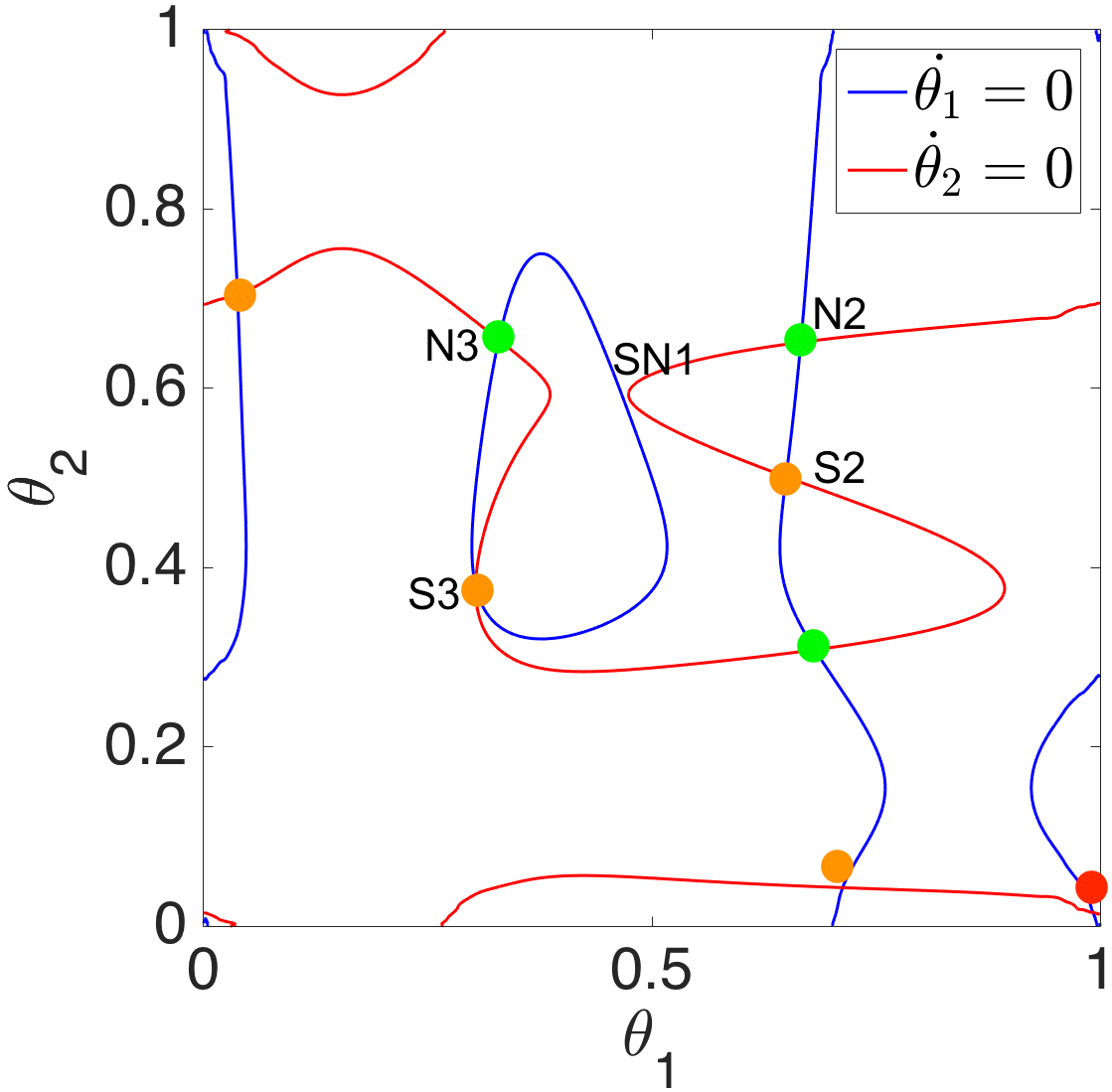}\;
\includegraphics[scale=.1]{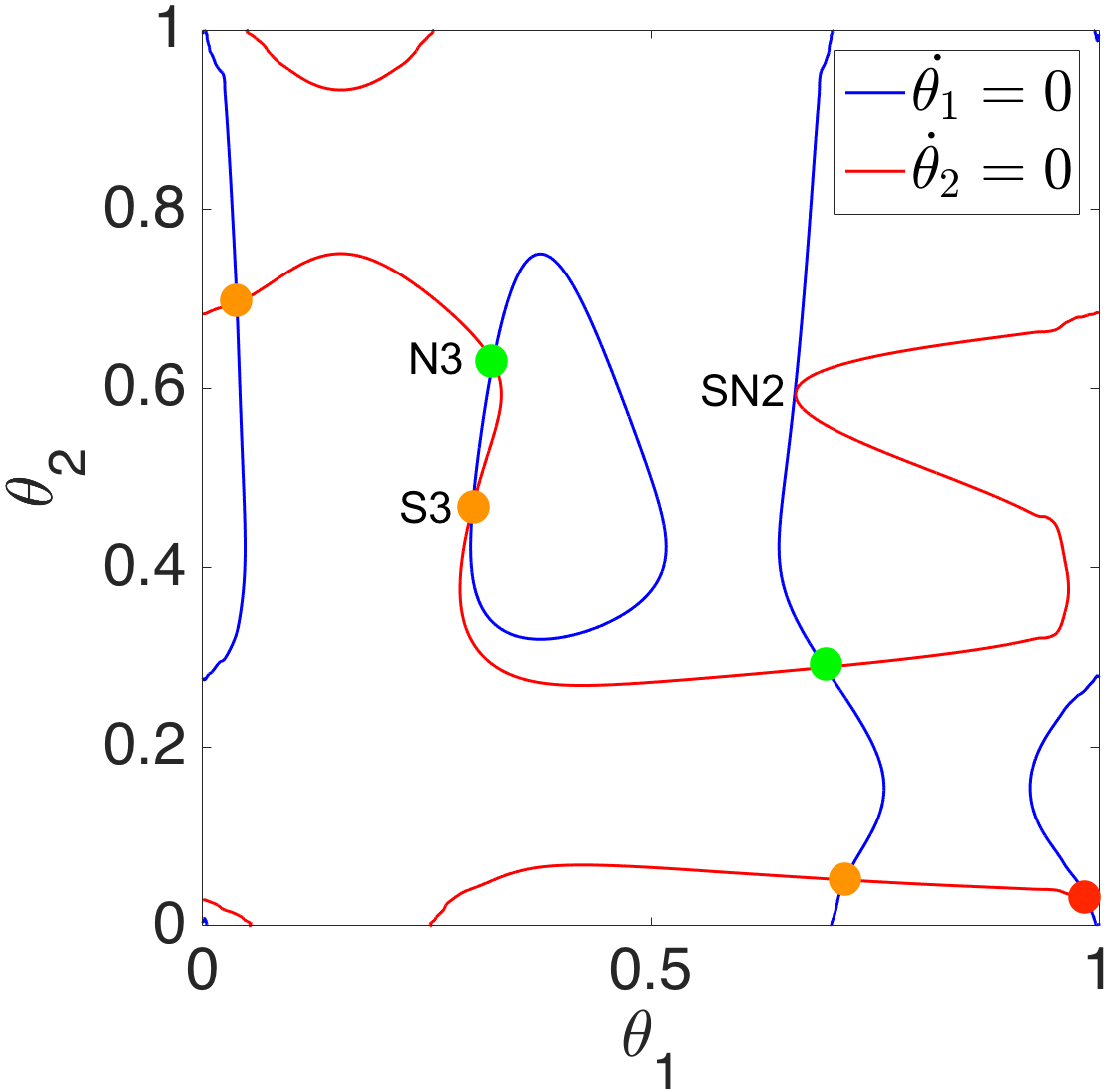}\;
\includegraphics[scale=.1]{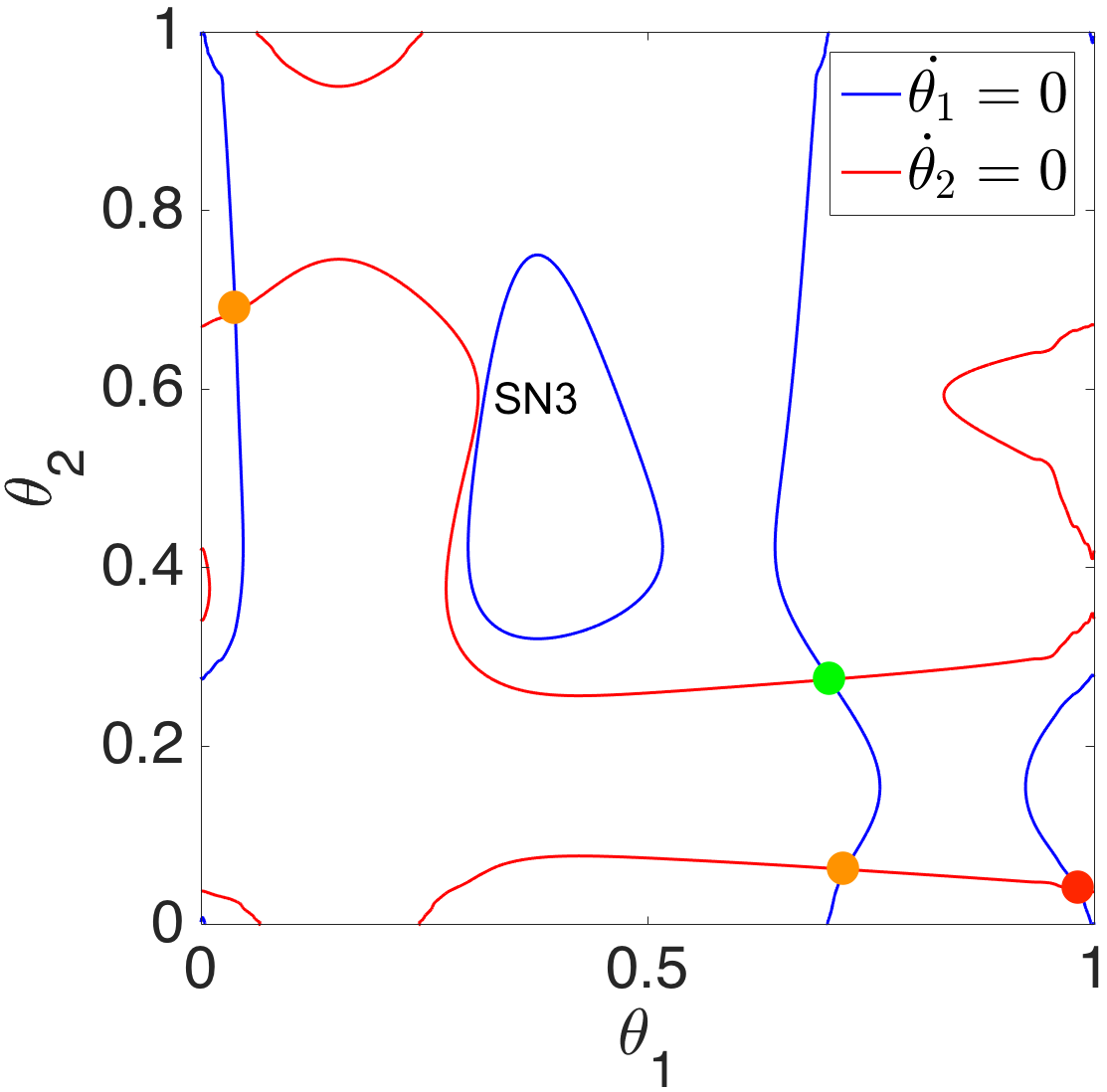}
\end{center}
\caption{(Left to right) Nullclines of 
Equations~(\ref{eq.osc.simplified}) with $\alpha=1/2$, 
$I_{ext} = 35.65$, and $\delta I_f \approx 0, 0.011, 0.025 , 0.037$,
respectively. At $\delta I_f =0$, there exist ten fixed points. 
At $\delta I_f \approx 0.011$, the first saddle-node bifurcation (shown by SN1) occurs 
and the unstable tripod gait (shown by N1) disappears together with a saddle point (shown by S1).
At $\delta I_f \approx 0.025$, the second saddle-node bifurcation (shown by SN2) occurs and a stable 
fixed point (shown by N2) disappears together with a saddle point (shown by S2).
 Finally, at $\delta I_f \approx 0.037$, the third saddle-node bifurcation (shown by SN3) occurs 
 and the stable backward tetrapod gait (shown by N3) disappears together with a saddle point  (shown by S3).
 A single stable approximate forward tetrapod gait remains together with a source and two saddle points.  
Nullclines and fixed points are indicated as in Figure~\ref{NC_hetero_random_bifurcations}.} 
\label{NC_hetero_alpha_12_bifurcations}
\end{figure}

Consider Equations~(\ref{eq.osc.simplified:forward}) with $\alpha_{min}<\alpha<\alpha_{max}$. 
Since the qualitative behavior of the solutions of Equations~(\ref{eq.osc.simplified:forward})  
with $\alpha_{min}<\alpha<\alpha_{max}$ are all similar, we show the results in an example with $\a=1/2$. 
As is clear from Equations~(\ref{eq.osc.simplified:forward}) and illustrated in
Figure \ref{NC_hetero_alpha_12_bifurcations}, choosing the heterogeneity of 
Equation~(\ref{external_Ii}) maintains the $\dot\theta_1 = 0$ nullclines and only
perturbs the $\dot\theta_2 = 0$ nullclines. This perturbation causes the topology
of the $\dot\theta_2 = 0$ nullclines to change, combining the isolated closed curve of the $\dot\theta_2=0$ nullcline 
with a nullcline that encircles the torus and thereafter reducing the number of fixed points.
In this case, where $\alpha=1/2$,  at $\delta I_f=0$, there exist 3 stable sinks, 2 unstable sources 
and 5 saddle points. As $\delta I_f$ increases, 3 saddle-node bifurcations occur at
approximately $\delta I_f \approx  0.011, 0.025 , 0.037$ and
one stable fixed point remains, which
 corresponds to a stable approximate forward tetrapod gait. The other 3
remaining fixed points are a source and 2 saddle points.

\subsubsection*{2 saddle-node bifurcations: \bm{$0<\alpha < \alpha_{min}$}}\label{forward_2SN}

\begin{figure}[h!]
\begin{center}
\includegraphics[scale=.1]{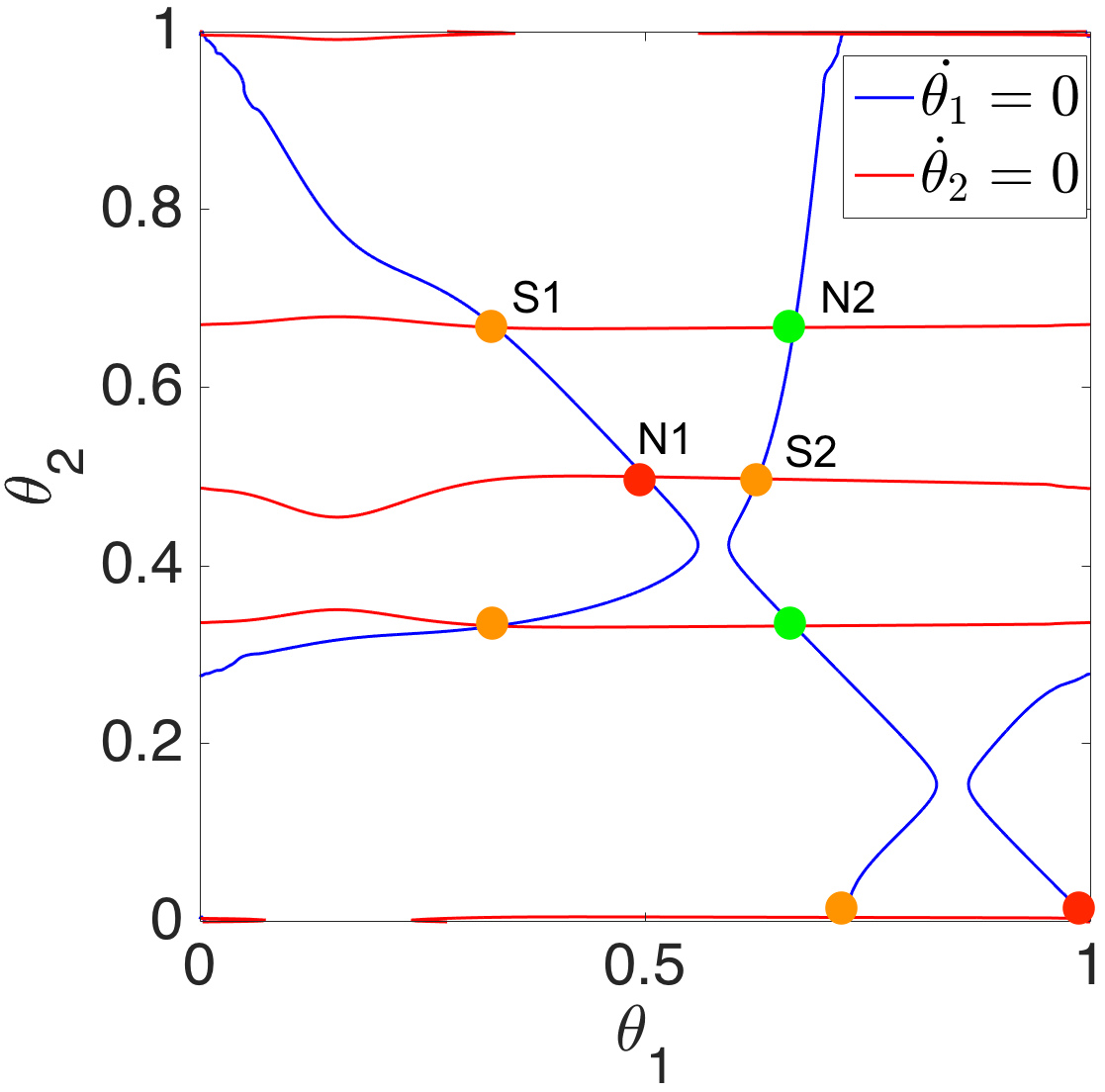}\;
\includegraphics[scale=.1]{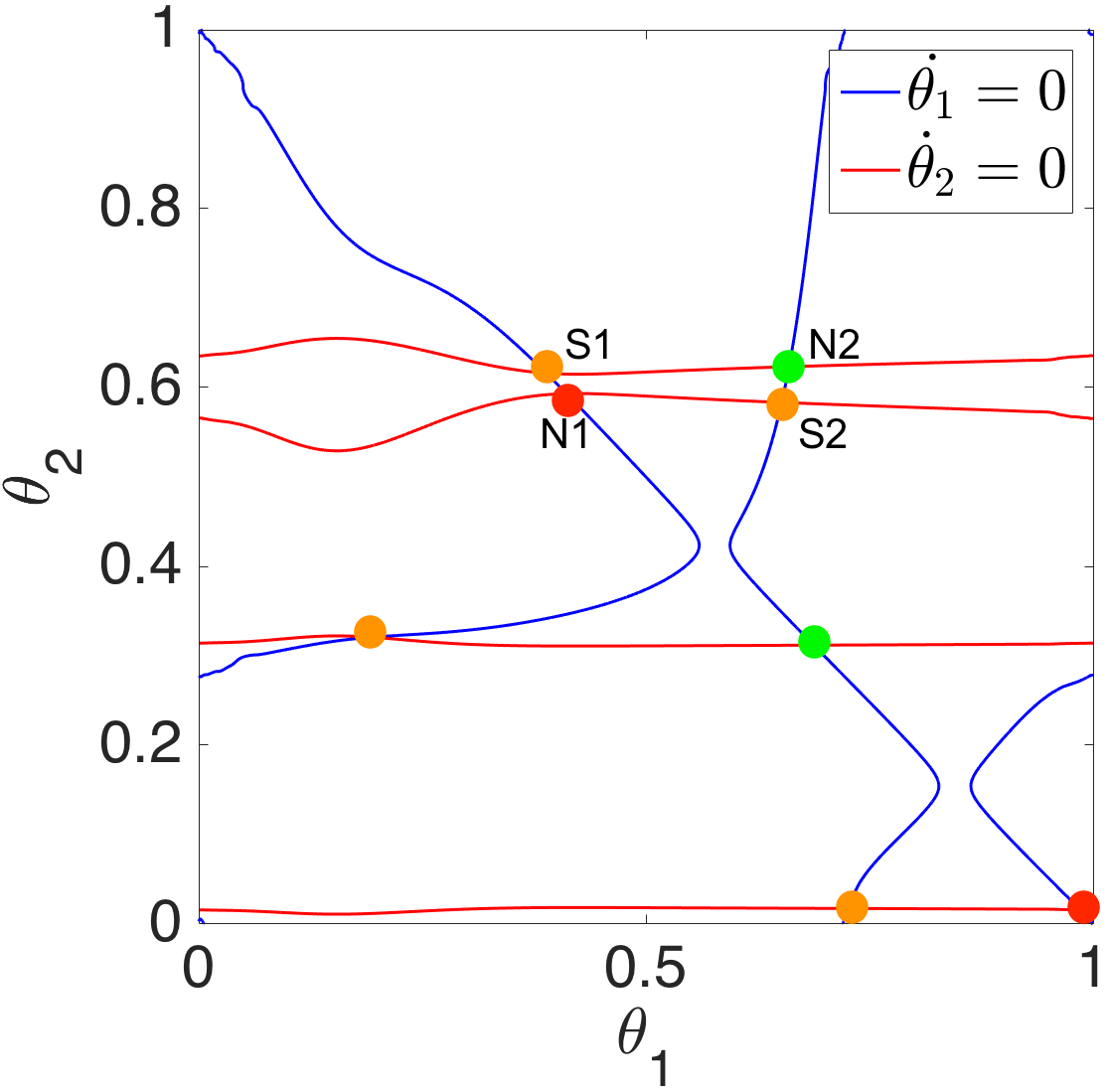}\;
\includegraphics[scale=.1]{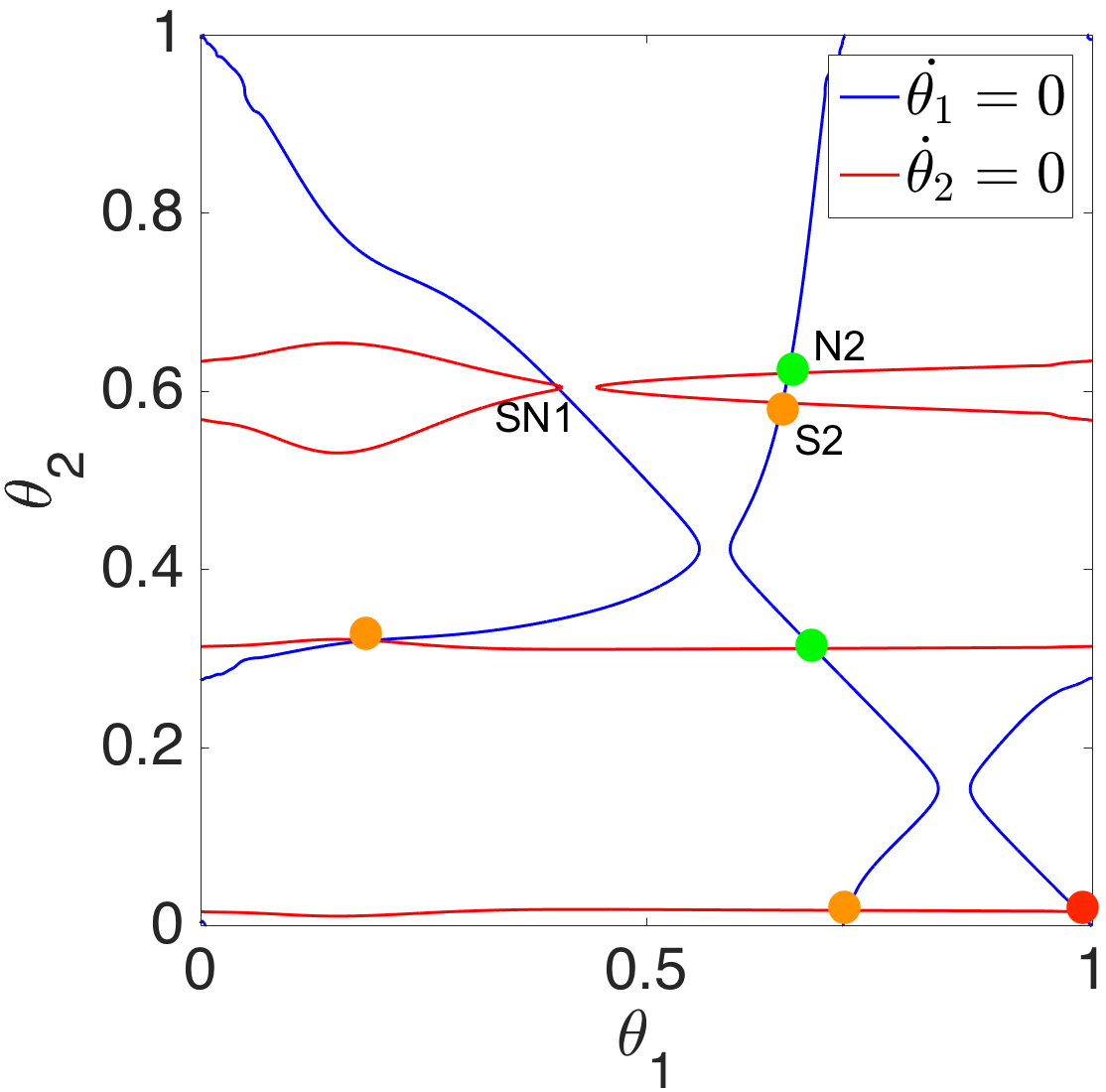}\;
\includegraphics[scale=.1]{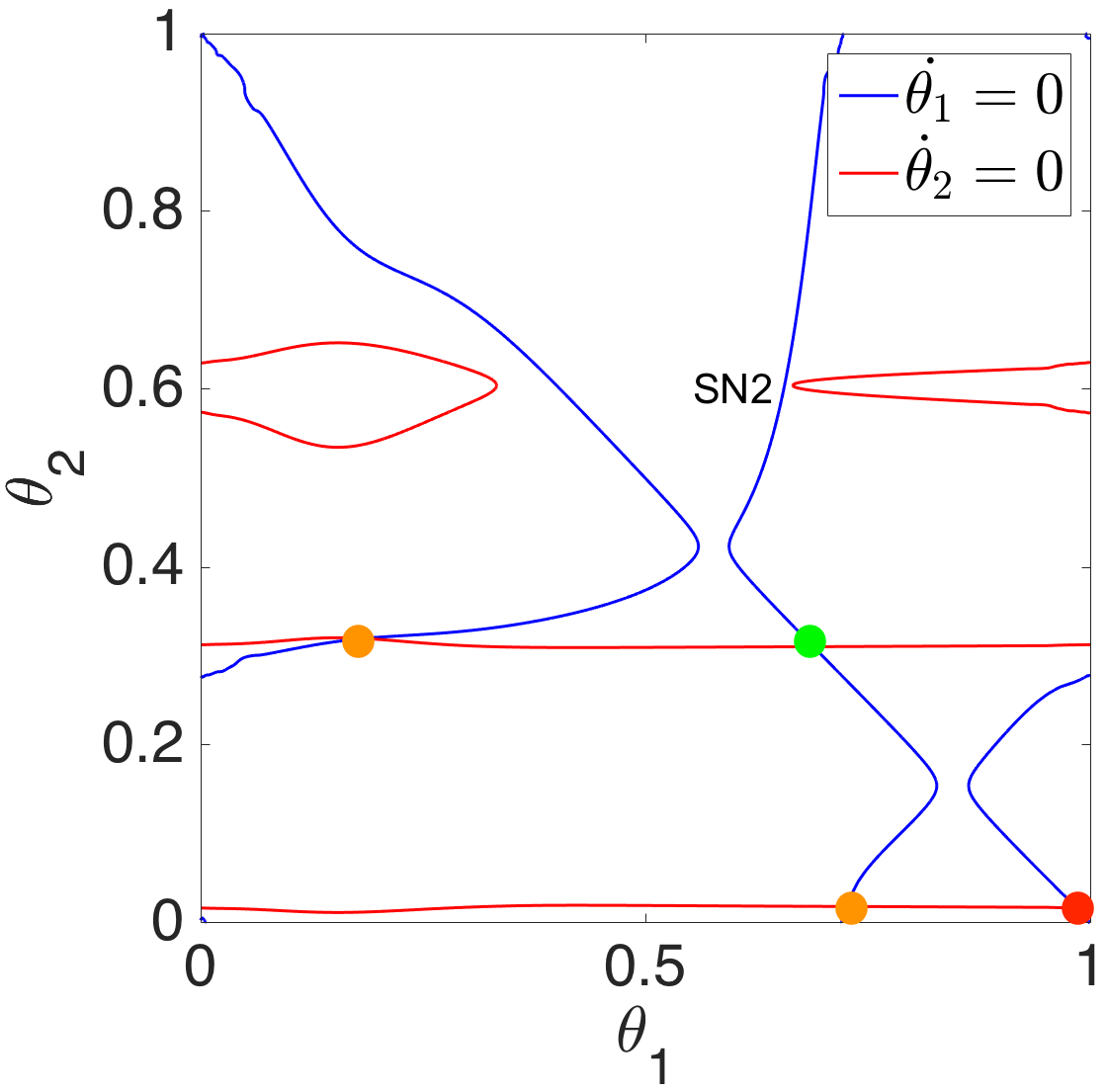}
\end{center}
\caption{(Left to right) Nullclines of 
Equations~(\ref{eq.osc.simplified})  with $\alpha\approx 0.03$ , 
$I_{ext} = 35.65$, and $\delta I_f \approx 0, 0.029, 0.03, 0.032$, respectively. 
At $\delta I_f=0$, there exist one stable forward tetrapod gait, one unstable 
(saddle) backward tetrapod gait, and one unstable tripod gait. One stable
fixed point on $\theta_1=\theta_2$ and several saddle points exist. 
At $\delta I_f \approx 0.03$, through a saddle-node bifurcation (shown by SN1), the unstable tripod gait (shown by N1) 
 and the unstable backward tetrapod gait (shown by S1) disappear. 
 At $\delta I_f \approx 0.032$, through another saddle-node bifurcation (shown by SN2),
  a stable fixed point (shown by N2) and a saddle point (shown by S2) disappear.
 A single stable approximate forward tetrapod gait remains,  together with a source and two saddle points.  
Nullclines and fixed points are indicated as in Figure~\ref{NC_hetero_random_bifurcations}.
 Note that the second figure only shows how the nullclines move and it is topologically equivalent to the first figure.}
\label{NC_hetero_alpha_near0_bifurcations}
\end{figure}

We now consider Equations~(\ref{eq.osc.simplified:forward})
with $\alpha < \alpha_{min}$.  
 Since the qualitative behavior of the solutions of Equations~(\ref{eq.osc.simplified:forward}) with 
  $0<\alpha < \alpha_{min}$ are all similar, we only show the results for $\alpha\approx 0.03$. 
  As illustrated in Figure \ref{NC_hetero_alpha_near0_bifurcations}, choosing the heterogeneity of 
Equation~(\ref{external_Ii}) maintains the $\dot\theta_1 = 0$ nullclines and,  by combining two nullclines
that encircle the torus,  changes the topology
of the $\dot\theta_2 = 0$ nullclines and thereafter, through two saddle-node bifurcations, reduces the number of fixed points. 
As $\delta I_f$ increases from 0 to $0.03$, one  saddle-node bifurcation occurs in which the unstable 
   tripod gait and the backward tetrapod gait disappear; as  $\delta I_f$ increases further to $0.032$, 
   another saddle-node bifurcation occurs and the stable  fixed point on $\theta_1=\theta_2$ disappears 
   and a unique stable approximate forward tetrapod gait remains at $(0.69,0.31)$. The nullclines at  
   $\delta I_f\approx 0.029$ are shown to illustrate how the nullclines move toward each other and 
   cause the saddle-node bifurcations. 
 
So far, we assumed that the forward tetrapod gait is always stable
 and  chose the control parameters $\Ii $  to get a unique stable approximate forward tetrapod gait. 
In the following section, we assume that the backward tetrapod gait is always stable 
and show how to choose $\Ii$ to get a unique stable approximate backward tetrapod gait. 
As discussed earlier, when $\a_{min} < \a$, the backward tetrapod gait is always stable. 

\subsection {Emergence of a unique backward tetrapod gait  at low speed}
\label{backward}
We assume $\a_{min} < \a$ so that the backward tetrapod gait, $(1/3,2/3)$, 
is stable while the forward tetrapod gait can be either stable or a saddle.  

For any $t$, let
 \be{external_Ii_backward}
 \bal
\delta I_b \;:=\; I_{ext}^2(t) = I_{ext}^3(t) \;\geq\; 0, \quad I_{ext}^1(t) = 0, 
 \eal
 \ee 
 and consider  $\delta I_b$ as a bifurcation parameter. 
 
 Choosing $\Ii$ as in Equation~(\ref{external_Ii_backward}) implies $\tilde\omega_3 -  \tilde\omega_2  = 0$, and  
$ \tilde\omega_1 -  \tilde\omega_2 = -\delta I_b \bar Z\leq0 $. 
Therefore, Equations~(\ref{eq.osc.simplified}) become
\begin{subequations} \label{eq.osc.simplified:backward}
\begin{align}
&\dot\theta_1=  - \alpha\delta I_b \bar Z + H(-\theta_1; \x) - \alpha  H(\theta_1; \x)- (1-\alpha)  H(\theta_2; \x), \\ 
&\dot\theta_2 =  H(-\theta_2; \x) -  \alpha H(\theta_1; \x) - (1-\alpha) H(\theta_2; \x). 
\end{align}
\end{subequations}

We will show that when $\alpha_{min}<\alpha<\alpha_{max}$ (resp. $\alpha > \alpha_{max}$), 
as the bifurcation parameter $\delta I_b$ increases, 
Equations~\eqref{eq.osc.simplified:backward} lose 6 (resp. 4) fixed points 
through 3 (resp. 2) saddle-node bifurcations and keep only one stable approximate backward tetrapod gait. 
To show this, we consider two topologically different cases:

\begin{description}[leftmargin=*]

\item[\bm{$\alpha_{min}<\alpha<\alpha_{max}$}] At $\delta I_b = 0$, 
Equations~\eqref{eq.osc.simplified:backward} admit 10 fixed points (5 saddle points, 3 sinks, and 2 sources). 
As $\delta I_b$ increases, 3 saddle-node bifurcations occur and one sink (corresponding to the approximate backward tetrapod gait), 
one source and 2 saddle points remain.  

\item[\bm{$\alpha > \alpha_{max}$}] At $\delta I_b = 0$, 
Equations~\eqref{eq.osc.simplified:backward} admit 8 fixed points (4 saddle points, 2 sinks, and 2 sources). 
As $\delta I_b$ increases, 2 saddle-node bifurcations occur and one sink (corresponding to the approximate backward tetrapod gait), 
one source and 2 saddle points remain.  

\end{description}

\begin{figure}[h!]
\begin{center}
\includegraphics[scale=.1]{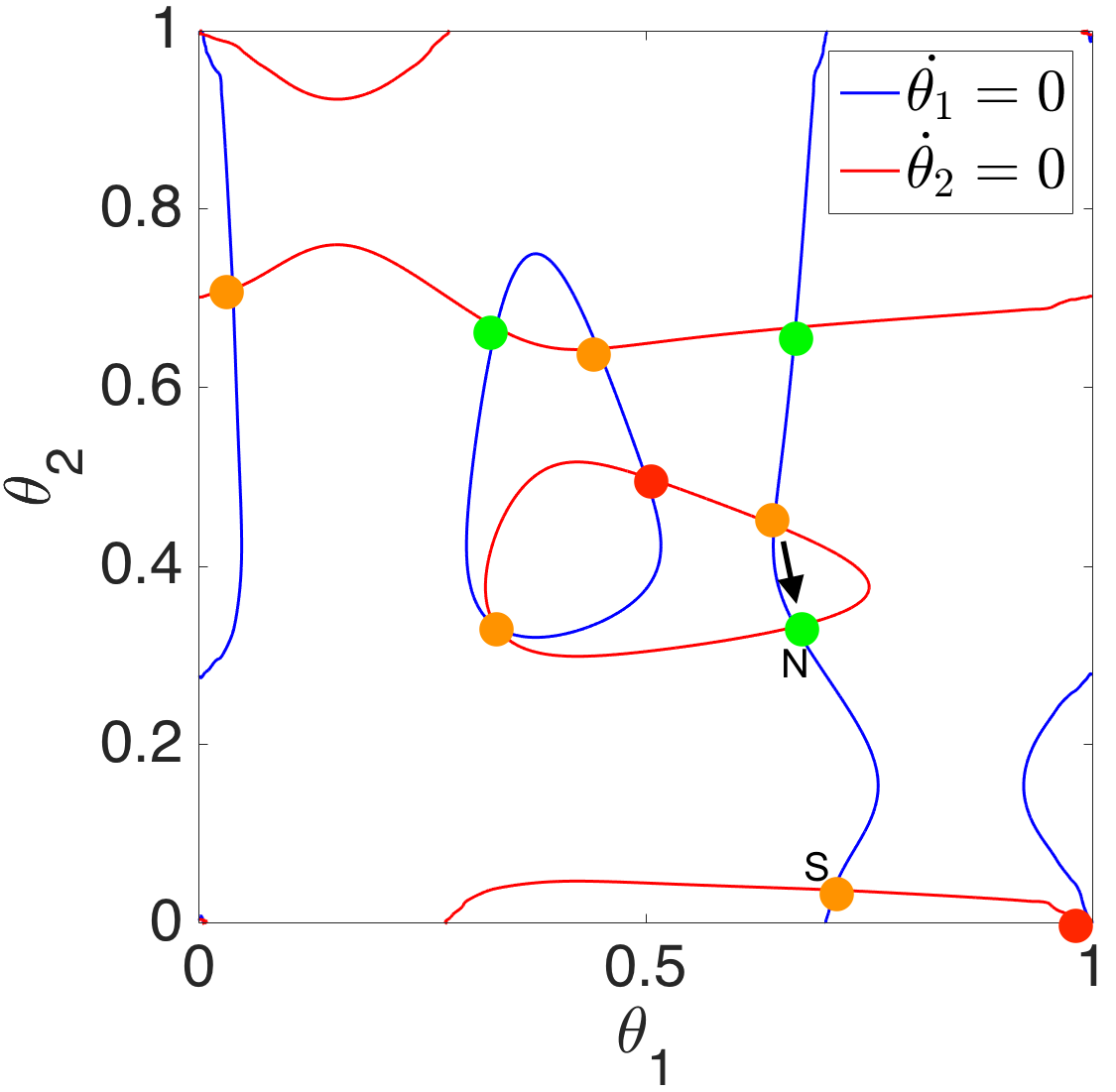}\;
\includegraphics[scale=.1]{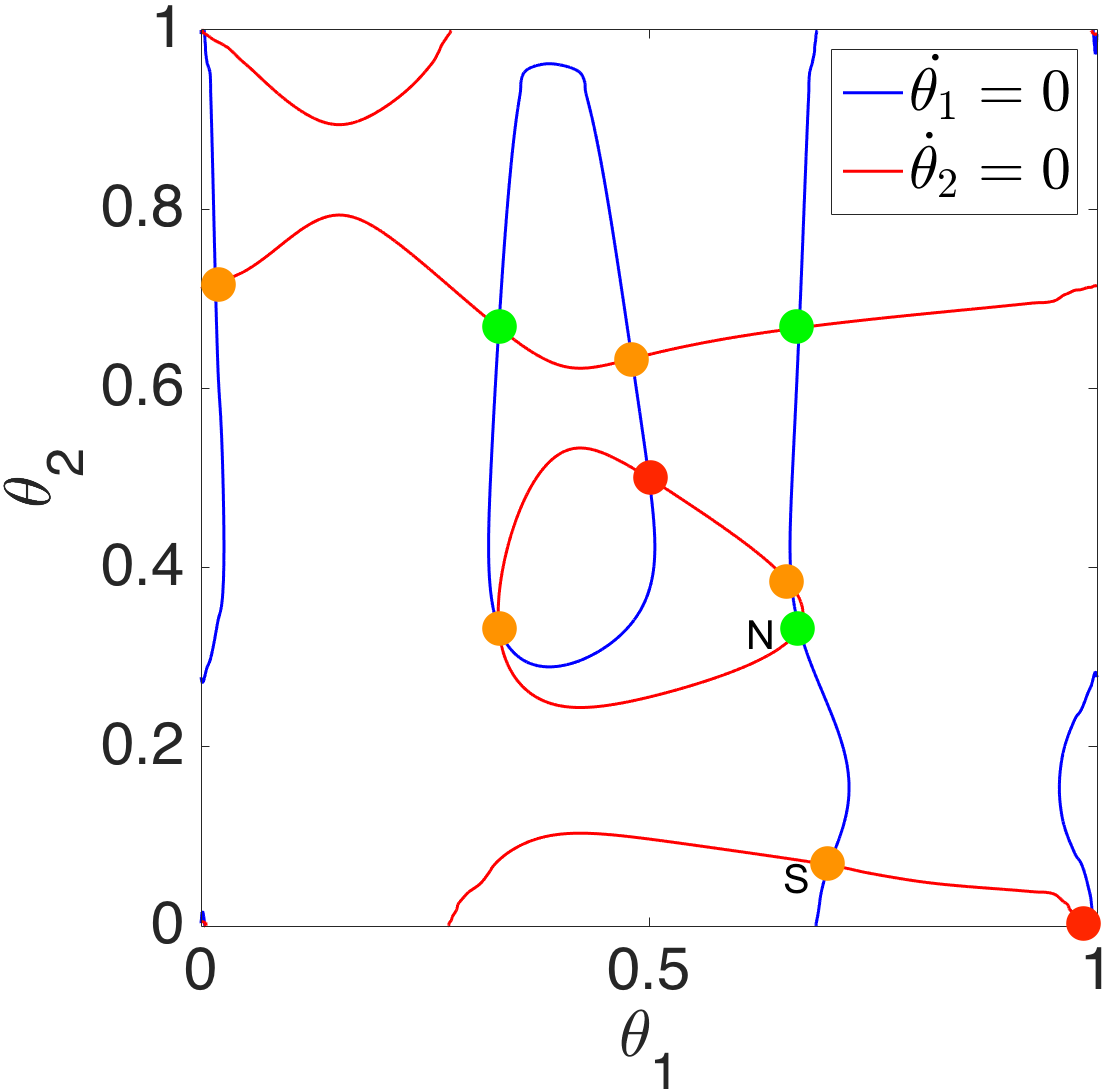}\;
\includegraphics[scale=.1]{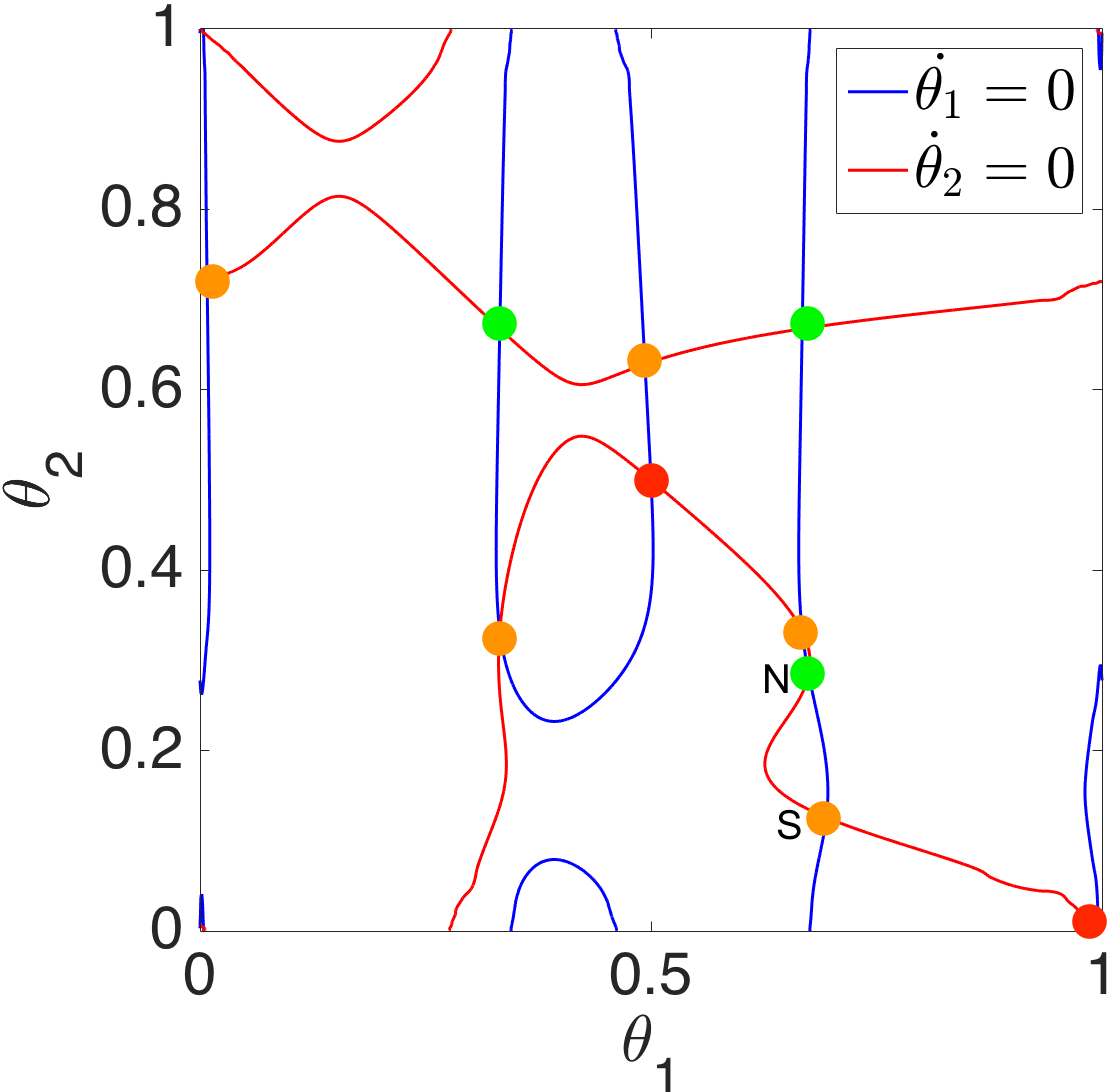}\;
\includegraphics[scale=.1]{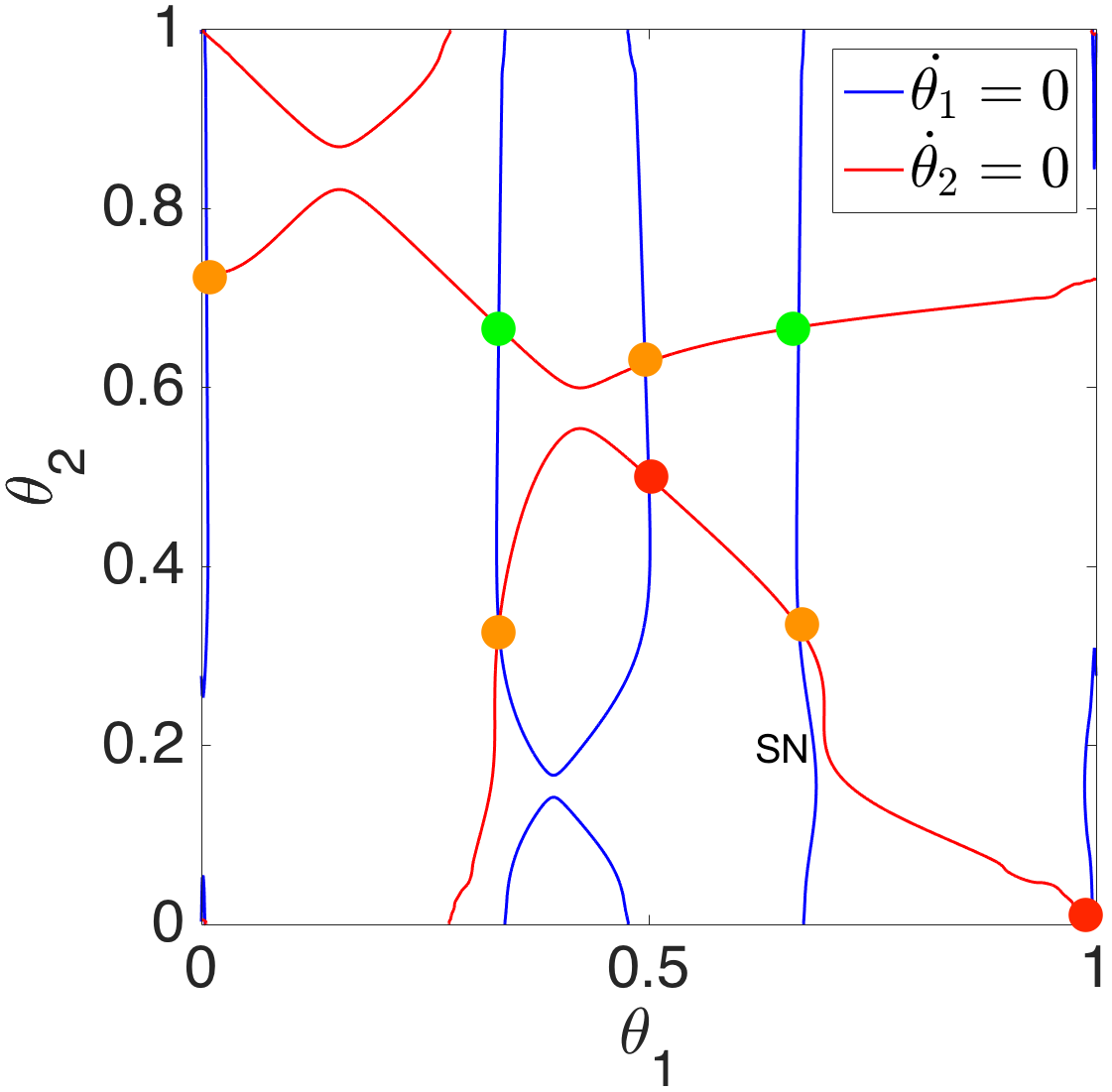}
\end{center}
\caption{(Left to right) Nullclines of 
Equations~(\ref{eq.osc.simplified:backward})  with $\alpha = 1/2, 1/1.3, 1/1.1,1/1.06$, 
$I_{ext} = 35.65$, and $\delta I_b =0$, respectively. 
At $\alpha = 1/2$, there exist 10 fixed points including one stable forward tetrapod gait shown by N. 
As $\alpha$ increases, through a saddle-node bifurcation  (shown by SN)  the stable forward tetrapod gait disappears together with 
a saddle point shown by S and another saddle point becomes 
 an unstable forward tetrapod gait. The stable backward tetrapod remains throughout.
Nullclines and fixed points are indicated as in Figure~\ref{NC_hetero_random_bifurcations}. }
\label{saddle_node_alpha_backward}
\end{figure}

 Figure~\ref{saddle_node_alpha_backward} shows that as $\a$ increases and approaches $\a_{max}$, 
 through a saddle-node bifurcation, the stable forward tetrapod gait disappears 
 and one of the  saddle points moves toward the position of the  forward tetrapod gait, as shown by an arrow in 
 Figure~\ref{saddle_node_alpha_backward} (left). 
 As $\alpha$ increases,  the isolated $\dot\theta_2=0$ nullcine combines with the $\dot\theta_2=0$ nullcline
that encircles the torus and thereafter, at $\a=\a_{max}$,  the number of fixed points reduces to 8 from 10.
 Therefore, when $\a>\a_{max}$, there are only 8 fixed points. 

\subsubsection*{3 saddle-node bifurcations: \bm{$\alpha_{min}<\alpha<\alpha_{max}$}}\label{backward_3SN}

\begin{figure}[h!]
\begin{center}
\includegraphics[scale=.1]{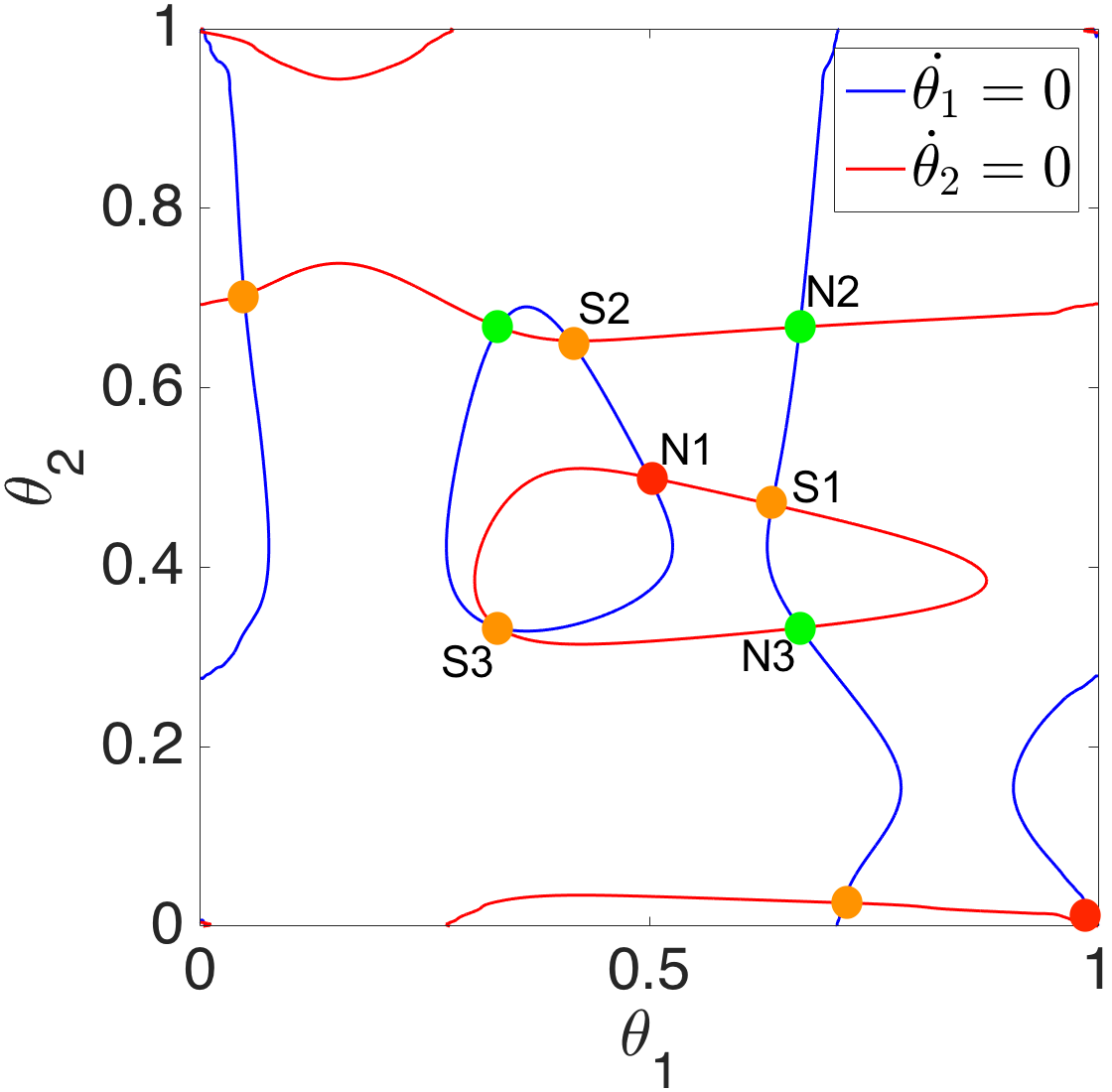}\;
\includegraphics[scale=.1]{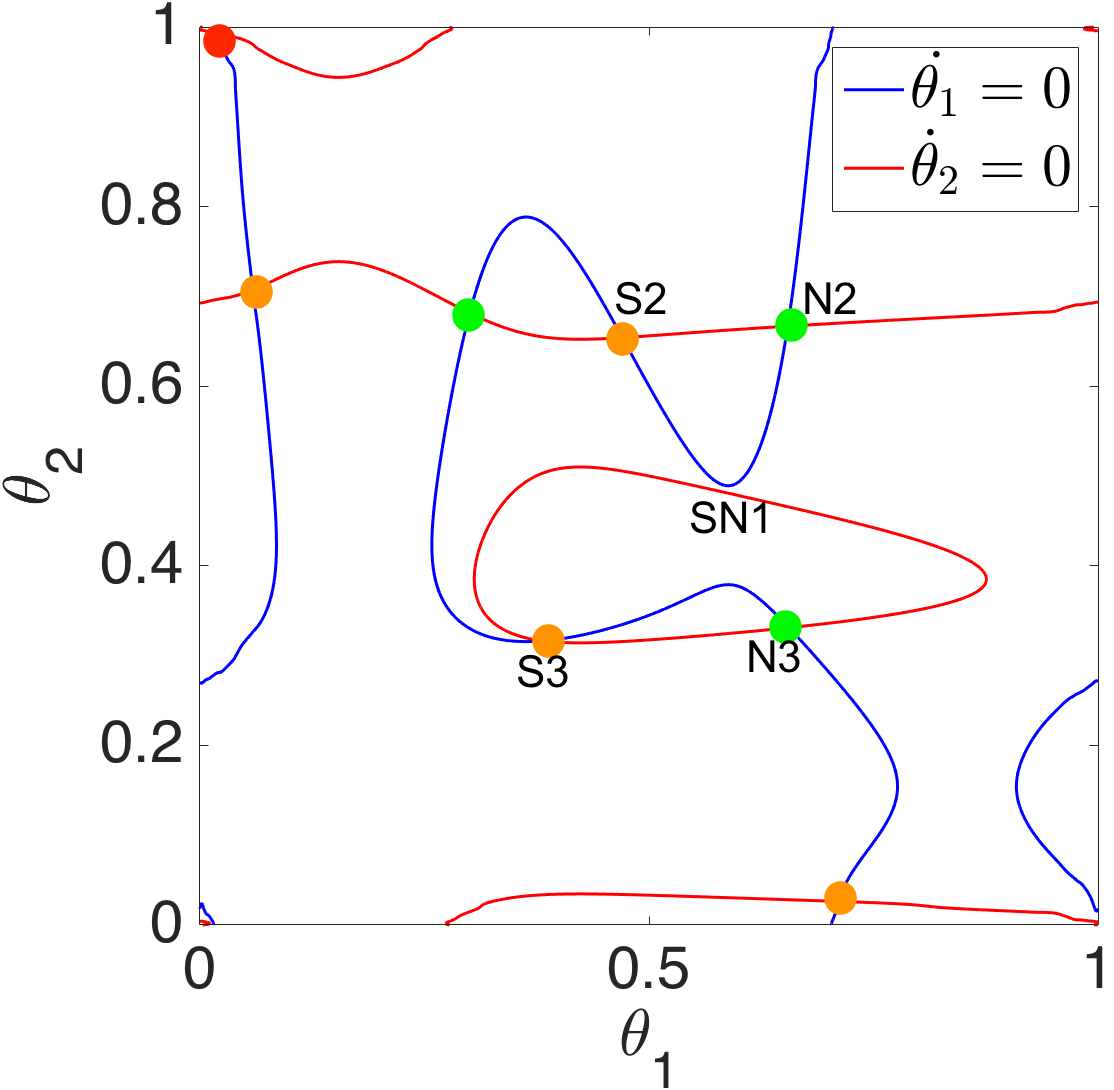}\;
\includegraphics[scale=.1]{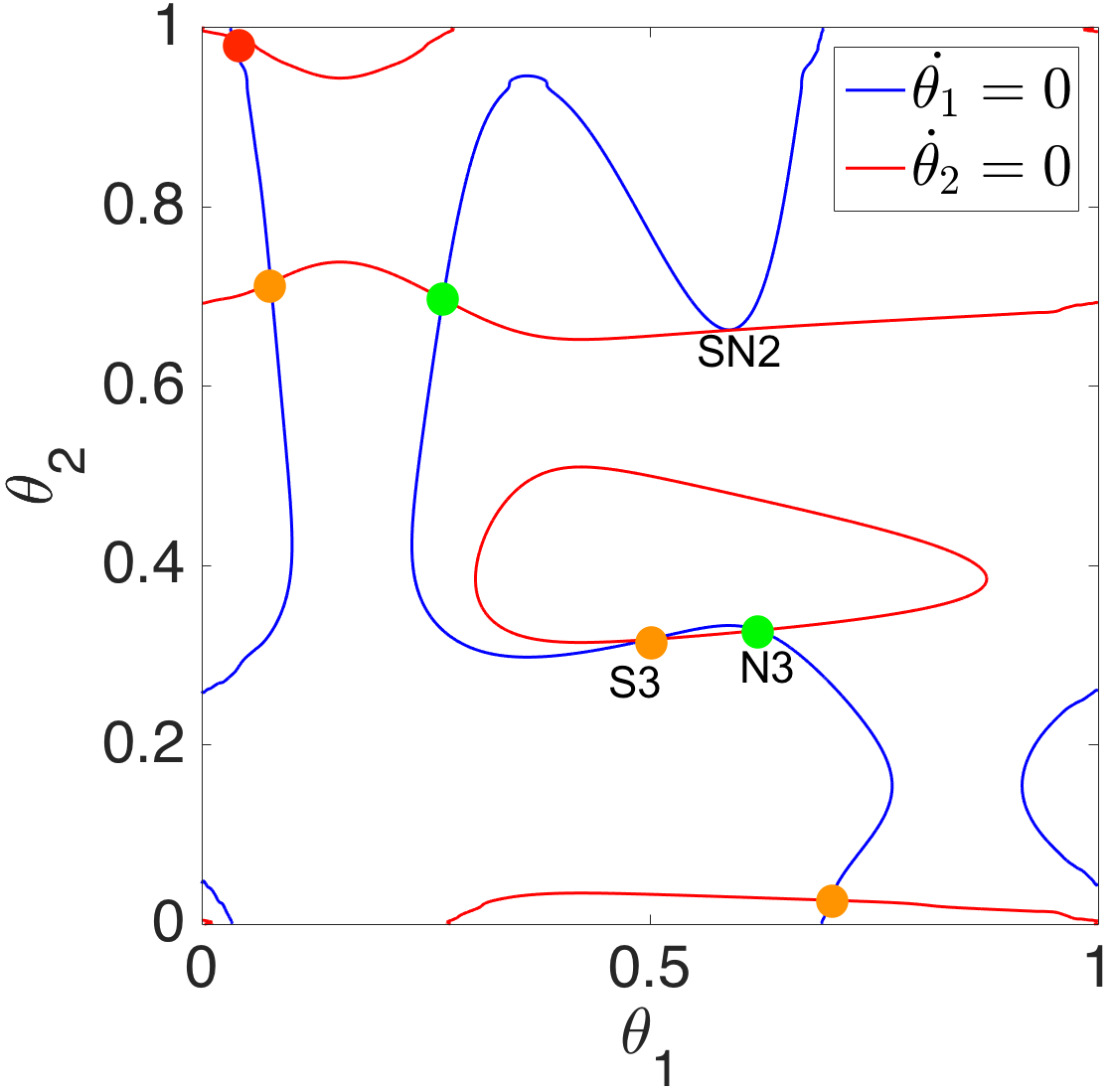}\;
\includegraphics[scale=.1]{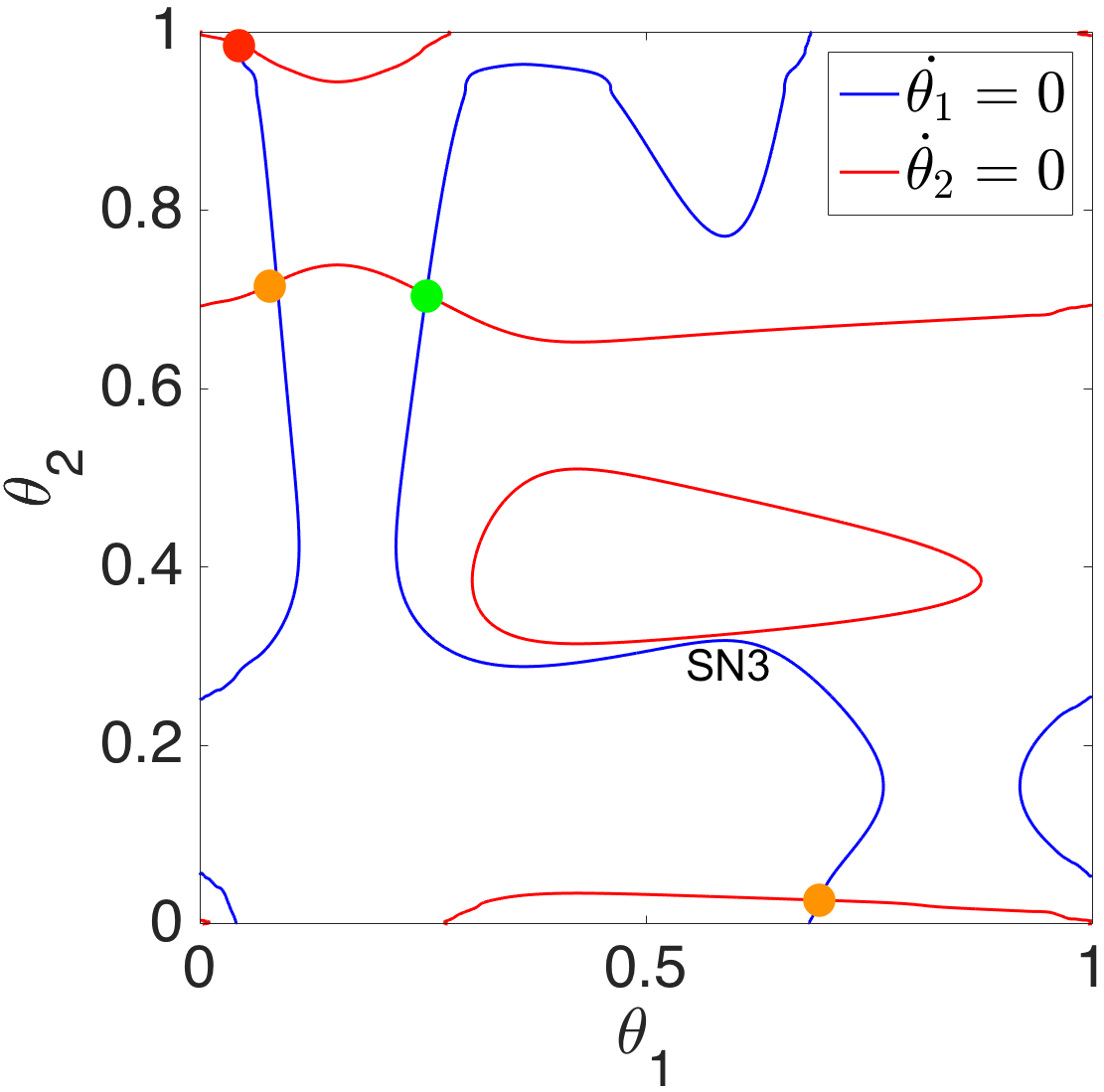}
\end{center}
\caption{(Left to right) Nullclines of 
Equations~(\ref{eq.osc.simplified:backward}) with $\alpha=1/3$, 
$I_{ext} = 35.65$, and $\delta I_b \approx 0, 0.015 , 0.04 , 0.056$, respectively.
At $\delta I_b =0$, there exist ten fixed points. 
At $\delta I_b \approx 0.015$, the first saddle-node bifurcation (shown by SN1) occurs 
and the unstable tripod gait (shown by N1) disappears together with a saddle point (shown by S1).
At $\delta I_b \approx 0.04$, the second saddle-node bifurcation (shown by SN2) occurs and a stable 
fixed point (shown by N2) disappears together with a saddle point (shown by S2).
 Finally, at $\delta I_b \approx 0.056$, the third saddle-node bifurcation (shown by SN3) occurs 
 and the stable forward tetrapod gait (shown by N3) disappears together with a saddle point  (shown by S3).
A single stable approximate backward tetrapod gait remains,  together with a source and two saddle points.  
Nullclines and fixed points are indicated as in Figure~\ref{NC_hetero_random_bifurcations}.} 
\label{NC_hetero_backward}
\end{figure}

Consider Equations~(\ref{eq.osc.simplified:backward}) with $\alpha_{min}<\alpha<\alpha_{max}$. 
Since the qualitative behavior of the solutions of Equations~(\ref{eq.osc.simplified:backward})  
with $\alpha_{min}<\alpha<\alpha_{max}$ are all similar, we show the results in an example with $\a=1/3$. 
As is clear from Equations~(\ref{eq.osc.simplified:backward}) and illustrated in
Figure \ref{NC_hetero_backward}, choosing the heterogeneity of 
Equation~(\ref{external_Ii_backward}), maintains the $\dot\theta_2 = 0$ nullclines and only
perturbs the $\dot\theta_1 = 0$ nullclines. This perturbation causes the topology
of the $\dot\theta_1 = 0$ nullclines to change, combining the isolated circle with a nullcline
that encircles the torus and thereafter reducing the number of fixed points. 

In Figure~\ref{NC_hetero_backward}, we show the nullclines of
Equations~(\ref{eq.osc.simplified:backward}) with $\alpha=1/3$, 
and increase $\delta I_b$ from 0 to $0.015$, where the first saddle-node bifurcation
occurs and the unstable tripod gait disappears.
We further increase $\delta I_b$ to $0.04$ where the second saddle-node bifurcation
occurs and the stable $\approx(2/3,2/3)$
fixed point disappears. Finally, when $\delta I_b$ reaches $0.056$, the third
saddle-node bifurcation occurs and the stable forward tetrapod
gait disappears and only one stable fixed point  remains, which corresponds to
the approximate backward tetrapod gait $\approx (0.25,0.7)$,  as we desired. 
  
\subsubsection*{2 saddle-node bifurcations: \bm{$\alpha > \alpha_{max}$}}\label{backward_2SN}

\begin{figure}[h!]
\begin{center}
\includegraphics[scale=.1]{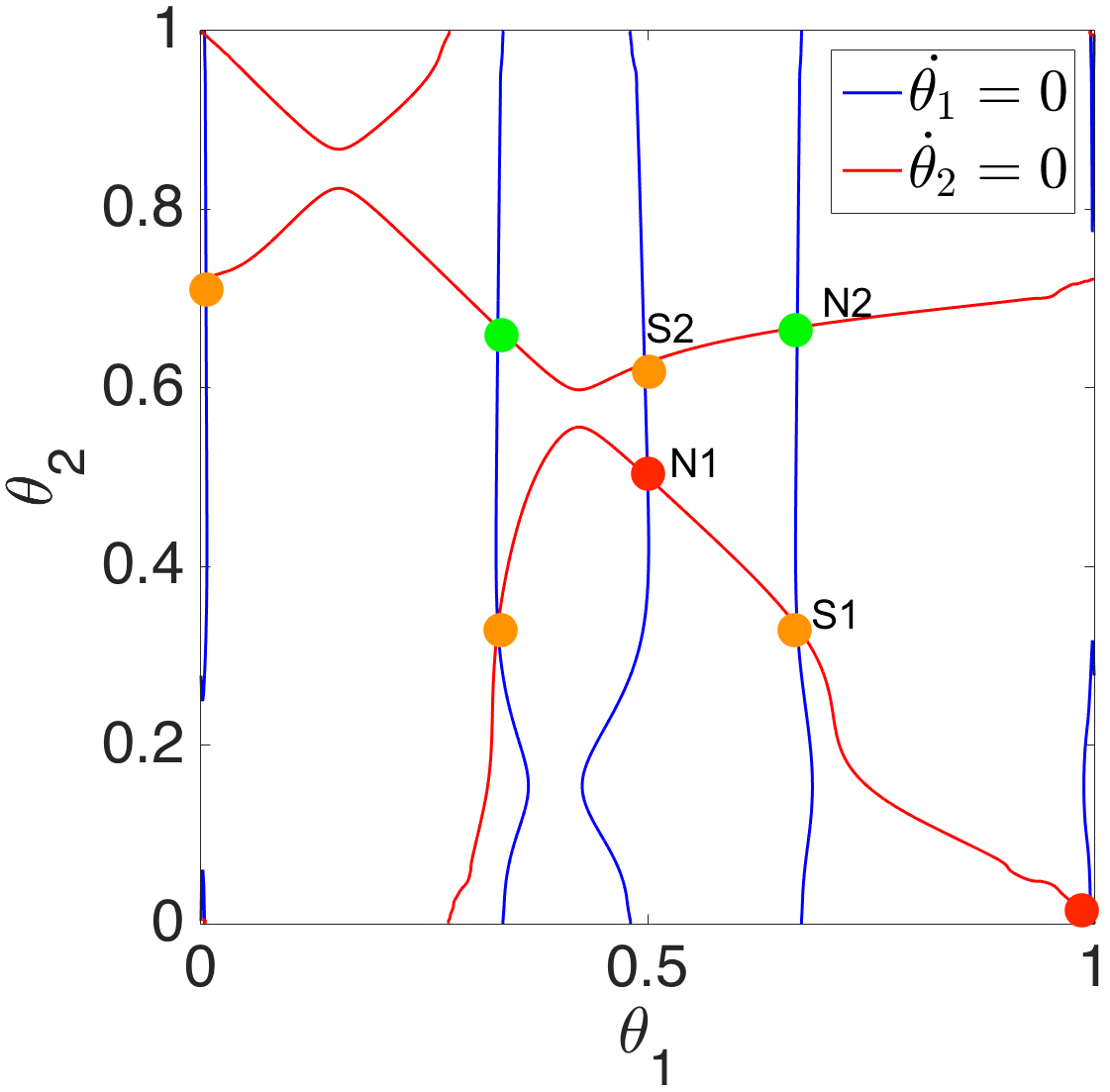}\;
\includegraphics[scale=.1]{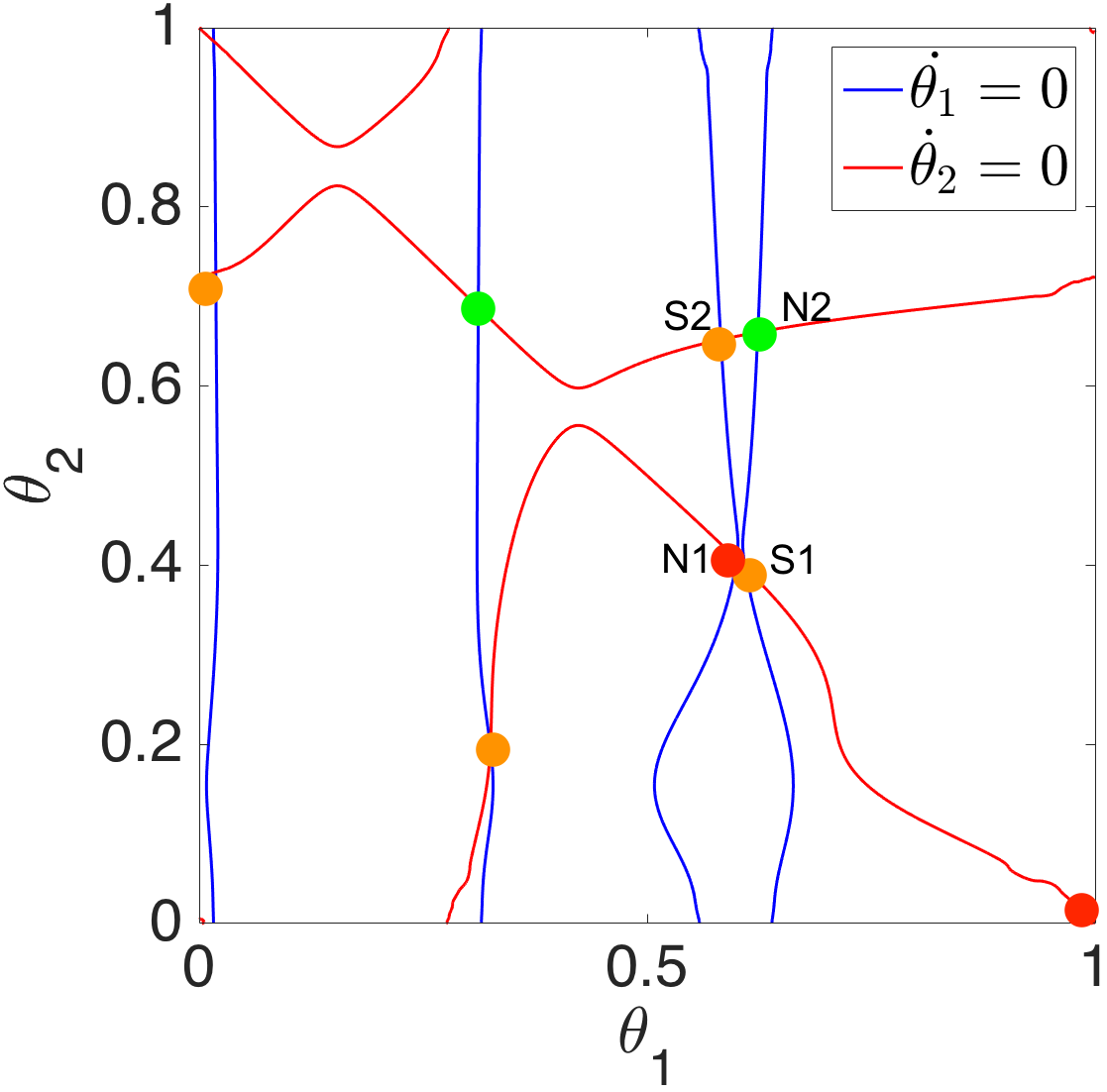}\;
\includegraphics[scale=.1]{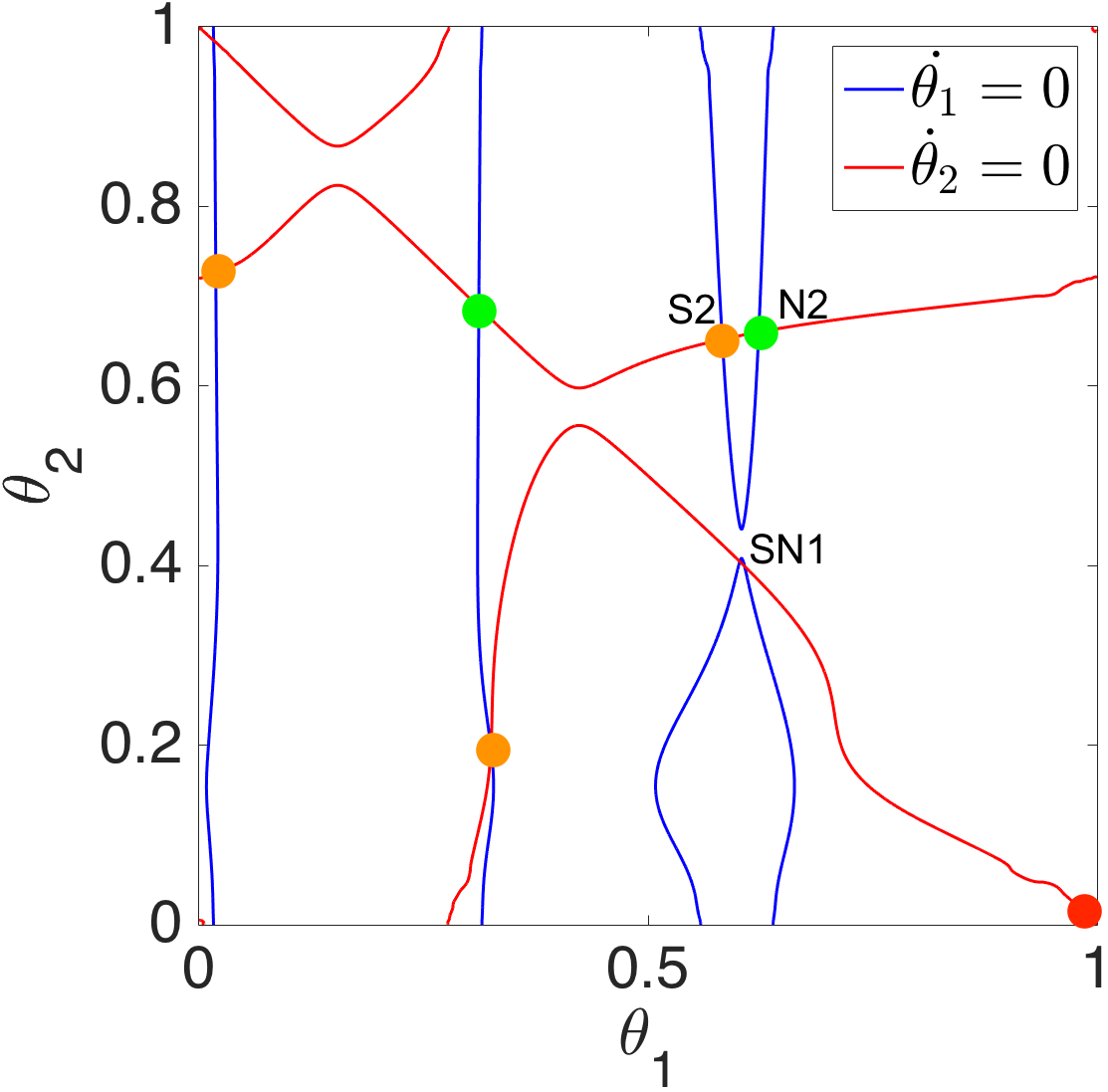}\;
\includegraphics[scale=.1]{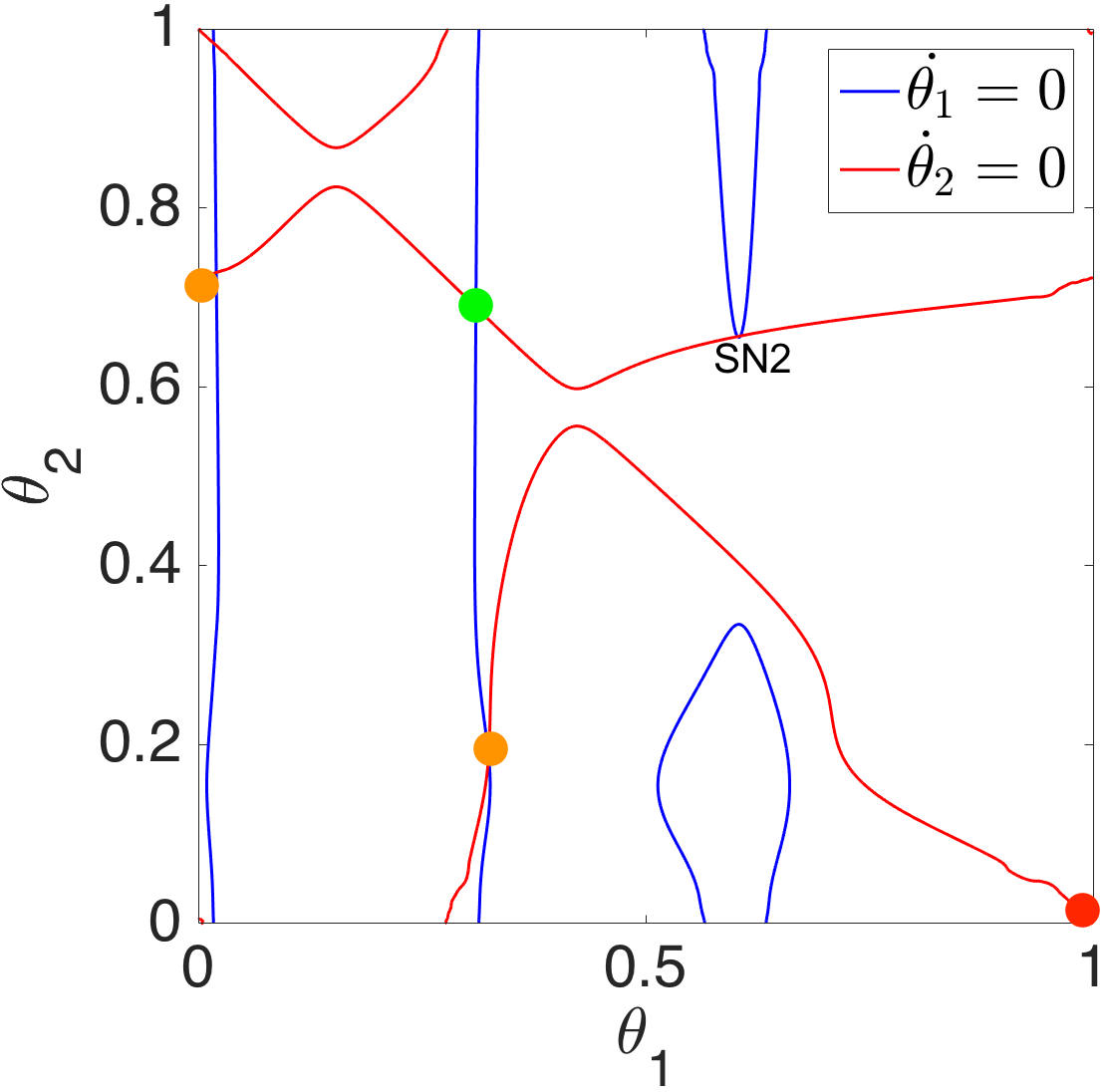}
\end{center}
\caption{(Left to right) Nullclines of 
Equations~(\ref{eq.osc.simplified:backward})  with $\alpha\approx 0.95$ , 
$I_{ext} = 35.65$, and $\delta I_b \approx 0, 0.012, 0.0121, 0.013$, respectively. 
At $\delta I_b=0$, there exist one stable backward tetrapod gait, one unstable 
(saddle) forward tetrapod gait, and one unstable tripod gait. One stable
fixed point on $\theta_1=\theta_2$ and four saddle points exist. 
At $\delta I_b\approx 0.0121$, through a saddle-node bifurcation (shown by SN1),
the unstable tripod gait (shown by N1) 
and the unstable forward tetrapod gait (shown by S1) disappear. 
At $\delta I_b\approx 0.013$, through another saddle-node bifurcation (shown by SN2),
a stable fixed point (shown by N2) and a saddle point (shown by S2) disappear.
A single stable approximate backward tetrapod gait remains together
with a source and two saddle points.  
Nullclines and fixed points are indicated as in Figure~\ref{NC_hetero_random_bifurcations}.
Note that the second panel only shows how the $\dot{\theta}_1=0$ nullclines move
and it is topologically equivalent to the first panel. }
\label{NC_hetero_alpha_near1_bifurcations}
\end{figure}

We now consider Equations~(\ref{eq.osc.simplified:backward})
with $\alpha > \alpha_{max}$. 
 Since the qualitative behavior of the solutions of Equations~(\ref{eq.osc.simplified:backward}) with 
  $\alpha > \alpha_{max}$ are all similar, we only show the results for $\alpha\approx0.95 $. 
 As is  illustrated in Figure \ref{NC_hetero_alpha_near1_bifurcations}, choosing the heterogeneity of 
Equation~(\ref{external_Ii_backward}) maintains the $\dot\theta_2 = 0$ nullclines and,   by combining two $\dot\theta_1=0$ nullclines
that encircle the torus, changes the topology
of the $\dot\theta_1=0$ nullclines and thereafter, through two saddle-node bifurcations,  reduces the number of fixed points. 
As $\delta I_b$ increases from 0 to $0.0121$, one  saddle-node bifurcation occurs in which the unstable 
   tripod gait and the backward tetrapod gait disappear; as  $\delta I_b$ increases further to $0.013$, 
   another saddle-node bifurcation occurs and the stable  fixed point shown by N2 disappears 
   and a unique stable approximate backward tetrapod gait remains at $(0.31,0.69)$. The nullclines at  
   $\delta I_b\approx 0.012$ are shown to illustrate how the $\dot\theta_1=0$ nullclines move toward each other 
   and cause the saddle-node bifurcations. 

\bremark 

In Table \ref{summary:figs}, we summarize the main results shown in 
Figures \ref{saddle_node_alpha_forward}--\ref{NC_hetero_alpha_near1_bifurcations}. 

\begin{table}[ht]
	{\scriptsize{
			\begin{center}
				\begin{tabular}{|c| c c c c c |}
					\hline
					Figure \# & $\a$ & $I^1_{ext}$ & $I^2_{ext}$ & $I^3_{ext}$ & bifurcation type \\
					\hline
                                           5 & varies & 0 & 0 & 0 & 1 saddle-node\\ 
					\hline
                                           6 & $\a_{min} < \a <\a_{max}$ & ($=\delta I_f$) varies & ($=\delta I_f$) varies & 0& 3 saddle-node \\
                                           \hline
                                           7 & $0<\a<\a_{min}$ &($=\delta I_f$) varies &($=\delta I_f$) varies & 0 & 2 saddle-node \\
					\hline
                                           8 & varies & 0 & 0 & 0 & 1 saddle-node \\ 
					\hline
                                           9 & $\a_{min} < \a <\a_{max}$ & 0 &($=\delta I_b$) varies & ($=\delta I_b$) varies & 3 saddle-node \\
                                           \hline
                                           10 & $\a_{max}<\a$ & 0 &($=\delta I_b$) varies & ($=\delta I_b$) varies & 2 saddle-node  \\
					\hline
				\end{tabular}
			\end{center}
			\caption{A summary of Figures \ref{saddle_node_alpha_forward}--\ref{NC_hetero_alpha_near1_bifurcations}. }
			\label{summary:figs}
		}}
	\end{table}%

\eremark 

\subsection{Transition from the approximate tetrapod to the approximate tripod gait}
\label{section:transition}

 In \cite{SIADS2018} we studied gait transition from 
multiple tetrapod gaits  (e.g., Figure~\ref{NC_hetero_random_bifurcations}(left)) 
to a unique tripod gait, as speed increases. 
Here, we introduced approximate transition gaits $(1/3+\et, 2/3-\et)$ (resp. $(2/3-\et, 1/3+\et)$) and discussed that 
for suitable heterogeneous systems, as $\et$ changes from 0 to $1/6$, 
 they connect a single stable approximate forward (resp. backward) tetrapod gait to a single stable approximate tripod gait. 
 
\begin{figure}[h!]
\begin{center}
\includegraphics[scale=.1]{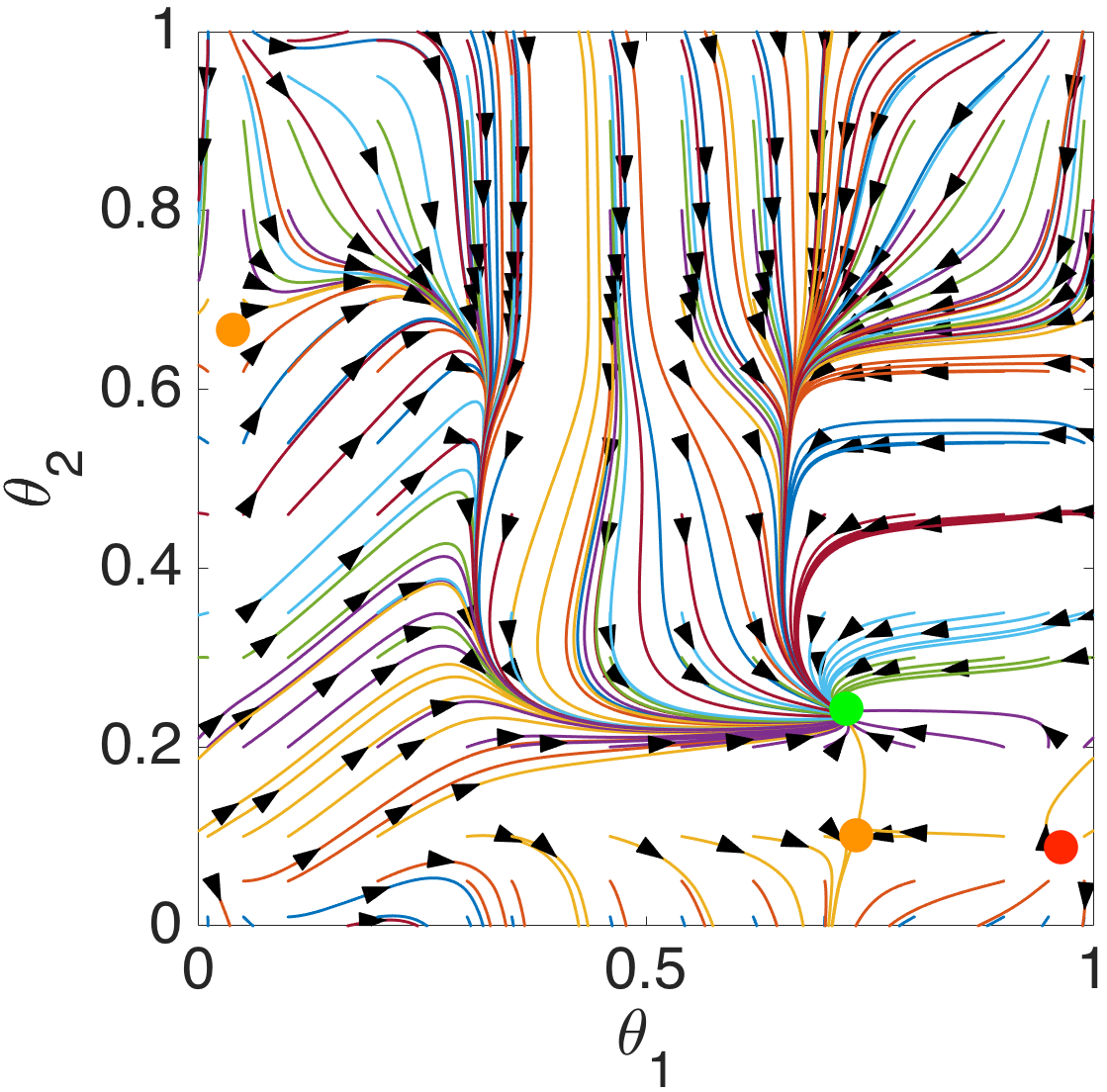}\;
\includegraphics[scale=.1]{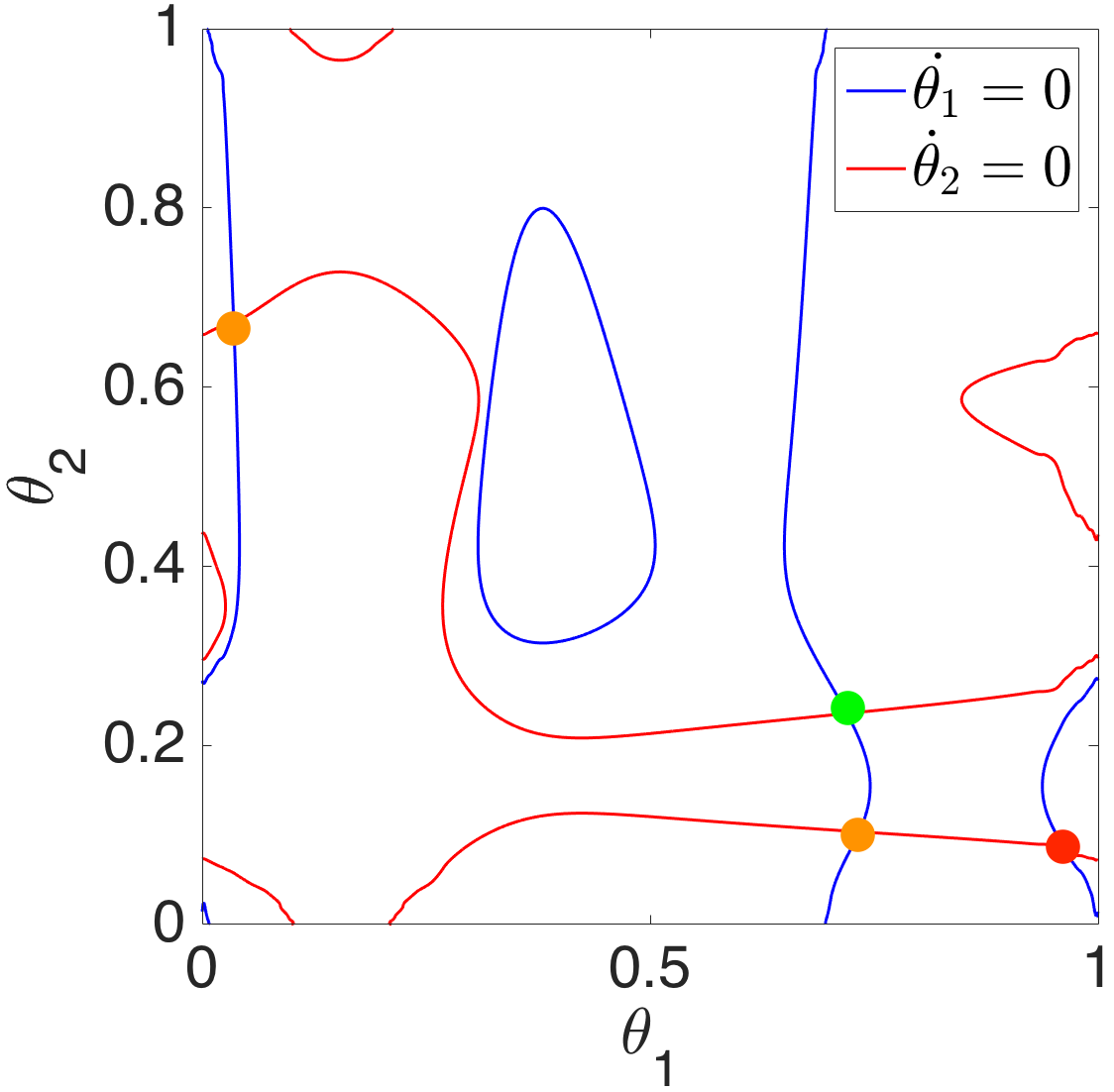}\;
\includegraphics[scale=.1]{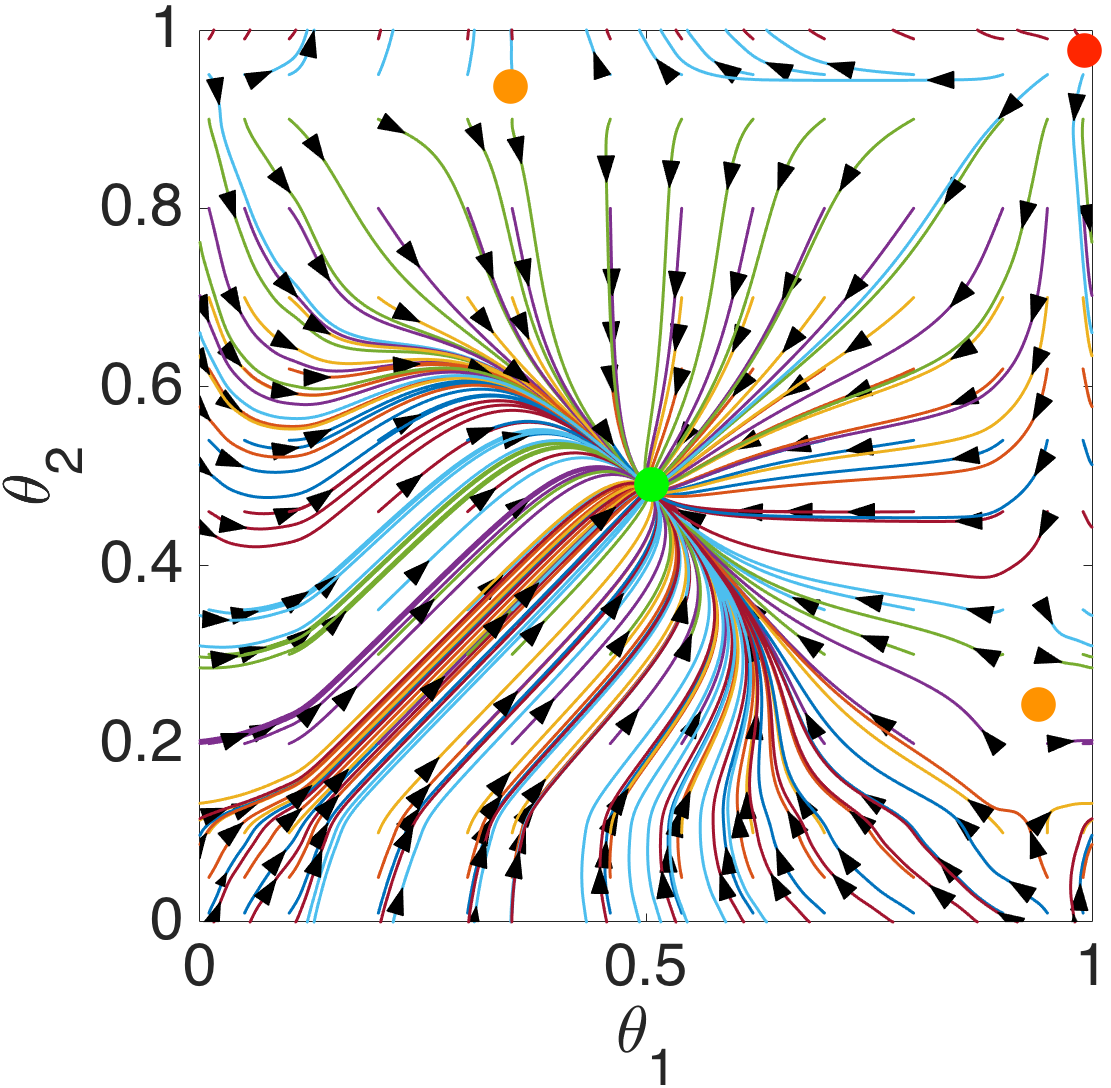}\;
\includegraphics[scale=.1]{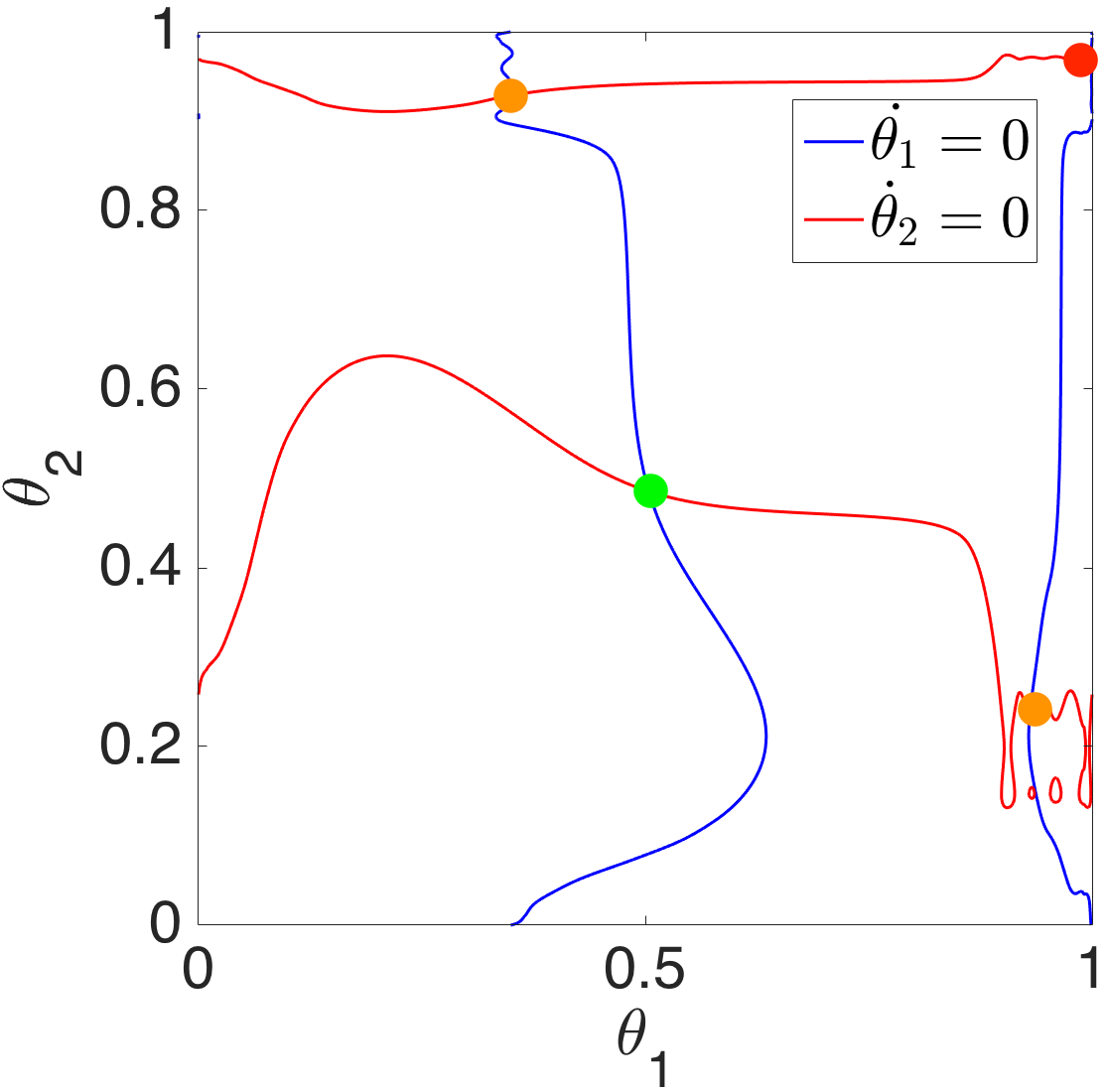}
\end{center}
\caption{(Left to right) Phase planes and nullclines of  Equations~(\ref{torus:equation})
when $\c_i$'s satisfy Equation~(\ref{random_balance_coupling}), 
$\delta I \approx 0.038$, and $I_{ext} = 35.65$ (the nullclines of which are 
also shown in Figure~\ref{NC_hetero_random_bifurcations} (right panel)), and $37.5$, 
respectively. Nullclines and fixed points are indicated as in 
Figure~\ref{NC_hetero_random_bifurcations}. } 
\label{NC_hetero_random_transition}
\end{figure}

As an illustration, we show the gait transition in Equations~(\ref{torus:equation})
when the $\c_i$'s satisfy Equation~(\ref{random_balance_coupling}), $\delta I \approx 0.038$, 
and $I_{ext}$ increases from $35.65$ to $37.5$. 
In Figure~\ref{NC_hetero_random_transition}, we observe that as $I_{ext}$ increases,  
the unique stable approximate forward tetrapod gait 
becomes a stable approximate tripod gait. 
The second and fourth figures show the nullclines and hence the positions 
of the fixed points for $I_{ext} = 35.65, 37.5$, respectively; and the first and 
third figures show the corresponding phase planes.
 
\section{Equivalent perturbations}
\label{Equivalent_Perturbations}

In this section, we will show that perturbing the intrinsic dynamics of each
unit of the CPG can be equivalent to 
perturbing the coupling function $g$ or the coupling strengths $\c_i$. 

Recalling Equation~(\ref{24ODE_closed}), we show that, under 
an appropriate condition on the $\c_i$'s, derived below,
adding $I^i_{ext}$ to each neuron $i$ is equivalent to  adding $dI_j$ 
to the coupling function $g(x_i,x_j)$ that 
connects neuron $i$ to its neighbor $j$, where  $dI_j$ is the unique solution of
\be{equi:pert:g}
\left(\begin{array}{c}I_{ext}^1 \\ \vdots \\I_{ext}^6\end{array}\right) = 
\mathcal{C}\left(\begin{array}{c}dI_1 \\dI_2 \\dI_3 \\dI_4 \\dI_5 \\dI_6\end{array}\right) := 
\left(\begin{array}{cccccc}
0 & \c_5 & 0 & \c_1 & 0 & 0 \\
\c_4 & 0 & \c_7 & 0 & \c_2 & 0 \\
0 & \c_6 & 0 & 0 & 0 & \c_3 \\
\c_1 & 0 & 0 & 0 & \c_5 & 0 \\
0 & \c_2 & 0 & \c_4 & 0 & \c_7 \\
0 & 0 & \c_3 & 0 & \c_6 & 0
\end{array}\right)
\left(\begin{array}{c}dI_1 \\dI_2 \\dI_3 \\dI_4 \\dI_5 \\dI_6\end{array}\right).  
\ee

For example
if the above equation has a unique solution, 
then since  $I_{ext}^1 = \c_5 dI_2 + \c_1 dI_4$,  
 adding $I_{ext}^1$ to unit $1$  is equivalent to adding $dI_2$ to $g(x_1, x_2)$ and  $dI_4$ to $g(x_1, x_4)$, 
 i.e., 
 \beqn
 \dot{x}_1  &=&  f(x_1) +I_{ext}^1 + \c_1g(x_1,x_4) + \c_5g(x_1,x_2)\\
                 &=&  f(x_1) +\c_5 dI_2 + \c_1dI_4+ \c_1g(x_1,x_4) + \c_5g(x_1,x_2)\\
                 &=&  f(x_1) +  \c_1(g(x_1,x_4) +dI_4)+ \c_5(g(x_1,x_2)+dI_2). 
\eeqn

Equation~(\refeq{equi:pert:g}) has a unique solution if the matrix $\mathcal{C}$ is non-singular,
i.e., $\det\mathcal{C} \neq0$. The matrix $\mathcal{C}$ can be written as 
\[
\mathcal{C} = \left(\begin{array}{cc}\mathcal{A} & \mathcal{B} \\\mathcal{B} & \mathcal{A}\end{array}\right),
\]
 where 
 $\mathcal{A} =
 \left(\begin{array}{ccc}
0 & \c_5 & 0  \\
\c_4 & 0 & \c_7 \\
0 & \c_6 & 0
\end{array}\right)$
 and $\mathcal{B} =  \diag (\c_1,\c_2,\c_3)$. 
Since
  \be{}
 \bal
  \left(\begin{array}{cc}I &0 \\I &I\end{array}\right)
    \left(\begin{array}{cc}\mathcal{A} & \mathcal{B} \\\mathcal{B} & \mathcal{A}\end{array}\right)
  \left(\begin{array}{cc}I &0 \\-I &I\end{array}\right)
     = 
      \left(\begin{array}{cc}\mathcal{A} - \mathcal{B} & \mathcal{B} \\ 0 & \mathcal{A}+\mathcal{B}\end{array}\right), \nonumber
 \eal
 \ee
 where $I$ and $0$ are identity and zero matrices of appropriate sizes, 
and as shown in  \cite{silvester2000}, 
\[
\det\left(\begin{array}{cc}\mathcal{A} - \mathcal{B} & \mathcal{B} \\ 0 & \mathcal{A}+\mathcal{B}\end{array}\right)
=  \det (\mathcal{A} - \mathcal{B}) \det (\mathcal{A} + \mathcal{B}),
\]
 we have 
   \be{}
 \bal
 \det\mathcal{C} = \det (\mathcal{A} - \mathcal{B}) \det (\mathcal{A} + \mathcal{B})
&=  \det  \left(\begin{array}{ccc}
-\c_1 & \c_5 & 0  \\
\c_4 & -\c_2& \c_7 \\
0 & \c_6 & -\c_3
\end{array}\right)
  \det  \left(\begin{array}{ccc}
\c_1 & \c_5 & 0  \\
\c_4 & \c_2& \c_7 \\
0 & \c_6 & \c_3
\end{array}\right) \\ 
&= - (\c_1\c_2\c_3 - \c_1\c_6\c_7 - \c_3\c_4\c_5)^2. \nonumber
 \eal
 \ee
Hence, $\mathcal{C}$ is non-singular if and only if 
$\c_1\c_2\c_3 - \c_1\c_6\c_7 - \c_3\c_4\c_5 \neq 0. $

Next, we show that, 
perturbing each CPG unit by an external current can be equivalent to
perturbing the coupling strengths. Recalling Equations~(\ref{eq.osc1}) 
and Assumptions~\ref{external_input_assumption} and \ref{constant_contralateral_phase_diff}, 
adding $I^i_{ext}$ to each unit $i$ is equivalent to  adding 
$\Delta_i := \tilde\omega_i/H(2/3-\et; \x)$ to the contralateral coupling $\c_i$, $i=1,2,3$, 
 and keeping the other coupling strengths unchanged. 
  Note that $\tilde\omega_i$ is of order $\epsilon$ and $H$ is of order 1, therefore $\Delta_i$ is of order $\epsilon$. 
 
 For example,  adding $I_{ext}^1$ to unit $1$  is equivalent to adding  $\tilde\omega_1$ to the corresponding phase equation, 
 therefore by Assumption \ref{constant_contralateral_phase_diff}, we get
 \be{}\nonumber
\bal
 \dot{\phi}_1  &=\omega + \tilde\omega_1+ \c_1H(\phi_4 - \phi_1;\x) + \c_5H(\phi_2 - \phi_1;\x)\\
                     &=\omega + \tilde\omega_1+ \c_1H(2/3-\et; \x) + \c_5H(\phi_2 - \phi_1;\x) \\
                     &=\omega + \frac{\tilde\omega_1}{H(2/3-\et; \x)} H(2/3-\et; \x)+ \c_1H(2/3-\et; \x) + \c_5H(\phi_2 - \phi_1;\x) \\
                     &=\omega + \lt(\frac{\tilde\omega_1}{H(2/3-\et; \x)} + \c_1\rt)H(2/3-\et; \x) + \c_5H(\phi_2 - \phi_1;\x) \\
                     &=\omega + (\Delta_1+ \c_1) H(2/3-\et; \x) + \c_5H(\phi_2 - \phi_1;\x). 
\eal
\ee

\section{Discussion}
\label{Conclusion}

In \cite{SIADS2018} we studied a homogeneous 
interconnected phase oscillator model for insect locomotion, 
and showed that the cyclic motion of each leg can be described by 
an oscillator, and that the insect's speed increases with the 
common external input, $I_{ext}$, that each leg receives. 
At high speeds, when $I_{ext}$ is large, the model
generates a unique stable tripod gait, as observed experimentally
in cockroaches and fruit flies. However, for small $I_{ext}$, the model's low speed
dynamics include both stable forward and backward tetrapod gaits and a stable
gait that has not been observed in insects, in which triple, double and single
swing phases occur \cite[Figure 29]{SIADS2018}. 
While fruit flies exhibit forward and backward tetrapod gaits at low speeds,
the latter have only been seen in backward walking \cite{Bidaye-Science14},
and we therefore propose that brain or central nervous system inputs are
likely used to switch among and select particular gaits.

In the present paper, we relax the assumption of homogeneous oscillators
and allow heterogeneous ipsilateral external inputs denoted by $I_{ext}+\Ii$ 
for $i=1,2,3$. We observe that, at low speed with small $I_{ext}$, and
for appropriate choices of small heterogeneities $\Ii$, $i=1,2,3$, 
the heterogeneous model generates only one stable approximate forward
or backward tetrapod gait, as is expected experimentally. The selection
of a stable gait is accomplished via sequences of saddle-node bifurcations
in which all but one of the stable gaits disappears as particular currents $\Ii$ increase.
See Table~\ref{summary:figs} for a summary of the behaviors presented in 
Section~\ref{Main_Result}.

At high speeds the single stable solution of the heterogeneous model
is a tripod gait, as in the homogeneous case, and the model exhibits
a transition from a forward or a backward tetrapod to a tripod gait as $I_{ext}$
increases (see Figure~\ref{NC_hetero_random_transition} for the former case).

In future work, we propose to allow the heterogeneous external inputs to be \textit{noisy} 
and to study the resulting effects on the existence of gaits and their transitions.

\section*{Acknowledgements}

This work was supported by the National Science Foundation under NSF-CRCNS grant DMS-1430077. 


\small{
}


\begin{thebibliography}{10}

\bibitem{SIADS2018}
Z.~Aminzare, V.~Srivastava, and P.~Holmes.
\newblock Gait transitions in a phase oscillator model of an insect central
  pattern generator.
\newblock {\em SIAM J. Appl. Dyn. Sys.}, 17(1):626--671, 2018.

\bibitem{SIAM2}
R.M.~Ghigliazza and P.~Holmes.
\newblock A minimal model of a central pattern generator and motoneurons for
  insect locomotion.
\newblock {\em SIAM J. Appl. Dyn. Sys.}, 3(4):671--700, 2004.

\bibitem{SIAM1}
R.M.~Ghigliazza and P.~Holmes.
\newblock Minimal models of bursting neurons: How multiple currents,
  conductances, and timescales affect bifurcation diagrams.
\newblock {\em SIAM J. Appl. Dyn. Sys.}, 3(4):636--670, 2004.

\bibitem{pearson_Iles_1973}
K.G.~Pearson and J.F.~Iles.
\newblock Nervous mechanisms underlying intersegmental co-ordination of leg
  movements during walking in the cockroach.
\newblock {\em J. Exp. Biol.}, 58:725--744, 1973.

\bibitem{pearson_Iles_1970}
K.G.~Pearson and J.F.~Iles.
\newblock Discharge patterns of coxal levator and depressor motoneurons of the
  cockroach, \emph{Periplaneta americana}.
\newblock {\em J. Exp. Biol.}, 52:139--165, 1970.

\bibitem{pearson_1972}
K.G.~Pearson.
\newblock Central programming and reflex control of walking in the cockroach.
\newblock {\em J. Exp. Biol.}, 56:173--193, 1972.

\bibitem{Fuchs14}
E.~Couzin-Fuchs, T.~Kiemel, O.~Gal, A.~Ayali, and P.~Holmes.
\newblock Intersegmental coupling and recovery from perturbations in freely
  running cockroaches.
\newblock {\em J. Exp. Biol.}, 218(2):285--297, 2015.

\bibitem{KukProc09}
R.P.~Kukillaya, J.L.~Proctor, and P.~Holmes.
\newblock Neuromechanical models for insect locomotion: {S}tability,
  maneuverability, and proprioceptive feedback.
\newblock {\em Chaos}, 19(2):026107, 2009.

\bibitem{Schwemmer2012}
M.A.~Schwemmer and T.J.~Lewis.
\newblock The theory of weakly coupled oscillators.
\newblock In N.W.~Schultheiss, A.A.~Prinz, and R.J.~Butera, editors, {\em Phase
  Response Curves in Neuroscience: Theory, Experiment, and Analysis}, pages
  3--31. Springer New York, New York, NY, 2012.

\bibitem{Kuramoto_book}
Y.~Kuramoto.
\newblock {\em Chemical Oscillations, Waves, and Turbulence}, 
volume~19 of {\em Springer Series in Synergetics}.
\newblock Springer-Verlag, Berlin, 1984.
  
 \bibitem{ginverse1971}
C.R.~Rao and S.~Kumar Mitra.
\newblock {\em {G}eneralized {I}nverse of {M}atrices and its {A}pplications}.
\newblock John Wiley \& Sons, Inc., New York-London-Sydney, 1971.

\bibitem{Yeldesbay2018}
A.~Yeldesbay, T.~T{\'o}th, and S.~Daun.
\newblock The role of phase shifts of sensory inputs in walking revealed by
  means of phase reduction.
\newblock {\em J.  Comput. Neurosci.}, 44(3):313--339, 
  2018.

\bibitem{silvester2000}
J.R.~Silvester.
\newblock {Determinants of Block Matrices}.
\newblock {\em {Mathematical Gazette}}, 84(501):460--467, 2000.

  \bibitem{kukillaya.09}
R.P.~Kukillaya and P.~Holmes.
\newblock A model for insect locomotion in the horizontal plane: {F}eedforward
  activation of fast muscles, stability, and robustness.
\newblock {\em J. Theor. Biol.}, 261 (2):210--226, 2009.

\bibitem{ProctorKH10}
J.~Proctor, R.P.~Kukillaya, and P.~Holmes.
\newblock A phase-reduced neuro-mechanical model for insect locomotion:
  feed-forward stability and proprioceptive feedback.
\newblock {\em Phil. Trans. Roy. Soc. A}, 368:5087--5104, 2010.

\bibitem{Bidaye-Science14}
S.S.~Bidaye, C.~Machacek, Y.~Wu, and B.J.~Dickson.
\newblock Neuronal control of {\it drosophila} walking direction.
\newblock {\em Science}, 344:97--101, 2014.

\end{thebibliography}
\end{document}